\documentclass[final]{siamart190516}

\usepackage[utf8]{inputenc}
\usepackage{amsmath}
\usepackage{hyperref}
\usepackage{cleveref}
\RequirePackage{amsfonts}
\RequirePackage{amsmath}
\RequirePackage{graphicx}

\usepackage{url}
\usepackage[space,noadjust,nocompress]{cite}

\newcommand{\eps} {\varepsilon}                            
\DeclareMathOperator*{\argmin}{arg\,min}                   
\renewcommand{\t} {^{\top}}                                
\newcommand{\norm} [2][]{\left\|#2\right\|_{#1}}           
\newcommand{\diag} [1]  {{\rm diag\!}\left( #1 \right)}    
\newcommand{\trace} [1]  {{\rm trace\!}\left( #1 \right)}    

\newcommand{\bfGamma}{{\boldsymbol{\Gamma}}}

\newcommand{\bfLambda}{{\boldsymbol{\Lambda}}}

\newcommand{\bfSigma}{{\boldsymbol{\Sigma}}}

\newcommand{\bfPsi}{{\boldsymbol{\Psi}}}

\newcommand{\bfmu}{{\boldsymbol{\mu}}}

\newcommand{\bfA}{{\bf A}}
\newcommand{\bfB}{{\bf B}}
\newcommand{\bfC}{{\bf C}}
\newcommand{\bfD}{{\bf D}}

\newcommand{\bfG}{{\bf G}}
\newcommand{\bfH}{{\bf H}}
\newcommand{\bfI}{{\bf I}}
\newcommand{\bfJ}{{\bf J}}

\newcommand{\bfL}{{\bf L}}
\newcommand{\bfM}{{\bf M}}

\newcommand{\bfP}{{\bf P}}
\newcommand{\bfQ}{{\bf Q}}
\newcommand{\bfR}{{\bf R}}
\newcommand{\bfS}{{\bf S}}
\newcommand{\bfT}{{\bf T}}
\newcommand{\bfU}{{\bf U}}
\newcommand{\bfV}{{\bf V}}
\newcommand{\bfW}{{\bf W}}

\newcommand{\bfZ}{{\bf Z}}

\newcommand{\bfb}{{\bf b}}
\newcommand{\bfc}{{\bf c}}
\newcommand{\bfd}{{\bf d}}
\newcommand{\bfe}{{\bf e}}

\newcommand{\bfg}{{\bf g}}

\newcommand{\bfp}{{\bf p}}
\newcommand{\bfq}{{\bf q}}
\newcommand{\bfr}{{\bf r}}
\newcommand{\bfs}{{\bf s}}
\newcommand{\bft}{{\bf t}}
\newcommand{\bfu}{{\bf u}}
\newcommand{\bfv}{{\bf v}}
\newcommand{\bfw}{{\bf w}}
\newcommand{\bfx}{{\bf x}}
\newcommand{\bfy}{{\bf y}}
\newcommand{\bfz}{{\bf z}}
\newcommand{\bfzero}{{\bf0}}

\newcommand{\calC}{\mathcal{C}}

\newcommand{\calE}{\mathcal{E}}

\newcommand{\calG}{\mathcal{G}}
\newcommand{\calH}{\mathcal{H}}

\newcommand{\calJ}{\mathcal{J}}
\newcommand{\calK}{\mathcal{K}}

\newcommand{\calN}{\mathcal{N}}

\newcommand{\calQ}{\mathcal{Q}}
\newcommand{\calR}{\mathcal{R}}
\newcommand{\calS}{\mathcal{S}}

\newcommand{\calV}{\mathcal{V}}
\newcommand{\calW}{\mathcal{W}}

\newcommand{\bbE}{\mathbb{E}}

\newcommand{\bbR}{\mathbb{R}}

\newcommand{\true}{\mathrm{true}}

\usepackage{tikz}
\usetikzlibrary{shapes,arrows}
\usetikzlibrary{decorations.pathreplacing}
\tikzstyle{decision} = [diamond, draw, fill=blue!20,
    text width=4.5em, text badly centered, node distance=3cm, inner sep=0pt]
\tikzstyle{blocky} = [rectangle, draw, fill=yellow!20,
    text width=13em, text centered, rounded corners, minimum height=4em]
\tikzstyle{blockb} = [rectangle, draw, fill=blue!20,
        text width=13em, text centered, rounded corners, minimum height=4em]
\tikzstyle{line} = [draw, -latex']
\tikzstyle{cloud} = [draw, ellipse,fill=red!20, node distance=3cm,
    minimum height=2em]

    \definecolor{matlab1}{RGB}{0,    114,  189} 
    \definecolor{matlab2}{RGB}{217,   83,   25}
    \definecolor{matlab3}{RGB}{237,  177,   32}
    \definecolor{matlab4}{RGB}{126,   47,  142}
    \definecolor{matlab5}{RGB}{119,  172,   48}
    \definecolor{matlab6}{RGB}{77,   190,  238}
    \definecolor{matlab7}{RGB}{162,   20,   47}

\usepackage{cprotect}

\newdimen\iwidth
\newdimen\iheight

\usepackage{bigstrut}
\usepackage{multirow}

\usepackage{amsmath,amssymb} 
\usepackage{algorithm,algorithmic}

\usepackage{color}

\newcommand{\bbfb}{\bar{\bfb}}
\newcommand{\IRtools}{\textsc{IR Tools}}

\newcommand{\bfAregpinv}{\bfA_{\text{{reg}}}^\dagger}
\newcommand{\bfxreg}{\bfx_{\text{{reg}}}}
\newcommand{\bfrreg}{\bfr_{\text{{reg}}}}

\newcommand{\TV}{\mbox{TV}}
\newcommand{\bfUA}{\bfU^\text{\tiny A}}
\newcommand{\bfSA}{ \bfSigma^\text{\tiny A}}
\newcommand{\bfLA}{ \bfLambda^\text{\tiny A}}

\newcommand{\bfUB}{\bfU^\text{\tiny B}}
\newcommand{\bfSB}{ \bfSigma^\text{\tiny B}}
\newcommand{\bfVB}{\bfV^\text{\tiny B}}

\numberwithin{algorithm}{section}

\usepackage{listings}

\newcommand{\TheTitle}{Computational methods for large-scale inverse problems: a survey on hybrid projection methods}

\title{{\TheTitle}\thanks{Current version: \today.}}
\author{
  Julianne Chung\thanks{Department of Mathematics, Computational Modeling and Data Analytics Division, Academy of Integrated Science, Virginia Tech, Blacksburg, VA, USA
    (jmchung@vt.edu, \url{http://www.math.vt.edu/people/jmchung/}).}
  \and
  Silvia Gazzola\thanks{Department of Mathematical Sciences, University of Bath, United Kingdom
    (S.Gazzola@bath.ac.uk, \url{http://people.bath.ac.uk/sg968/}).}
}

\begin{document}
\maketitle

\begin{abstract} This paper surveys an important class of methods that combine iterative projection methods and variational regularization methods for large-scale inverse problems. Iterative methods such as Krylov subspace methods are invaluable in the numerical linear algebra community and have proved important in solving inverse problems due to their inherent regularizing properties and their ability to handle large-scale problems. Variational regularization describes a broad and important class of methods that are used to obtain reliable solutions to inverse problems, whereby one solves a modified problem that incorporates prior knowledge. \emph{Hybrid projection methods} combine iterative projection methods with variational regularization techniques in a synergistic way, providing researchers with a powerful computational framework for solving very large inverse problems. Although the idea of a hybrid Krylov method for linear inverse problems goes back to the 1980s, several recent advances on new regularization frameworks and methodologies have made this field ripe for extensions, further analyses, and new applications. In this paper, we provide a practical and accessible introduction to hybrid projection methods in the context of solving large (linear) inverse problems.
\end{abstract}


\section{Introduction}
\label{sec:introduction}
We provide a gentle introduction to hybrid projection methods for regularizing inverse problems by answering three essential questions: (1) what is an inverse problem?, (2) what is regularization and why do we need it?, and (3) why should we use hybrid projection methods?

\subsection{What is an inverse problem?}
\label{subsec:IP}
Inverse problems arise in various scientific applications, including astronomy, geoscience, biomedical sciences, mining engineering, and medicine (see,  e.g., \cite{AnHu, BoReHaMo96, HaNaOL06, Gon82, JeScMa02, Eps08, ChHaNa06, Ha14, hansen2018air, BoBo98, Wa15, Bu08}). In these and other applications, one formulates an inverse problem for the purpose of extracting important information from available, noisy measurements. In this survey, we are mainly interested in large, linear inverse problems of the form
\begin{equation}
  \label{eq:linearmodel}
  \bfb = \bfA \bfx_\true + \bfe\,,
\end{equation}
where a certain unknown quantity of interest is stored in 
$\bfx_\true \in \bbR^n,$ the observed measurements are contained in $\bfb \in \bbR^m$, the matrix $\bfA \in \bbR^{m \times n}$ represents the forward data acquisition process, and $\bfe \in \bbR^m$ represents inevitable errors (noise) that arise from measurement, discretization, or floating point arithmetic.
Given measured data, $\bfb$, and knowledge of the forward model, $\bfA,$ the goal is to compute an approximation of $\bfx_\true$, i.e., a solution to the inverse problem.

It is worth mentioning that, in stating \cref{eq:linearmodel}, we have made three simplifying assumptions. First, we assume that the matrix $\bfA$ is known exactly while, in practice, one should take into account the error between the assumed mathematical model and the underlying physical model.
Second, many inverse problems have an underlying mathematical model that is not linear, i.e., $\bfb = F(\bfx_\true) + \bfe$, where $F(\cdot):\bbR^n \to \bbR^m$ is a nonlinear operator that is sometimes referred to as the `parameter-to-observation' map.  For these cases, more sophisticated methods should be used to solve the nonlinear inverse problem; however, many of these methods require solving a subproblem with an approximate linear model, and the methods described herein have been successfully used for this purpose. Third, we have assumed an additive noise model. Other realistic assumptions (e.g., multiplicative noise or mixed noise models) may need to be incorporated in the model. These simplifying assumptions will be dropped in \cref{sec:extensions}.

\subsubsection*{Two model applications}
In this paper, we focus on two model inverse problems from image processing: image deconvolution and tomographic reconstruction. These two problems will be used throughout the paper for illustrative purposes, so we briefly describe them here. All of the test problems presented in this survey can be generated using the MATLAB toolbox \IRtools\ \cite{IRtools}\footnote{The MATLAB programs used to produce the illustrations and experiments reported herein are available at the website: {\url{https://github.com/juliannechung/surveyhybridprojection}}. The commented lines in the MATLAB files are designed to make the various commands comprehensible to readers with basic programming experience.}; see also \cref{sec:software}.\smallskip

\emph{Image deconvolution (or deblurring)} problems are very popular in the image and signal processing literature, and are of core importance in fields such as astronomy, biology, and medicine. More precisely, the mathematical model of this problem can be expressed in the continuous setting as an integral equation
\begin{equation}\label{eq:deconv}
b(\bfs) = \int a(\bfs,\bft)x(\bft)d\bfs \quad+ e\,,
\end{equation}
where $\bfs,\bft\in\bbR^2$ represent spatial locations.
The kernel or point spread function (PSF) $a(\bfs,\bft)$ defines the blur, and, if the kernel has the property that $a(\bfs, \bft) = a(\bfs-\bft)$, then the blur is spatially invariant and the integration in \cref{eq:deconv} is a convolution.
From equation \cref{eq:deconv} it can be observed that each pixel in the blurred image is formed by integrating the PSF with pixel values of the true image scene, and then further corrupted by adding a random perturbation $e$.

In a realistic setting, images are collected only at discrete points (pixels), and are only available in a finite region (i.e., in a viewable region).
Thus, the basic image deconvolution problem is of the form \cref{eq:linearmodel} where $\bfx_\true$ represents the vectorized sharp image, $\bfA$ represents the blurring process
which specifies how the points in the image are distorted, and $\bfb$ contains the observed, vectorized, blurred and noisy image. Here we assume that the true and corrupted images have the same size, so that $\bfA\in \bbR^{n \times n}$.
Generally the integration operation is local, and so pixels in the center of the viewable region are well defined. This results in a sparse, structured matrix $\bfA$. However, pixels of the blurred image near the boundary of the viewable region are affected by information outside the viewable region. Therefore, in constructing the matrix $\bfA$, one needs to incorporate boundary conditions to model how the image scene extends beyond the boundaries of the viewable region. Typical boundary conditions include zero, periodic, and reflective. We highlight that the discrete problem associated to \cref{eq:deconv} often gives rise to matrices $\bfA$ with a well-defined structure: for instance, if the blur is spatially invariant and periodic boundary conditions are assumed, then $\bfA$ is a block circulant matrix with circulant blocks and efficient implementations of any image deconvolution algorithm can be obtained by exploiting structure of the matrix $\bfA$. We refer to \cite{HaNaOL06} for a discussion of the many fundamental modeling aspects of the image deconvolution problem and relevant features of the discrete problem.

\begin{figure}
\begin{center}
\begin{tabular}{ccc}
  \includegraphics[height=3cm]{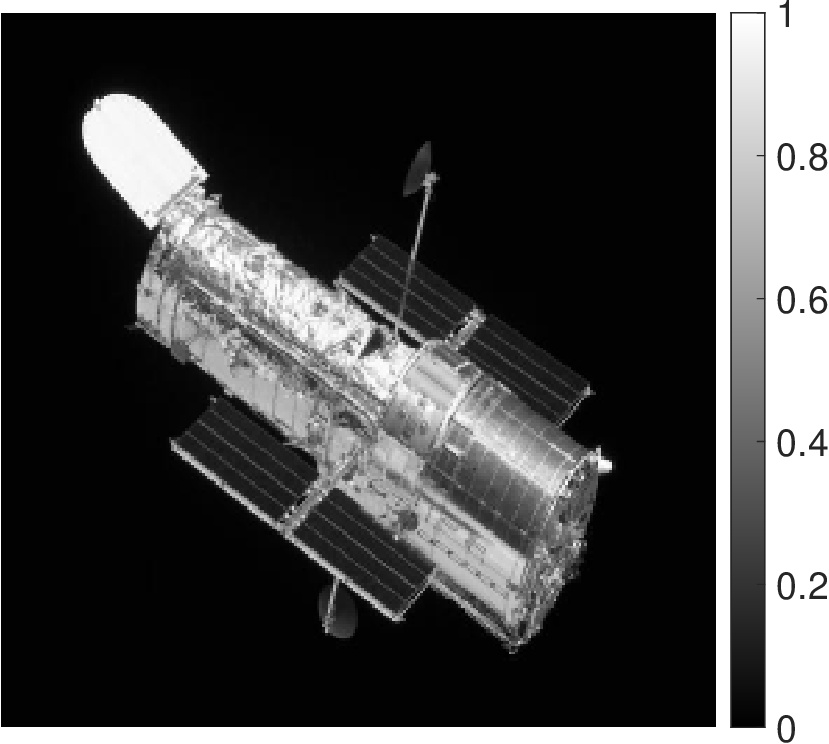} &
  \includegraphics[height=3cm]{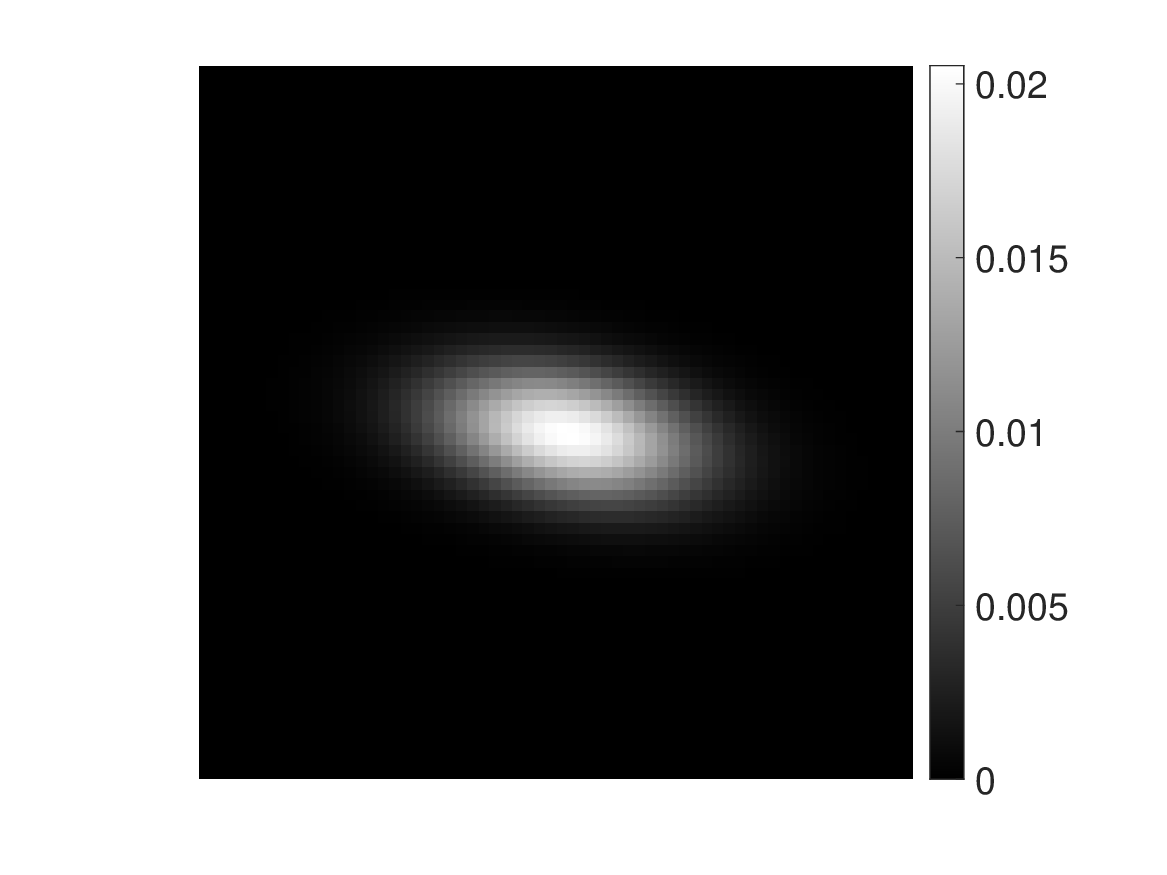} &
  \includegraphics[height=3cm]{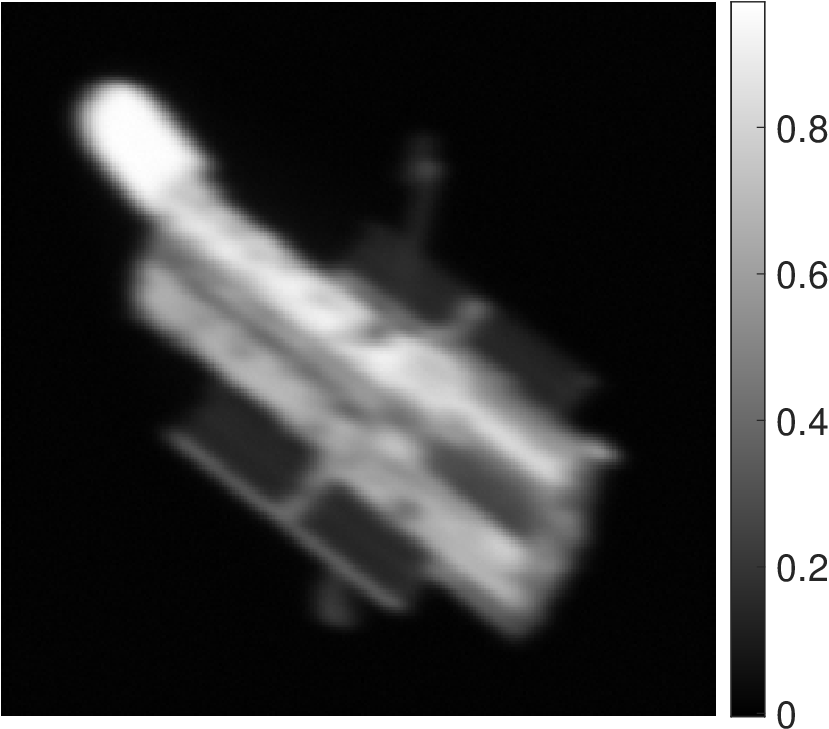}\\
\hspace{-1.2cm}\small{\textbf{(a)} true image} & \hspace{-1cm}\small{\textbf{(b)} PSF} & \hspace{-1cm}\small{\textbf{(c)} observed image}
\end{tabular}
\end{center}
\caption{\label{fig:deblur_testpb} Image deblurring test problem. \textbf{(a)} Sharp image of $256\times 256$ pixels. \textbf{(b)} Zoomed image (400\%) of an anisotropic Gaussian PSF. \textbf{(c)} Blurred and noisy image.}
\end{figure}

In \cref{fig:deblur_testpb} we display the test data that will be used through this paper, and that can be generated running the MATLAB script \texttt{generate\_blur}.
The size of the images are $256\times 256$ pixels, so $n=256^2=65,536$. The matrix $\bfA$ models a spatially invariant blur, with
an anisotropic Gaussian PSF (whose analytical expression is given in \cite[Chapter 3]{HaNaOL06}) and reflective (sometimes referred to as reflexive) boundary conditions.
\smallskip

\emph{Tomographic reconstruction (or computed tomography, CT)} is a critical tool in many applications such as nondestructive evaluation {and} electron microscopy; the advent of newer technologies (e.g., spectral CT and photoacoustic tomography) prompts more efficient and more accurate reconstruction methods. CT consists in computing reconstructions of (parameters of) objects from projections, i.e., data obtained by integrations along rays (typically straight lines) that penetrate a domain.
More precisely, in a (semi)continuous setting, the unknown parameters are modeled as a function $x(\bft)$, $\bft\in\bbR^2$,
and it is assumed that the damping of a infinitesimally small portion of length $ds$ of a ray penetrating the object at position $\bft$ is proportional to the product $x(\bft)ds$.
The $i$th observation $b_i$, $i=1,\dots,m$, is the damping of a signal that penetrates the object along the $i$th ray $\bft_i$ (parameterized by a ray length $s\in\bbR$), and can be represented by the integral
\begin{equation}\label{eq:CT}
b_i = \int x(\bft_i(s))ds \quad + e_i,
\end{equation}
where $e_i$ is a random perturbation corrupting the measurement.
Employing a classical discretization scheme, which subdivides the object into an array of pixels and assumes that the function $x(\bft)$ is constant in each pixel, the above integral can be expressed as a discrete sum and the following expression for the $(i,j)$th entry $a_{ij}$ of the sparse matrix $\bfA$ can be readily derived as
\[
a_{ij}=\begin{cases}
L_{ij} & \text{if $j\in\calS_i$}\\
0 & \text{otherwise}
\end{cases},
\]
where $\calS_i$ is the set of indices of those pixels that are penetrated by the $i$th ray, $L_{ij}$ is the length of the $i$th ray through the $j$th pixel.
Thus, the basic CT problem \cref{eq:CT} is of the form \cref{eq:linearmodel}, where $\bfx_\true$ represents the unknown material parameters, $\bfA$ represents a physical attenuation (damping) process, and $\bfb$ contains the vectorized projections (as {measured} by a detector) and is often referred to as a sinogram. In most cases of tomography, $\bfA$ is a rectangular matrix, where $m < n$ if fewer measurements are collected than the number of unknowns, or $m > n$ if many projections can be obtained, or parameterization reduces the number of unknowns.
We refer to \cite{Bu08} for more details on the tomography reconstruction problem.

\begin{figure}
\begin{center}
\begin{tabular}{ccc}
\includegraphics[width=0.25\textwidth]{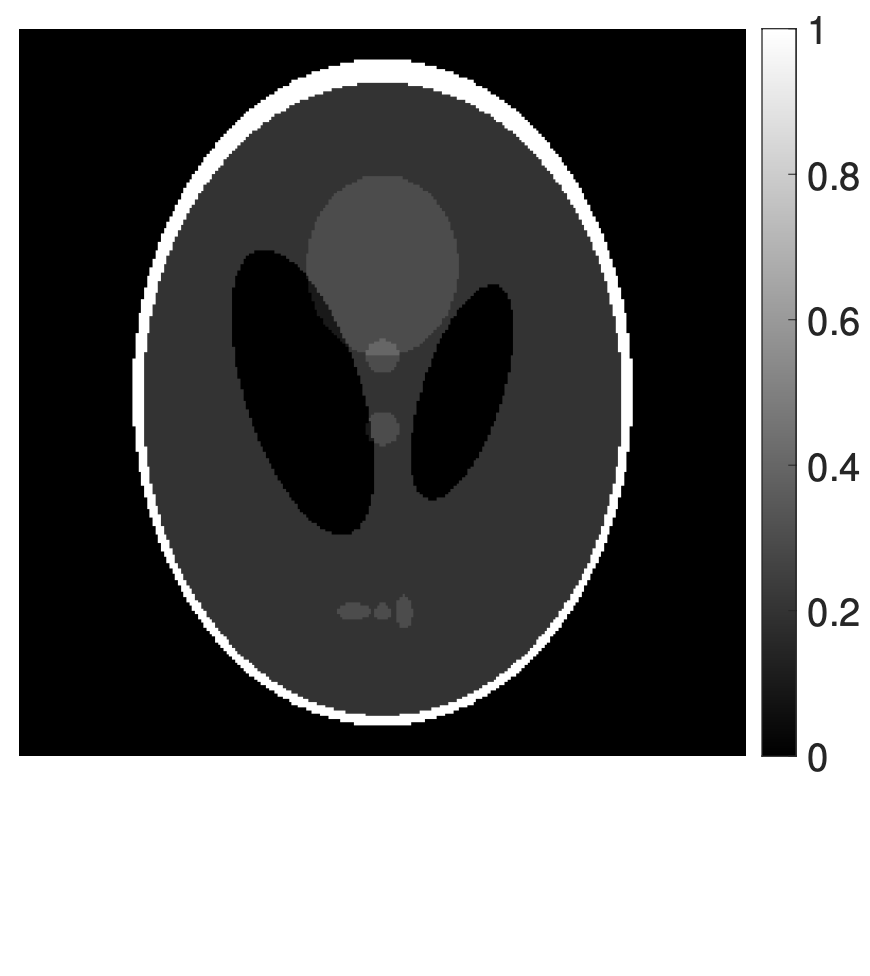} &
\resizebox{.35\linewidth}{!}{\input{figures/tomographyFig.tex}} &
\includegraphics[width=0.2\textwidth]{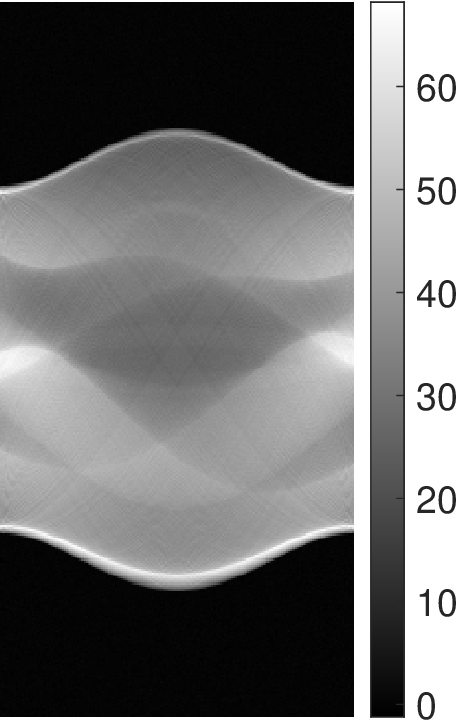} \\
\small{\textbf{(a)} true image} &  \small{\textbf{(b)} tomography setup} & \small{\textbf{(c)} observations}
\end{tabular}
\end{center}
\caption{\label{fig:tomo_testpb} Parallel ray CT test problem. \textbf{(a)} Phantom of size $256\times 256$ pixels.  \textbf{(b)} Illustration of 2D x-ray parallel-beam tomography setup, modified from \cite{ruthotto2018optimal}. \textbf{(c)} Observed sinogram containing data.}
\end{figure}

In \cref{fig:tomo_testpb} we display the test data for tomography that will be used through this paper, and that can be generated running the MATLAB script \texttt{generate\_tomo}: this program still draws on \IRtools, which itself calls functions available within the MATLAB package \textsc{AIR Tools II} \cite{hansen2018air}. The object we wish to recover is the Shepp-Logan `medical' phantom of size $256\times 256$ pixels, so that $n=65,536$. The damping matrix is modeled after a parallel tomography process, where $p$ parallel rays are generated from a source ideally placed infinitely far from a flat detector; moreover, the source-detector pair is rotated around the object, and measurements are recorded for $N_\beta$ angles. Hence the number of observations is $m=pN_\beta$. One can vary the parameters that define the measurement geometry, such as the number of angles or the number of rays.  For this example, we take measurements from $0$ to $179$ degrees in intervals of $1$ degree, resulting in a set of data with $m=64798$ measurements, so that $m>n$.

Now that we have addressed what is an inverse problem and looked at some examples of inverse problems, our focus turns to computing solutions to inverse problems, and for this we need regularization.

\subsection{What is regularization and why do we need it?}
\begin{figure}
\begin{center}
\begin{tabular}{cc|cc}
    \multicolumn{2}{c|}{deblurring example} &   \multicolumn{2}{c}{tomography example} \\
\includegraphics[width=0.2\textwidth]{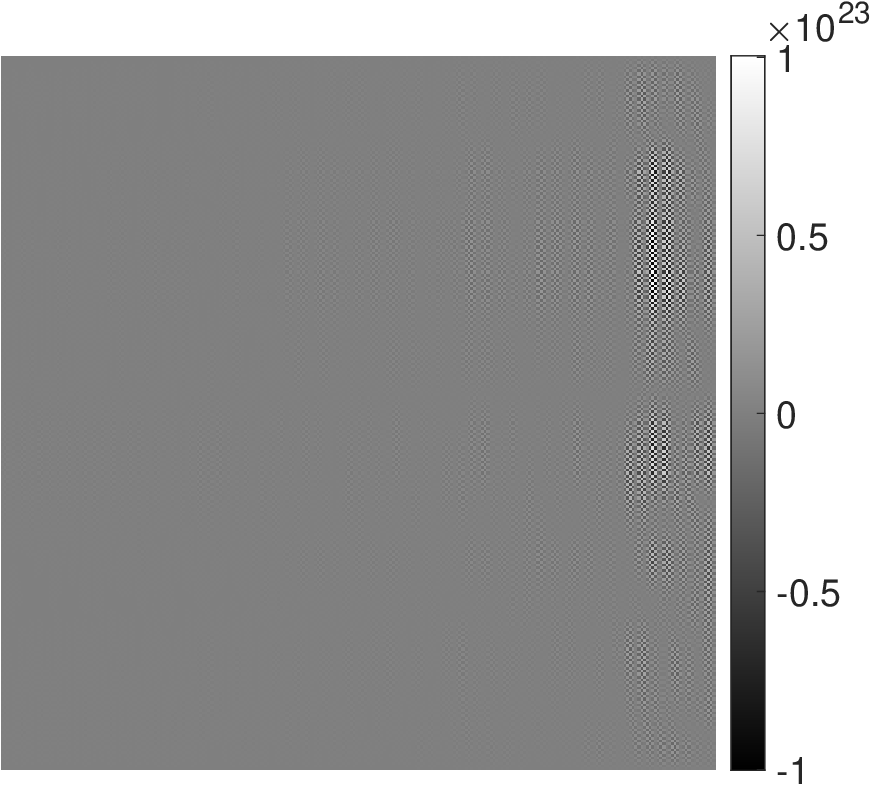} &
\includegraphics[width=0.2\textwidth]{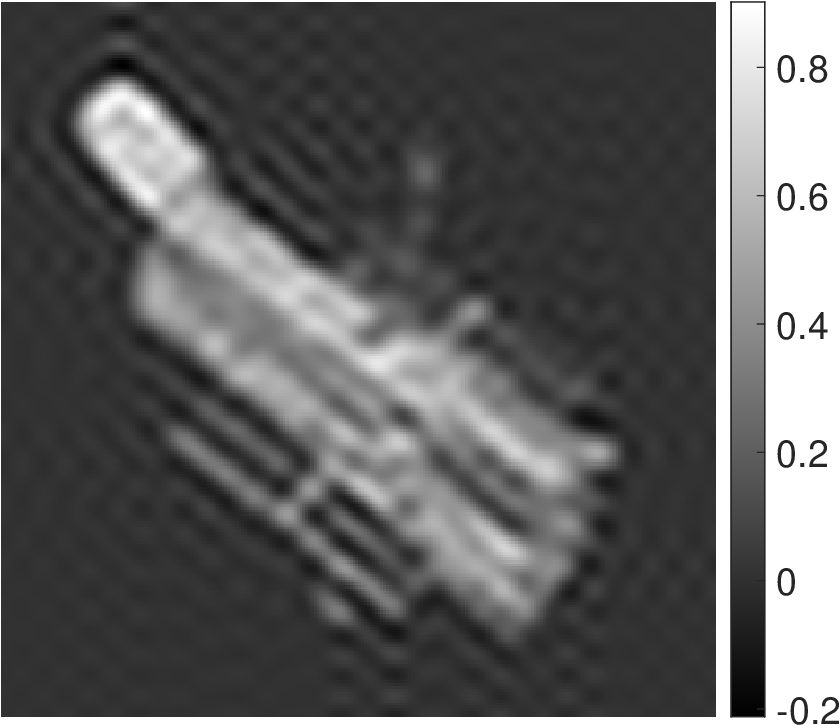} &
\includegraphics[width=0.2\textwidth]{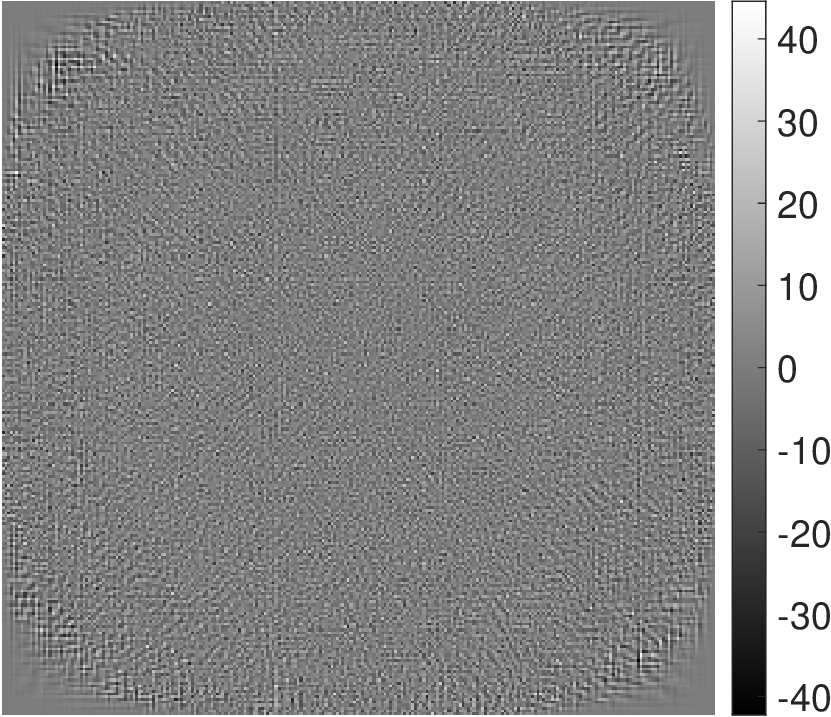} &
\includegraphics[width=0.2\textwidth]{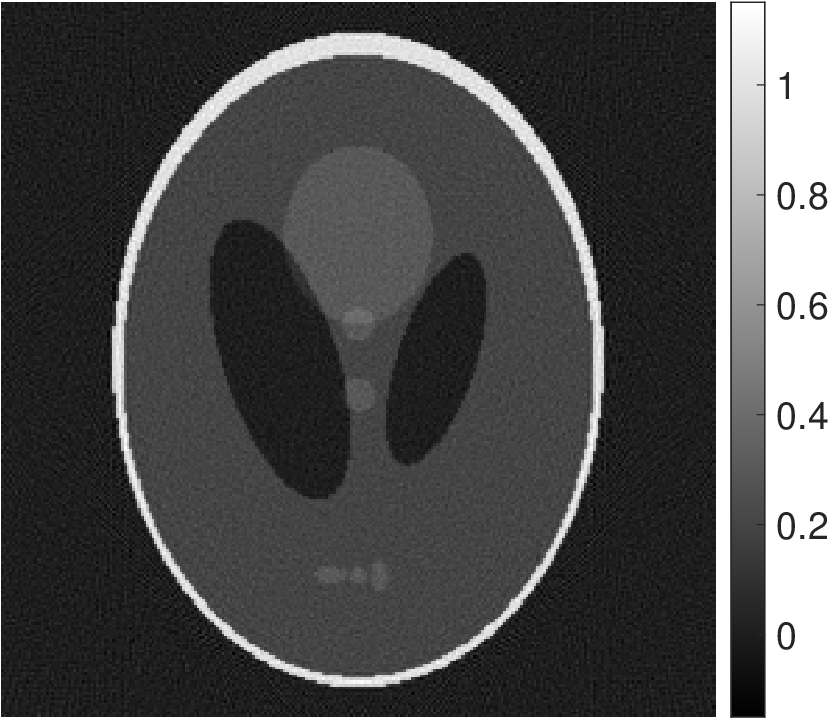}\\
\small{inverse solution} &
\small{regularized solution} &
\small{inverse solution} &
\small{regularized solution}
\end{tabular}
\end{center}
\label{fig:naive}
\caption{For each of the image deblurring and tomography test problems, we provide the computed inverse solution $\bfA^{-1}\bfb$ or $\bfA^{\dagger}\bfb$ on the left and a regularized solution on the right. Clearly the inverse solution is heavily corrupted with noise and errors, while the regularized solution can provide a good approximation to the true object.}
\end{figure}

Solving inverse problems is notoriously difficult due to \emph{ill-posedness}\footnote{Here we consider the discretized problem, but point the reader to \cite{engl1996regularization,hansen2010discrete} for discussions on the full implications of ill-posedness in continuous formulations.}, whereby data acquisition noise and computational errors can lead to large changes in the computed solution.
For instance, if $\bfA$ is invertible, the solution $\bfA^{-1}\bfb$ is dominated by noise and is useless for all practical purposes. The same happens for a more general $\bfA$ and for the solution expressed as
\begin{equation}\label{xpseudoinv}
\bfx = \bfA^\dagger\bfb,
\end{equation}
where $\bfA^\dagger$ is the generalized inverse of $\bfA$; see \cite{golub2012matrix}.
We illustrate this phenomenon for the described test problems, where a tiny amount of noise (here, $\|\bfe\|_2/\|\bfA \bfx_{\rm true}\|_2=10^{-4}$) enters the data collection process. Naive inverse solutions for the image deblurring and tomography test problems are displayed in the first and third frames of \cref{fig:naive}, respectively. It is clear that these solutions are unacceptable approximations to $\bfx_\true$ (so much that, for image deblurring, the available corrupted image $\bfb$ even looks better!).

\emph{Regularization} can be used to overcome the inherent instability of ill-posed problems in order to compute a meaningful solution.  The basic idea of regularization is to replace problem \cref{eq:linearmodel} by a modified problem that incorporates other information about the solution.  If proper regularization is applied, a regularized solution (e.g., the ones shown in the second and fourth frame of \cref{fig:naive}) should be close to the true solution. Regularization techniques come in many forms, and more details will be unfolded in the coming sections. However, for now, we consider two types of methods: variational regularization methods and iterative regularization methods.

\emph{Variational regularization methods} for \cref{eq:linearmodel}
involve solving optimization problems of the form,
\begin{equation}
	\label{eq:directreg}
	\min_{\bfx \in \calC} \calJ(\bfb - \bfA \bfx) + \lambda \calR(\bfx)\,,
	\end{equation}
where $\calJ$ is a loss (or fit-to-data) function,
$\calR$ is a regularization operator\footnote{Commonly-used regularization operators are based on smoothness constraints or assumptions of Gaussian priors.},
$\lambda>0$ is a regularization parameter that controls the amount of regularization, thereby determining how faithful the modified problem is to the original problem, and $\calC$ denotes the set of feasible solutions (e.g., those that satisfy some constraints). The main advantage of formulation \cref{eq:directreg} is that general constraints and priors can be easily incorporated in the problem. Furthermore, well-known optimization methods can be used to solve \cref{eq:directreg}.  For the particular case where $\calJ(\cdot)=\calR(\cdot)=\|\cdot\|_2^2$, direct factorization methods can, in principle, be used to compute a solution (see \cref{sub:direct}), although this is often computationally infeasible for large-scale problems. In this setting, the main disadvantage is the need to select the regularization parameter $\lambda$, often prior to solution computation.

On the other hand, \emph{iterative regularization methods} for \cref{eq:linearmodel} typically consist in applying an iterative solver to
\begin{equation}\label{eq:noreg}
	\min_\bfx \calJ(\bfb - \bfA \bfx),
\end{equation}
and computing a regularized solution by early termination of the iterations. In practice, when $\calJ$ is expressed in the 2-norm, the iterative methods of choice are often subspace projection methods, whereby the problem is projected onto increasing (but {relatively} small) subspaces, and a projected subproblem is solved at each iteration by imposing certain optimality conditions on the approximated solution. Although an explicit choice of the regularization parameter $\lambda$ is no longer required (contrary to \cref{eq:directreg}), the stopping iteration essentially serves as a regularization parameter, as it balances how `faithful' the projected problem is to the original problem. In this setting, the main disadvantage is that general constraints cannot be easily incorporated, and, even when they can, they can only be handled through complicated nested iterative schemes.

\emph{Hybrid regularization methods} combine variational and iterative regularization methods, leveraging the best features of each class, to provide a powerful computational framework for solving large-scale inverse problems. In this paper, we focus specifically on \emph{hybrid projection methods} that start off as iterative regularization methods, i.e., the original problem \cref{eq:linearmodel} is projected onto a subspace of increasing dimension, and then the projected subproblem is solved using a variational regularization method.

\subsection{Why should we use hybrid projection methods?}
The main motivations for using hybrid projection methods to solve large inverse problems can be grouped as follows.
\begin{itemize}
  \item For many multidimensional inverse problems (such as the model ones described in \cref{subsec:IP}), the matrix $\bfA$ is very large (to the point that it cannot be explicitly stored). In this setting, only matrix-vector multiplications with $\bfA$ (and possibly $\bfA\t$) can be performed, most often taking advantage of high-performance computing or tools such as GPUs. Therefore, methods that work explicitly with $\bfA$ (e.g., methods that require some factorizations of $\bfA$, such as the singular value decomposition) are ruled out.  Hybrid regularization methods can be implemented without explicitly constructing and storing the matrix $\bfA$, by treating $\bfA$ and $\bfA\t$ as linear operators.
	\item Even when adopting a variational regularization method \cref{eq:directreg},
  one may not know a good regularization parameter $\lambda$ in advance, and standard methods to estimate $\lambda$ can be expensive: indeed, these typically require first approximating problem \cref{eq:linearmodel}, or solving many instances of \cref{eq:linearmodel} with different regularization parameters.
 Hybrid projection methods provide a natural setting to estimate an appropriate value of $\lambda$ during the solution computation process: often in a heuristic way, but in some cases with theoretical guarantees.
 \item As mentioned above, many iterative solvers for \cref{eq:noreg}, such as projection methods onto Krylov subspaces, have inherent regularizing properties and avoid data overfitting by early termination of the iterations. By using hybrid methods,
it is possible to stabilize and enhance the regularized solutions computed by these methods. Moreover, hybrid methods enjoy nice theoretical properties.
For instance we know that, for many hybrid frameworks, iterates obtained by first projecting the problem and then regularizing are mathematically equivalent to iterates obtained by first regularizing and then projecting the regularized problem.
\item Hybrid methods provide a convenient framework for the development of extensions to more general problems and regularization terms.
For instance, one could incorporate a variety of regularization terms of the form $\calR(\bfx)=\|\bfL\bfx\|_2^2$, where $\bfL\in\bbR^{p\times n}$, or also consider functions $\calR(\bfx)$ that are not expressed in the 2-norm.
As we will see in \cref{sec:extensions}, this amounts to modifying
the projection subspace that can be used within hybrid methods. By doing so, other popular (but expensive) iterative regularization schemes (which are often based on nested cycles of iterations) can be avoided.
\end{itemize}

\subsection{Outline of the paper}
\label{subsec:outline}
The remaining part of this survey is organized as follows. In \cref{sec:background}, we describe the general problem setup to provide some background for hybrid methods, including an overview of the most relevant direct and iterative regularization techniques. Then, in \cref{sec:hybrid}, we discuss hybrid projection methods.  After providing a brief summary of the historical developments of hybrid methods, we address the two main building blocks of any hybrid method: namely, generating a subspace for the solution (\cref{sub:buildsubspace}) and solving the projected problem (\cref{sub:solvingprojected}). Important numerical and theoretical aspects will also be covered: these include strategies to efficiently set the Tikhonov regularization parameter and stop the iterations (\cref{sec:param}) and available convergence proofs and approximation properties (\cref{sec:theory}). In recent years, we have witnessed the extension of hybrid methods to solve a larger scope of problems and to cover broader scientific applications. In \cref{sec:extensions}, we provide an overview of some of these extensions:
\begin{itemize}
 \item Beyond standard-form Tikhonov: Hybrid projection methods for general-form Tikhonov; see \cref{sub:generalform}
  \item Beyond standard projection subspaces: Enrichment, augmentation, and  recycling; see \cref{sub:recyclEnrich}
  \item Beyond the 2-norm: Sparsity-enforcing hybrid projection methods for {$\ell_p$ regularization}; see \cref{sub:flexible}
  \item Beyond deterministic inversion: Hybrid projection methods in a Bayesian setting; see  \cref{sub:bayesian}
  \item Beyond linear forward models: Hybrid projection methods for nonlinear inverse problems; see \cref{sub:nonlinear}
\end{itemize}
 Pointers to relevant software packages are provided in \cref{sec:software}.  Conclusions and outlook on future directions are provided in \cref{sec:outlook}.

\paragraph{Notations} Boldface lower-case letters denote vectors: e.g., $\bfc\in\bbR^n$; $[\bfc]_i$, \linebreak[4]$1\leq i\leq n$ denotes the $i$th component of the vector $\bfc$. Boldface upper-case letters denote matrices: e.g., $\bfC\in\bbR^{m\times n}$; $c_{i,j}$, $1\leq i\leq m$, $1\leq j\leq n$, denotes the $(i,j)$th entry of the matrix $\bfC$. Indexed matching boldface lower-case letters denote matrix columns: e.g., $\bfc_j\in\bbR^{m}$, $j=1,\dots,n$ denotes the $j$th column of the matrix $\bfC$. $\bfI_p$ denotes the identity matrix of size $p$, whose columns $\bfe_j$, $j=1,\dots,p$ are the canonical basis vectors of $\bbR^p$. The column space of a matrix $\bfC$ is denoted by ${\rm ran}(\bfC)$. Lower-case Greek letters, sometimes indexed, are often used to denote scalars. In the following, when there is no ambiguity (essentially until \cref{sec:extensions}), the shorthand notation $\norm[]{\cdot}=\norm[2]{\cdot}$ will be used.

\section{General problem setup and background}
\label{sec:background}
For the simplified case where
$\calJ(\bfb-\bfA\bfx) = \norm[]{\bfb-\bfA\bfx}^2$, $\calR(\bfx) = \norm[]{\bfx}^2$, and $\calC=\bbR^n$ in \cref{eq:directreg}, we get the so-called standard-form Tikhonov problem,
\begin{equation}
	\label{eq:Tikhonov}
	\min_\bfx \norm[]{\bfb - \bfA \bfx}^2 + \lambda \norm[]{\bfx}^2.
	\end{equation}
We stress that the standard-form Tikhonov problem has been widely studied in both the mathematics and statistics communities and has been used in many scientific applications.  However, computing Tikhonov-regularized solutions can still be challenging if the size of $\bfx$ is very large or if $\lambda$ is not known a priori.  In fact, being able to efficiently compute solutions to \cref{eq:Tikhonov} was a main motivation for much of the early works on hybrid projection methods. For a discussion on extensions of hybrid methods to solve more general problems including the general-form Tikhonov problem, see \cref{sec:extensions}.

\subsection{SVD-based direct regularization methods}
\label{sub:direct}
In this section, we begin with the standard-form Tikhonov problem \cref{eq:Tikhonov}, and describe a direct approach based on the singular value decomposition (SVD) of $\bfA$ to compute a solution.  For the problems of interest, direct application of this approach is not computationally feasible.  However, the formulations briefly explored here will be useful for analysis and interpretation later in this paper. Furthermore, in the hybrid framework, the direct methods described here can be used to solve the projected problem. We point the interested reader to \cite{hansen2010discrete} for a more thorough exposition of direct methods.

Without loss of generality, let us assume that $\bfA\in \bbR^{m\times n}$ with ${\rm rank}(\bfA) = n \leq m$. The SVD of $\bfA$ is defined as
\begin{equation}\label{eq:svdA}
  \bfA = \bfU^\text{\tiny A} \bfSigma^\text{\tiny A} (\bfV^\text{\tiny A})\t
\end{equation}
where $\bfU^\text{\tiny A} = \begin{bmatrix} \bfu_1^\text{\tiny A} & \ldots & \bfu_m^\text{\tiny A} \end{bmatrix} \in \bbR^{m \times m}$ and $\bfV^\text{\tiny A}= \begin{bmatrix} \bfv_1^\text{\tiny A} & \ldots & \bfv_n^\text{\tiny A} \end{bmatrix} \in \bbR^{n \times n}$ are orthogonal and $\bfSigma^\text{\tiny A} = \diag{\sigma_1^\text{\tiny A}, \ldots, \sigma_n^\text{\tiny A}} \in \bbR^{m \times n}$ contains the singular values, $\sigma_1^\text{\tiny A}\geq\sigma_2^\text{\tiny A}\geq\dots\geq\sigma_n^\text{\tiny A}>0$.

The solution to \cref{eq:Tikhonov} can be written as
\[\bfx(\lambda) =\sum_{i=1}^{n} \phi_i(\lambda) \frac{(\bfu_i^\text{\tiny A})\t \bfb}{\sigma_i^\text{\tiny A}} \bfv_i^\text{\tiny A},\quad
	\mbox{where}\quad \phi_i(\lambda) = \frac{(\sigma_i^\text{\tiny A})^2}{(\sigma_i^\text{\tiny A})^2 + \lambda}
\]
are the Tikhonov filter factors.  Tikhonov regularization is one example from a wider class of spectral filtering methods. SVD-filtering methods compute regularized solutions by imposing suitable filtering on the SVD components of the solution: namely, the filter factors $\phi_i$ should be close to $1$ for small $i$, and should approach $0$ for large $i$. The Tikhonov filter factors $\phi_i(\lambda)$ have this property, with the amount of filtering prescribed by the regularization parameter $\lambda$. Indeed, for small values of $\lambda$, little weight is put on the regularization term and the filter factors are more prone to approach $1$ (or, at least, a good portion of filter factors are close to $1$, depending on the location of $\lambda$ within the range of the singular values of $\bfA$). Thus, regularized solutions may be `under-regularized', i.e., still contaminated with noise.  However, for large values of $\lambda$, too much weight is put on the regularization term and the filter factors are more prone to approach $0$, so solutions may be `over-regularized' or `over-smoothed'.  We refer the reader to the images provided in \cref{fig:param_iterative} for an illustration of the impact of the regularization parameter $\lambda$ on the solution.

Also the truncated SVD (TSVD) method, which regularizes \cref{eq:linearmodel} by computing
\[
\bfx(k)=\sum_{i=1}^{k} \frac{(\bfu_i^\text{\tiny A})\t \bfb}{\sigma_i^\text{\tiny A}} \bfv_i^\text{\tiny A}=\sum_{i=1}^{n} \phi_i(k) \frac{(\bfu_i^\text{\tiny A})\t \bfb}{\sigma_i^\text{\tiny A}} \bfv_i^\text{\tiny A},
\]
 is a filtering method, with filter factors obtained by first setting a truncation index $k\in\{1,\dots,n\}$, and by then considering $\phi_i(k)=1$, if $i\leq k$, and $\phi_i(k)=0$ otherwise. Note that choosing $k=n$ returns the unregularized solution of \cref{eq:linearmodel}. Therefore, $k$ plays the role of regularization parameter.

Every spectral filtering method requires the selection of a regularization parameter.
Common strategies to do so include the discrepancy principle \cite{Mor66}, the generalized cross-validation (GCV) method \cite{GoHeWh79}, the unbiased predictive risk estimation (UPRE) method \cite{UPREoriginal}, and the $L$-curve \cite{LaHa95}. Since there is not one method that will work for all problems, it is usually a good idea to try a variety of methods. For problems where the SVD of $\bfA$ is available, one can efficiently apply these regularization parameter strategies. For problems where computing the SVD is not feasible, parameter choice strategies are limited. However, we will see in \cref{sec:param} that many of these existing regularization parameter selection techniques are not only feasible but also can be successfully integrated within hybrid projection methods.

The SVD of $\bfA$ \cref{eq:svdA}, besides being an essential building block of SVD-filtering methods, is a pivotal tool for the analysis of discrete inverse problems \cite{engl1996regularization,Han97,hansen2010discrete}. For instance, looking at the decay of the singular values, one may infer different degrees of ill-posedness: typically, following \cite{Hof93}, a polynomial decay of the form $\sigma_i^\text{\tiny A}=ci^{-\alpha}$, $c>0$, is classified as mild (if $0<\alpha\leq 1$) or moderate (if $\alpha>1$) ill-posedness, while an exponential decay of the form $\sigma_i^\text{\tiny A}=\exp(-\alpha i)$, $\alpha>0$, is classified as severe ill-posedness. The Picard condition, which in the continuous setting provides necessary and sufficient conditions for the existence of a solution of the form \cref{xpseudoinv} \cite{engl1996regularization}, can be adopted in a discrete setting to gain insight into the relative behavior of the magnitude of $(\bfu_i^\text{\tiny A})\t \bfb$ and $\sigma_i^\text{\tiny A}$, which appear in the expression of the (filtered) solutions. Employing tools such as the so-called `Picard plot' \cite{hansen2010discrete} 
can inform us on the existence of a meaningful solution to \cref{eq:linearmodel} as well as the effectiveness of any (filtering) regularization method applied to \cref{eq:linearmodel}.

\subsection{Iterative regularization methods}
\label{sub:iterative}
Next we describe iterative regularization, where an iterative method is used to approximate the solution of the \emph{unregularized} least-squares problem,
\begin{equation}\label{LS}
	\min_\bfx \norm[]{\bfb - \bfA \bfx}^2\,,
\end{equation}
and early termination of the iterative process results in a regularized solution.  We focus on projection methods where the underlying concept is to constrain the solution at the $k$th iteration to lie in a $k$-dimensional subspace spanned by the (typically orthonormal) columns of some matrix $\bfV_k = \begin{bmatrix}
	\bfv_1 & \cdots & \bfv_k
\end{bmatrix}$, where $\bfv_i\in \bbR^n$. That is, the regularized solution is given as
\begin{equation}
\bfx_k = \bfV_k \bfy_k\,,\quad\mbox{where}\quad
	\label{eq:projectedproblem}
\bfy_k =	\argmin_{\bfy\in\bbR^k} \norm[]{\bfb - \bfA \bfV_k \bfy}^2.
\end{equation}

We are interested in projection methods using Krylov subspaces \cite{saad2003iterative}. Given $\bfC\in \bbR^{n\times n}$ and $\bfd\in \bbR^n$, a Krylov subspace is defined as
\begin{equation*}
	\calK_k(\bfC,\bfd) = {\rm span}\{\bfd, \bfC \bfd, \bfC^2 \bfd, \ldots , \bfC^{(k-1)}\bfd\}\,.
\end{equation*}
Here and in the following, we assume that the dimension of $\calK_k(\bfC,\bfd)$ is $k$. When Krylov methods are employed to solve problem \cref{LS}, $\bfC$ and $\bfd$ are defined in terms of $\bfA$ and $\bfb$, respectively.
In iterative regularization methods based on the Arnoldi process (such as GMRES), the columns of $\bfV_k$ form an orthonormal basis for $\calK_k(\bfA,\bfb)$, where $\bfA$ must be square.
However, for problems where $\bfA$ is not square (e.g. in tomography), iterative regularization methods based on the Lanczos or Golub-Kahan bidiagonalization process (such as LSQR \cite{PaSa82a,PaSa82b} or LSMR \cite{fong2011lsmr}) are used, and the columns of $\bfV_k$ form an orthonormal basis for $\calK_k(\bfA\t\bfA,\bfA\t\bfb)$. Note that LSQR is mathematically equivalent to CGLS, i.e., the conjugate gradient (CG) method applied to the normal equations.
All the Krylov methods mentioned so far are minimum residual methods, so that the $k$th iterate can be written as $\bfx_k = \bfV_k \bfy_k$, where $\bfy_k$ solves \cref{eq:projectedproblem}. More details regarding the generation of these basis vectors will be addressed in \cref{sub:buildsubspace}.

\begin{figure}[bt!]
\begin{center}
\includegraphics[width=.9\textwidth]{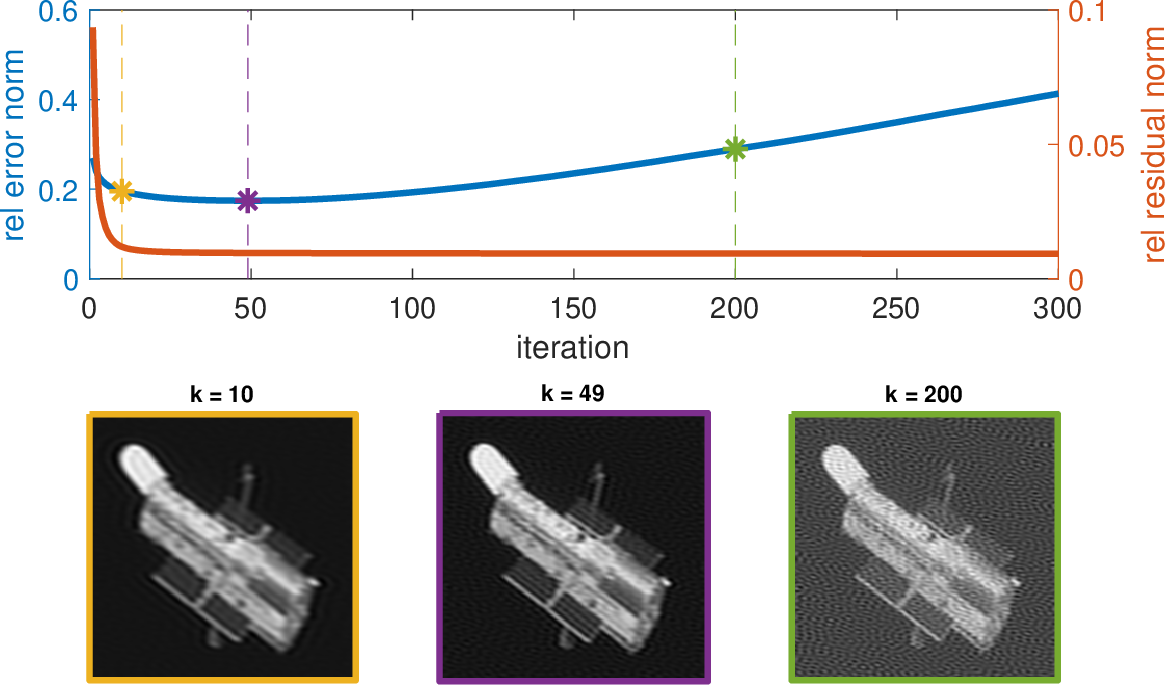}
\end{center}
\caption{Illustration of the semi-convergence phenomenon for the image deblurring test problem.  In the plot we provide the relative reconstruction error norms $\|\bfx_k-\bfx_\true\|/\|\bfx_\true\|$ per iteration $k$ and the relative residual norms $\|\bfb-\bfA\bfx_k\|/\|\bfb\|$ per iteration $k$.  Image reconstructions at iterations $10$, $49$ and $200$ are provided and correspond to the stars in the plot.}
\label{fig:semiconvergence}
\end{figure}

A crucial feature of classical Krylov methods when applied to ill-posed inverse problems is that, oftentimes, these iterative methods exhibit a regularizing effect in that the projection subspace in early iterations provides a good basis for the solution.
One of the first Krylov methods that was proven to have regularizing properties is CGLS
\cite{nemirovskii1986regularizing, hanke1995minimal}.
Indeed, the CGLS iterates can be expressed as filtered SVD solutions \cite{hansen2010discrete}.
Since these iterative methods converge to the least-squares solution, we get a phenomenon commonly referred to as `semi-convergence' (after \cite{Nat01}), whereby the relative reconstruction errors decrease at early iterations but increase at later iterations due to noise manifestation and amplification.  See \cref{fig:semiconvergence} for an illustration.
Because of this, for iterative regularization, the stopping iteration plays the role of the regularization parameter. There have been many investigations into developing stopping criteria for iterative methods for inverse problems.  If a good estimate of the amount of noise is available, the most widely used and intuitive approach is to stop iterating as soon as the value of the objective function in \cref{eq:projectedproblem}
reaches the magnitude of the noise: this is the so-called discrepancy principle.
GCV and modifications of GCV have also been used for stopping criteria \cite{Bjo88,BjGrDo94,ChNaOL08}.  We again refer to \cref{sec:param} for more details about these strategies.
The main challenge is that for some problems, the reconstruction can be very sensitive to the choice of the stopping iteration.  Thus, if the method stops a little too late, the reconstruction is already contaminated by noise (i.e., under-regularized). On the contrary, if the iterations are stopped too early (i.e., over-regularization), a potentially better reconstruction is precluded. We remark that it has been observed that semi-convergence may appear somewhat `prematurely' and that it is sometimes important to have a larger approximation subspace (which would otherwise be beneficial for the solution; see \cite{huang2017some,OlSi81}).

In some circumstances, even if direct regularization methods are feasible, one may prefer to adopt iterative regularization methods. Some reasons are highlighted in \cite{Hanke01}, where the following heuristic motivation is given: contrary to TSVD, Krylov methods (such as LSQR) generate an approximation subspace that is tailored to the current right-hand-side vector. Therefore, the basis vectors for the solution may be better `adapted' to the given problem than the right singular vectors.

In many fields from numerical linear algebra to differential equations, iterative methods and, in particular, preconditioned Krylov methods, have been immensely successful in solving large, sparse systems of equations efficiently \cite{saad2003iterative}. Iterative methods have also gained widespread use in the inverse problems community.  One of the main reasons for this is that neither the matrix $\bfA$ nor its factorization need to be constructed, and thus, these methods are ideal for large-scale problems. Furthermore, it has been observed many times that the generated Krylov subspaces are rich in information for representing the solution (i.e., corresponding to the large singular vectors).  Thus, a reasonable solution can be obtained in only a few iterations. The greatest caveat for inverse problems is semi-convergence, so great care must be taken to find good stopping criteria. In \cref{sec:theory} we will dwell more on the properties of the approximation subspaces generated by different projection methods.

\subsection{Iterative methods for solving the Tikhonov problem}\label{sub:IterTikh}

As described in \cref{sub:direct}, the Tikhonov problem \cref{eq:Tikhonov} can be easily solved if the SVD of $\bfA$ is available and, in this case, many well-regarded parameter choice strategies can be applied to compute a suitable value of the regularization parameter $\lambda$. Nonetheless, for large-scale problems where $\bfA$ cannot be constructed but matrix-vector products with $\bfA$ and $\bfA\t$ can be computed efficiently, iterative projection methods \cite{PaSa82a,PaSa82b} can be used to solve the equivalent Tikhonov problem,
\begin{equation}\label{eq:TikhonovLS}
\min_\bfx \norm[]{\left[\begin{array}{c}
\bfA \\
\sqrt{\lambda}\bfI_n
\end{array}\right]
\bfx - \left[\begin{array}{c}
\bfb \\
\bfzero
\end{array}\right]}^2.
\end{equation}
However, this approach may not be convenient if a suitable value of $\lambda$ is not known a priori.
 In this case, often one must solve problem \cref{eq:TikhonovLS} from scratch for many different values of $\lambda$ and this eventually results in an expensive approach (note that for specific iterative solvers, one may adopt smart strategies to reduce computations; see \cite{FrMa99}).

\begin{figure}[bt!]
\begin{center}
\includegraphics[width=.9\textwidth]{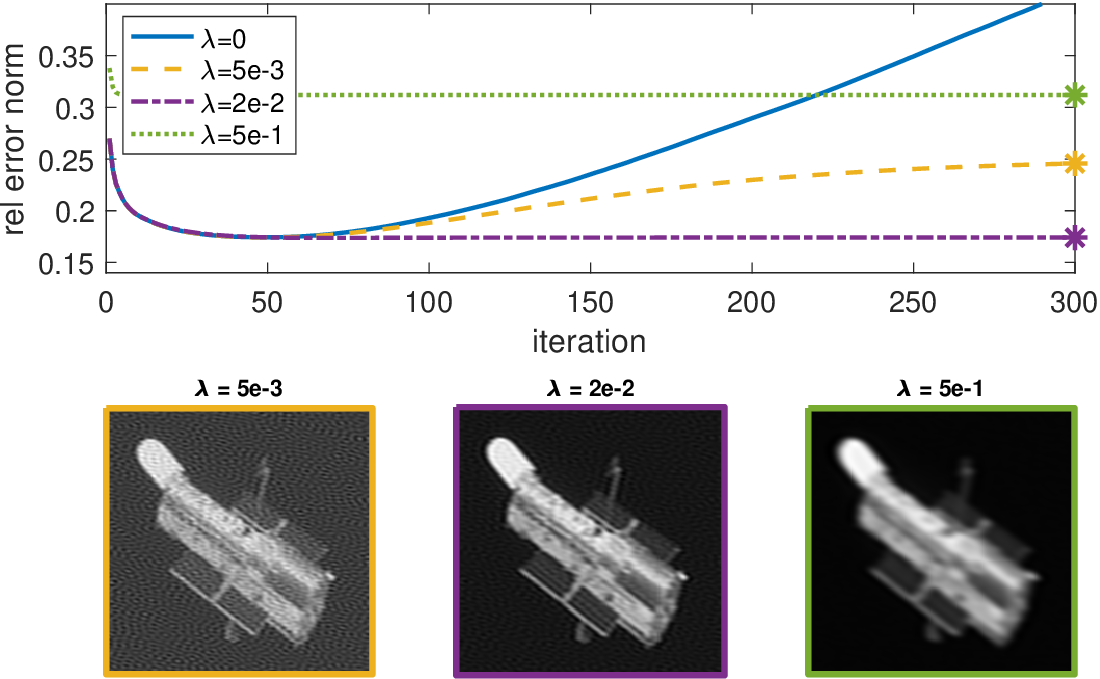}
\end{center}
\caption{Illustration of the interplay of the regularization parameter and the number of iterations when using an iterative method to solve the Tikhonov problem with a \emph{fixed} regularization parameter. We provide the history of the relative reconstruction error norms for \cref{eq:TikhonovLS} versus the number of (CGLS) iterations, for different values of the regularization parameter $\lambda$. Image reconstructions at iteration $300$ are provided for different choices of $\lambda$ (these correspond to the stars on the error plots).  Notice that for too small values of $\lambda$, semi-convergence can be still detected (i.e., the behavior is  similar to the `unregularized' $\lambda=0$ case, also displayed in \cref{fig:semiconvergence}). For too large values of $\lambda$, the iterative method can only produce smooth reconstructions that never achieve high accuracy.}
\label{fig:param_iterative}
\end{figure}

Similar to the discussion in \cref{sub:direct}, the value of $\lambda$ can have a considerable impact on the quality of the reconstructed solution. Moreover, when using iterative methods to solve the Tikhonov problem, one can to some extent leverage the number of iterations to enforce additional regularization. For illustration, we use CGLS to solve \cref{eq:TikhonovLS} for various choices of $\lambda$ for the image deblurring example and provide relative reconstruction error norms per iteration in \cref{fig:param_iterative}. Notice that if $\lambda$ is chosen too small, severe semi-convergence appears and a good stopping iteration is as crucial as for the `purely' iterative (i.e., $\lambda=0$) methods introduced in \cref{sub:iterative}; on the contrary, if $\lambda$ is chosen too large, the solution is over-regularized and additional iterations cannot mitigate this phenomenon. Thus, $\lambda$ seems to have a greater impact on the quality of the reconstructed solution than the number of iterations, which can be attributed to the modification of the spectral properties of the coefficient matrix in \cref{eq:TikhonovLS}.

We can draw the following conclusions about the successful working of a standard iterative method to solve a Tikhonov regularized problem. First, the user must be reasonably confident in the choice of the regularization parameter.
Second, the user must be cognizant about the interplay of the choice of the regularization parameter and the number of iterations for \cref{eq:TikhonovLS}, as both of these contribute to a suitable regularized solution. As we will see in the next section, hybrid projection methods provide an alternative approach of combining an iterative method and Tikhonov regularization, whereby a regularization parameter can be automatically and adaptively tuned during the iterative process.

\section{Hybrid projection methods}
\label{sec:hybrid}
A hybrid projection method is an iterative strategy that regularizes problem \cref{LS} by projecting it onto subspaces of increasing dimension and by solving the projected problem using a variational regularization method. Formally, recalling the framework for iterative projection methods unfolded in \cref{eq:projectedproblem}, the solution at the $k$th iteration of a hybrid projection method for standard Tikhonov can be represented as
\begin{equation}
	\label{eq:regprojproblem}
	\bfx_k(\lambda_k) = \bfV_k\bfy_k(\lambda_k)\,,\quad\mbox{where}\quad
	\bfy_k(\lambda_k)=\argmin_{\bfy \in \bbR^k} \calJ(\bfb - \bfA\bfV_k\bfy) + \lambda_k \norm[]{\bfy}^2.
\end{equation}
Notice that the solution incorporates both regularization from the projection subspace (determined by the choice of the projection method, i.e., $\bfV_k$ and $k$) and a potentially changing variational regularization term (defined by $\lambda_k$).
Again, for simplicity and for historical reasons, in this section we focus on the standard-form Tikhonov problem \cref{eq:Tikhonov}; extensions will be considered in \cref{sec:extensions}.

\begin{figure}[bt!]
\begin{center}
\includegraphics[width=.9\textwidth]{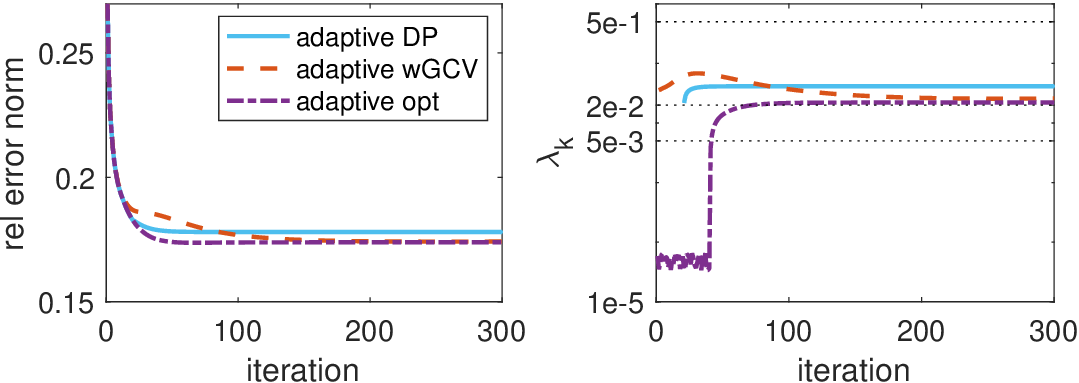}
\end{center}
\caption{Illustration of a hybrid projection method where an iterative method is run with an \emph{adaptive} regularization parameter selection method. In the left panel, we provide the history of the relative reconstruction error norms per iteration $k$, for $\lambda=\lambda_k$ set adaptively at each iteration using the discrepancy principle (DP) and weighted GCV (wGCV); see \cref{sec:param} for more details.  The optimal (opt) parameter corresponds to selecting the $\lambda_k$ that minimizes the relative error at each iteration $k$. In the right panel, we provide the computed regularization parameters, where the three horizontal lines correspond to $\lambda = 5\cdot 10^{-3}, 2\cdot 10^{-2}$ and $5\cdot 10^{-1}$.}
\label{fig:param_adaptive}
\end{figure}

\paragraph{An illustration} Before we get into the details of hybrid projection methods, we begin by illustrating the benefits of allowing \emph{adaptive} choices for the regularization parameter, where a different regularization parameter can be used at each iteration, i.e., $\lambda=\lambda_k$. In \cref{fig:param_adaptive} for the image deblurring problem, we provide relative reconstruction error norms per iteration in the left panel and computed regularization parameters per iteration in the right panel.  Here, `opt' refers to selecting at each iteration the regularization parameter that delivers the smallest relative reconstruction error. Parameter selection methods `DP' and `wGCV' and others will be discussed in \cref{sec:param}. Notice that relative reconstruction errors using (appropriately selected) adaptive regularization parameters can overcome semi-convergence behavior (see \cref{fig:semiconvergence}) and result in reconstruction errors that are close to a pre-selected optimal regularization parameter. For early iterations, a small $\lambda_k$ (corresponding to little variational regularization) is sufficient, but as iterations progress, the adaptive methods should select regularization parameters close to the optimal one.

\paragraph{An algorithm} Now that we have seen the potential benefits of a hybrid projection method, let us discuss the general structure of such algorithms.  A sketch of a hybrid projection method for standard Tikhonov regularization is provided in \cref{alg:hybridprojection}, with links to the appropriate sections of the paper for more details.  Notice that each iteration requires the expansion of the projection subspace and the solution of a projected regularized problem where a regularization parameter can be selected at each iteration.
\begin{algorithm}
\caption{Hybrid projection method for standard Tikhonov}
\label{alg:hybridprojection}
\begin{algorithmic}[1]
	\REQUIRE{$\bfA$, $\bfb$, $k=1$, projection method, regularization parameter selection method}
\WHILE {stopping criterion not satisfied}
\STATE Expand the projection subspace ${\rm ran}(\bfV_k)$;
see \cref{sub:buildsubspace}
\STATE Select a regularization parameter $\lambda_k$;
see \cref{sec:param}
\STATE Solve the projected regularized problem \cref{eq:regprojproblem};
see \cref{sub:solvingprojected}
 \STATE $k = k+1$
\ENDWHILE
\ENSURE{$\bfx_k(\lambda_k)$}
\end{algorithmic}
\end{algorithm}

An outline of the remaining part of this section is as follows. We begin in \cref{sub:history} with a brief historical overview of hybrid projection methods. Details and derivations for hybrid projection methods can be found in \cref{sub:hybriddetails}, where we describe some techniques for steps 2 and 4 of \cref{alg:hybridprojection}: namely, expanding the projection subspace (e.g., via the Arnoldi or Golub-Kahan process) and solving the projected regularized problem.  One of the main advantages and features of a hybrid projection method is the ability to select regularization parameters automatically and adaptively (i.e., during the iterative process as in step 3 of \cref{alg:hybridprojection}).  Thus, we dedicate \cref{sec:param} to providing an overview of regularization parameter selection methods, with a particular emphasis on methods that have been critical for the success of hybrid projection methods.  Although we still do not have a complete analysis of the regularizing properties of every hybrid projection method, important theoretical results and properties have nonetheless been established for specific methods, and we highlight some of these in \cref{sec:theory}.

\subsection{Historical development of hybrid methods}
\label{sub:history}
Since the seminal publication by O'Leary and Simmons in 1981 \cite{OlSi81}, there have been four decades of progress and developments in the field of hybrid projection methods.  A quick google scholar search shows that the number of citations for this paper alone nearly doubles in each subsequent decade. In this section, we provide a brief overview of the main contributions and highlights by decade.  This is by no means an exhaustive list of publications, and we realize that there may be bias in the selections, but our goal is to provide the reader with some historical context.  For a novice reader, this section can be skipped upon first reading.

\paragraph{1981-1990} To the best of our knowledge, the first journal publication to introduce a hybrid projection method was by O'Leary and Simmons in 1981 \cite{OlSi81}, in which they describe a `projection-regularization' method that combines the Golub-Kahan bidiagonalization for projection and TSVD for regularization.  In the paper, the authors mention that a different hybrid algorithm was proposed independently, with a reference to a technical report by Bj{\"o}rck in 1980 \cite{bjorck80}.
Bj{\"o}rck's work appeared in BIT in 1988 \cite{Bjo88}, where a key difference was that the right hand side vector $\bfb$ was used as the starting vector, since it was noted that allowing the algorithms to start with an arbitrary vector did not always perform well.  Other algorithmic contributions included the use of cross-validation to determine stopping criteria and transformations to standard form based on Eld{\' e}n's earlier work \cite{elden1977}. Example problems from time series deconvolution were used in \cite{OlSi81}, but the computational technology at the time was still quite limited \cite{simmons1986deconvolution}.

\paragraph{1991-2000} During this period, hybrid projection methods gained significant traction in the numerical linear algebra community as well as practical utility in various seismic imaging applications \cite{yao1999,sacchi1998}. As described in Bj{\"o}rck's book \cite{bjorck1996numerical}, various researchers were interested in characterizing the regularizing properties of Krylov methods (e.g., see Nemirovskii \cite{nemirovskii1986regularizing} and various works by Hanke and Hansen \cite{hanke1993regularization}), where the main motivation was to determine appropriate stopping criteria for iterative methods when applied to ill-posed problems.  Advancements in hybrid methods included investigations into stable computations (e.g., reorthogonalization) by Bj{\" o}rck, Grimme and Van Dooren \cite{BjGrDo94} as well as new methods for selecting regularization parameters.  For example, Calvetti, Golub and Reichel \cite{CaGoRe99} proposed a hybrid approach where a so-called `$L$-ribbon' was computed using a partial Lanczos bidiagonalization.
The first application of hybrid methods for large, nonlinear inversion was described in Haber and Oldenburg \cite{haber2000gcv}, and hybrid approaches based on projections with GMRES or Arnoldi were described by Calvetti, Morigi, Reichel, and Sgallari in \cite{CaMoReSg00}. In terms of software, Hansen laid some groundwork for hybrid methods in the \textsc{Regularization Tools} package \cite{Han94}, where routines for computing the lower bidiagonal matrix and the SVD of a bidiagonal matrix were provided.

Not all of the work on hybridized methods during this time followed the project-then-regularize framework.  Related work on estimating the regularization parameter for large-scale Tikhonov problems by exploiting relationships with CG also appeared during this time. For example, Golub and Von Matt \cite{golub1997generalized} estimated the GCV parameter for large problems using relations between Gaussian quadrature rules and Golub-Kahan bidiagonalization, and Frommer and Maass \cite{FrMa99} exploited the shift structure of CG methods to solve Tikhonov problems with multiple regularization parameters simultaneously.

\paragraph{2001-2010}
With significant improvements in computational resources and capabilities, researchers became interested in how to realize the benefits of hybrid projection methods in practice. In this decade, we saw a surge in the development of regularization parameter selection methods for hybrid projection methods \cite{KiOLe01,calvetti2002curve,ChNaOL08,mead2008newton,bazan2010gkb} and various extensions of hybrid methods to include general regularization terms (e.g., total variation \cite{calvetti2002hybrid} and general-form Tikhonov \cite{KiHaEs06, hochstenbach2010iterative}).  Still, the main focus was on TSVD and Tikhonov for regularization.  New interpretations of hybrid methods based on Lanczos and TSVD were described in \cite{Hanke01}, and noise level estimation from the Golub-Kahan bidiagonalization process, which can be used for determining the stopping iteration, were described in \cite{HnPlSt09}. A fully automatic MATLAB routine called `HyBR' was provided in \cite{ChNaOL08}, where a Golub-Kahan hybrid projection method with a weighted GCV parameter choice rule was implemented for standard Tikhonov regularization. Hybrid methods were implemented on distributed computing architectures \cite{ChStYa09} and were used for many new applications such as cryoelectron microscopy and electrocardiography \cite{jiang2008two}.

\paragraph{2011-present}  At present, hybrid projection methods have gained significant interest in many research fields, and contributions range from new methodologies and advanced theories to innovations in scientific applications. Many of the recent developments in hybrid methods will be elaborated on in \cref{sec:extensions}, so here we just provide some highlights.  In particular, there have been many papers on the Arnoldi-Tikhonov method \cite{GaNo14,GaNoRu14}.  Flexible and generalized hybrid methods based on state-of-the-art Krylov subspace methods have been introduced for including more general regularization terms and constraints \cite{lampe2012large, reichel2012tikhonov, gazzola2014generalized,ChSa17,chung2019, GaSL18}.  There have been connections made to the field of computational uncertainty quantification \cite{saibaba2020efficient} and new insights in the regularization parameter selection strategies \cite{ReVaAr17}, and regularization by approximate matrix functions \cite{fcnM1, cipolla2021regularization, fcnM2}.  There exists a plethora of papers in application and imaging journals, and a new software package called \textsc{IR Tools} \cite{IRtools}.  We also refer interested readers to an older survey paper on Krylov projection methods and Tikhonov regularization \cite{gazzola2015survey}.

\subsection{Algorithmic approaches to hybrid projection methods}
\label{sub:hybriddetails}
Following derivations of iterative methods for sparse rectangular systems, in \cref{sub:buildsubspace} we first describe
common techniques to iteratively build a solution subspace and to efficiently project problem \cref{LS} onto subspaces of increasing dimension. Then, in \cref{sub:solvingprojected}, we describe techniques to solve the projected regularized problem.

\subsubsection{Projecting onto subspaces of increasing dimension}
\label{sub:buildsubspace}
In this section, we address the first building block of a hybrid projection method, which is to use an iterative method to build a sequence of solution subspaces of small but increasing dimension, and to efficiently project the problem onto them. 
The first important question to address is: in practice, what makes a good solution subspace?
Indeed, a clever choice of basis vectors should satisfy multiple requirements. Since the number of iterations required for solution computation should be small, a good basis should be able to accurately capture the important information in a few vectors.  In this sense, the SVD basis vectors (i.e, the columns of  $\bfV^\text{\tiny A}$ in \cref{eq:svdA}) would be ideal (recall discussion about filtered SVD in \cref{sub:direct}).  Unfortunately, computing the SVD can be infeasible for very large problems.  Thus, we seek a subspace sequence that exhibits a rapid enough decrease in the residual norm
$\|\bfr\|^2=\norm[]{\bfb - \bfA \bfx}^2$ in \cref{LS} so that early termination provides good approximations \cite{golub2012matrix}.
Another desirable property is that the solution subspace can be generated efficiently.  Here we mean that the main computational cost per iteration is manageable, i.e, requires one matrix-vector multiplication with $\bfA$ and possibly one with $\bfA\t$. Lastly, for numerical stability it may be desirable to have an orthogonal basis for the solution subspace \cite{o2009scientific}.
Although a variety of projection methods and solution subspaces could be used, here we focus on two of the most common projection methods on Krylov subspaces, namely the Arnoldi process and the Golub-Kahan bidiagonalization process. {We will argue that these algorithms compute good solution spaces, according to the criteria listed above.}

\paragraph{Arnoldi process} The Arnoldi process generates an orthonormal basis for $\calK_k(\bfA,\bfb)$ at the $k$th iteration \cite{saad2003iterative}.
Given matrix $\bfA \in \bbR^{n \times n}$ and vector $\bfb \in \bbR^n,$ with initialization $\bfv_1 = \bfb/\beta_1$ where $\beta_1 = \norm[]{\bfb}$, the $k$th iteration of the Arnoldi process generates vector $\bfv_{k+1}$ and scalars
$h_{i,k}$ for $i=1,\ldots,(k+1)$
such that, in matrix form, the following Arnoldi relationship holds,
\begin{equation}
	\label{eq:Arnoldi}
\bfA \bfV_k = \>\bfV_{k+1} \bfH_k,
\end{equation}
where
\[ \bfH_k := \> \begin{bmatrix}
h_{1,1} & h_{1,2} & \cdots & h_{1,k} \\
h_{2,1} & h_{2,2} & \cdots & h_{2,k} \\
0 & \ddots & \ddots & \vdots\\ \vdots & \ddots & \ddots & h_{k,k} \\ 0 & \cdots& 0 & h_{k+1,k}
\end{bmatrix}\in\bbR^{(k+1)\times k}, 
\bfV_{k+1} := [\bfv_1,\dots,\bfv_{k+1}]\in\bbR^{n\times (k+1)},\]
and where, in exact arithmetic, $\bfV_{k+1}\t \bfV_{k+1} = \bfI_{k+1}.$ Notice that $\bfV_{k+1}$ contains an orthonormal basis for $\calK_{k+1}(\bfA,\bfb)$ and $\bfH_{k}$ 
is an upper Hessenberg matrix (i.e., upper triangular with the first lower diagonal).
The Arnoldi process is summarized in \cref{alg:Arnoldi}.

\begin{algorithm}
\caption{Arnoldi process}
\label{alg:Arnoldi}
\begin{algorithmic}[1]
	\REQUIRE{$\bfA\in \bbR^{n \times n}$, $\bfb \in \bbR^n$, $k$}
\STATE $\bfv_1=\bfb/\norm{\bfb}$ 
\FOR {$j=1, \dots, k$}
\STATE $\bfv = \bfA\bfv_j$
\STATE \textbf{for} $i=1,\dots,j$ \textbf{do} $h_{i,j} = \bfv\t \bfv_i$ \textbf{end for}
\STATE $\bfv = \bfv-\sum_{i=1}^j h_{i,j}\bfv_i$
\STATE $\bfv_{j+1} = \bfv/h_{j+1,j}$ \, where \, $h_{j+1,j} = \norm[]{\bfv}$
\ENDFOR
\ENSURE{$\bfV_k$, $\bfH_k$}
\end{algorithmic}
\end{algorithm}

Notice that the main computational cost of each iteration of the Arnoldi process is one matrix-vector multiplication with $\bfA$.  Breakdown of the iterative process may happen
(when $h_{j+1,j}=0$), but we are typically not concerned about this because these methods are meaningful only when the number of iterations is not too high.

One caveat, especially for inverse problems such as tomography, is that the Arnoldi method relies on the fact that $\bfA$ is square.  For rectangular problems, one could apply the Arnoldi process to the normal equations $\bfA\t \bfA \bfx = \bfA\t \bfb$, in which case the Arnoldi algorithm simplifies to the symmetric Lanczos tridiagonalization process where, thanks to symmetry, $\bfH_k$ is tridiagonal rather than upper Hessenberg. However, a numerically favorable approach would be to avoid the normal equations and to work directly with $\bfA$ and $\bfA\t$ separately. This is achieved by adopting the Golub-Kahan bidiagonalization process, which we describe next; we will also state more precisely the links between the Golub-Kahan bidiagonalization and the symmetric Lanczos process.

\paragraph{Golub-Kahan bidiagonalization process}
Contrary to the Arnoldi process that generates only one orthonormal basis (i.e, for $\calK_k(\bfA,\bfb)$), the Golub-Kahan bidiagonalization (GKB) process generates two sets of orthonormal vectors that span Krylov subspaces $\calK_k(\bfA\t \bfA,\bfA\t\bfb)$ and $\calK_k(\bfA\bfA\t,\bfb)$ \cite{golub1965calculating}. Given a matrix $\bfA \in \bbR^{m \times n}$ and a vector $\bfb \in \bbR^m,$ with initialization $\bfu_1 = \bfb/\beta_1$ where $\beta_1 = \norm[]{\bfb}$ and $\bfv_1 = \bfA\t \bfu_1 / \alpha_1$ where $\alpha_1 = \norm[]{\bfA\t\bfu_1}$, at the $k$th iteration of the GKB process, we generate vectors $\bfu_{k+1}$, $\bfv_{k+1}$, and scalars $\alpha_{k+1}$ and $\beta_{k+1}$ such that in matrix form, the following relationships hold,
\begin{equation}
	\label{eq:GKB}
\begin{array}{lcl}
\bfA\bfV_{k}&\!\!\!=\!\!\!&\bfU_{k+1}\bfB_{k}\\
&\!\!\!=\!\!\!&\bfU_{k}\bfB_{k,k} + \beta_{k+1}\bfu_{k+1}\bfe_k\t
\end{array},
\;
\begin{array}{lcl}
 \bfA\t\bfU_{k+1}&\!\!\!=\!\!\!&\bfV_{k}\bfB_{k}\t + \alpha_{k+1} \bfv_{k+1} \bfe_{k+1}\t\\
&\!\!\!=\!\!\!&\bfV_{k+1}\bfB_{k+1,k+1}\t
\end{array},
\end{equation}
where
\begin{equation*}
	\bfB_k := \> \begin{bmatrix}
 \alpha_{1} &0 & \cdots & 0 \\
 \beta_{2} & \alpha_{2} & \ddots & \vdots \\
 0 & \beta_3 & \ddots & 0\\ \vdots & \ddots & \ddots & \alpha_{k} \\ 0 & \cdots& 0 & \beta_{k+1}
 \end{bmatrix} = \begin{bmatrix}
 \bfB_{k,k}\\
\beta_{k+1}\bfe_k\t
 \end{bmatrix}\in\bbR^{(k+1)\times k}\,,
\end{equation*}
is bidiagonal, $\bfU_{k+1} := [\bfu_1,\dots,\bfu_{k+1}]\in\bbR^{m\times (k+1)}$ and $\bfV_{k} := [\bfv_1,\dots,\bfv_{k}]\in\bbR^{n\times k}$.
In exact arithmetic $\bfU_{k+1}$ has orthonormal columns that span $\calK_{k+1}(\bfA\bfA\t,\bfb)$ and $\bfV_{k}$ has orthonormal columns that span $\calK_{k}(\bfA\t\bfA,\bfA\t\bfb)$. The GKB process is summarized in \cref{alg:GKbidiag}.

\begin{algorithm}[!ht]
\begin{algorithmic}[1]
	\REQUIRE{$\bfA$, $\bfb$, $k$}
\STATE $\beta_1 \bfu_1 = \bfb,$ where $\beta_1 = \norm[]{\bfb}$
\STATE $\alpha_1 \bfv_1 = \bfA\t \bfu_1,$ where $\alpha_1 =  \norm[]{\bfA\t \bfu_{1}}$
\FOR {$j=1, \dots, k$}
\STATE $\beta_{j+1}\bfu_{j+1} = \bfA\bfv_j - \alpha_j \bfu_j$, where $\beta_{j+1} = \norm[]{\bfA\bfv_j - \alpha_j \bfu_j}$
\STATE $\alpha_{j+1}\bfv_{j+1} = \bfA\t \bfu_{j+1} - \beta_{j+1} \bfv_j$, where $\alpha_{j+1} = \norm[]{\bfA\t \bfu_{j+1} - \beta_{j+1} \bfv_j}$
\ENDFOR
\ENSURE{$\bfU_{k+1}$, $\bfV_{k+1}$, $\bfB_k$}
\end{algorithmic}
\caption{Golub-Kahan bidiagonlization (GKB) process}
\label{alg:GKbidiag}
\end{algorithm}

The main computational cost of each iteration of the GKB process is one matrix-vector multiplication with $\bfA$ and one with $\bfA\t$.  Breakdown in the GKB algorithm would occur when
$\beta_{j+1}=0$ or $\alpha_{j+1}=0$
but, similarly to the Arnoldi process, we will assume that breakdowns do not happen.  A potential concern with the GKB process is loss of orthogonality in the vectors of $\bfU_{k+1}$ and $\bfV_k$, for which it may be necessary to use a reorthogonalization strategy to preserve convergence of the singular values.  Indeed, it has been shown that reorthogonalization of only one set of vectors is necessary \cite{simon2000low,barlow2013reorthogonalization}; more details are provided in \cref{sec:theory}.

As mentioned above, GKB \cref{eq:GKB} is related to the symmetric Lanczos decomposition.
Notice that, multiplying the first equation in \cref{eq:GKB} from the left by $\bfA\t$, and using the second equation in \cref{eq:GKB}, we obtain
\[
\bfA\t \bfA\bfV_{k} \!=\! \bfA\t \bfU_{k+1} \bfB_{k}  \!=\! \bfV_{k+1}\bfB_{k+1,k+1}\t\bfB_k=\bfV_k \underbrace{\bfB_{k}\t\bfB_{k}}_{=:\hat{\bfT}_{k,k}} + \alpha_{k+1}\beta_{k+1}\bfv_{k+1}\bfe_{k}\t.
\]
The above chain of equalities is the symmetric Lanczos algorithm applied to $\bfA\t \bfA$ and $\bfA\t \bfb$, generating an orthonormal basis for $\calK_{k+1}(\bfA\t \bfA,\bfA\t \bfb)$. Similarly, multiplying the second equation in \cref{eq:GKB} (written in terms of $\bfU_k$) on the left by $\bfA$, and using the first equation in \cref{eq:GKB}, we obtain
\begin{equation}\label{SymLanczos1}
\bfA \bfA\t \bfU_{k} \!=\! \bfA \bfV_{k} \bfB_{k,k}\t  \!=\! \bfU_{k} \underbrace{\bfB_{k,k}\bfB_{k,k}\t}_{=:\bfT_{k,k}} + \alpha_{k}\beta_{k+1}\bfu_{k+1}\bfe_{k}\t.
\end{equation}
The above chain of equalities is the symmetric Lanczos algorithm applied to $\bfA \bfA\t$ and $\bfb$, generating an orthonormal basis for $\calK_{k+1}(\bfA\bfA\t,\bfb)$.

\paragraph{Comparison of Solution Subspaces}
Although the choice of a projection method relies heavily on the problem being solved, there have been some investigations into which projection methods might be better suited for certain types of problems, especially in light of some recent work on flexible methods, general-form Tikhonov regularization, and matrix transpose approximation. {We remark that the discussion here also relates to the material in \cref{sec:theory}.}

For Arnoldi based methods, which build the approximation subspace starting from $\bfb$, a potential concern is the explicit presence of the (rescaled) noise vector among the basis vectors \cite{JeHa07}. For problems with very low noise level, this is not a significant drawback \cite{calvetti2002gmres}.  However, various methods have been developed to remedy this. For example, the symmetric range-restricted minimum residual method MR-II \cite{Hank95a} and the range-restricted GMRES (RRGMRES) method \cite{calvetti2001choice} discard $\bfb$ and seek solutions in the Krylov subspace $\calK_k(\bfA,\bfA \bfb)$. For various projection methods, Hansen and Jensen \cite{HaJe08} studied the propagation of noise in both the solution subspace and the reconstruction, noting the manifestation of band-pass filtered white noise as `freckles' in the reconstructions.

For both the Arnoldi and the GKB process, notice that each additional vector in the Krylov solution subspace is generated by matrix-vector multiplication with $\bfA$ or $\bfA\t \bfA$ respectively.  Therefore the vectors defining the Krylov solution space are iteration vectors of the power method for computing the largest eigenpair of a matrix, and hence they become increasingly richer in the direction of the dominant eigenvector of $\bfA$ or $\bfA\t \bfA$. For many inverse problems (and due to the discrete Picard condition), the orthonormal basis solution vectors generated in early iterations tend to carry important information about the low-frequency components (i.e., the right singular vectors that correspond to the large singular values) \cite{hansen2010discrete}.  Therefore, the Krylov subspaces considered in this section can be appropriate choices for use in hybrid projection methods.

It is important to note that the choice of Krylov subspace may be problem dependent.  For some image deblurring problems, it has been observed that Arnoldi-based methods may not be suitable since multiplication with $\bfA$ corresponds to recursive blurring of the observation, resulting in a poor solution subspace. Various modifications have been proposed. For example, for image deblurring problems with nonsymmetric blurs and anti-reflective boundary conditions, `preconditioners' were used to improve the quality of the computed solution when using Arnoldi-Tikhonov methods \cite{donatelli2015arnoldi}. For tomography applications where the coefficient matrix is often rectangular, it is natural to consider Golub-Kahan based methods; even if the matrix $\bfA$ modeling a tomographic acquisition problem happens to be square (because of a specific scanning geometry), it should be stressed that standard Arnoldi methods would not generate a meaningful approximation subspace for the solution. This is because multiplication with $\bfA$ maps images to sinograms, and subsequent mappings to sinograms are neither physically meaningful nor fit for image reconstruction. Nonetheless, a specific `preconditioned' GMRES method, with `preconditioners' modeled as backprojection operators, has been used for tomography in \cite{coban2014regularised}.

It should be noted that preconditioning ill-posed problems can be a tricky business
\cite{HaNaPl93}. Indeed, besides improving the solution subspace, preconditioning can also be cautiously used with the classical goal of accelerating convergence (provided that it does not exacerbate the semi-convergence phenomenon).
Flexible (right) preconditioners have also been considered as means to improve the solution subspace \cite{gazzola2014generalized,chung2019}.
Sometimes additional vectors (that are hopefully meaningful to recover known features of the solution) can be added to the approximation subspace, where the goal is not necessarily to accelerate convergence, but rather to improve the solution subspace.  This is the idea behind the so-called enriched or recycling methods (see, for instance, \cite{dong2014r3gmres,hansen2019hybrid}); these are described in more detail in \cref{sub:recyclEnrich}.

\begin{figure}[bthp]
	\centering
	\includegraphics[width=.9\textwidth]{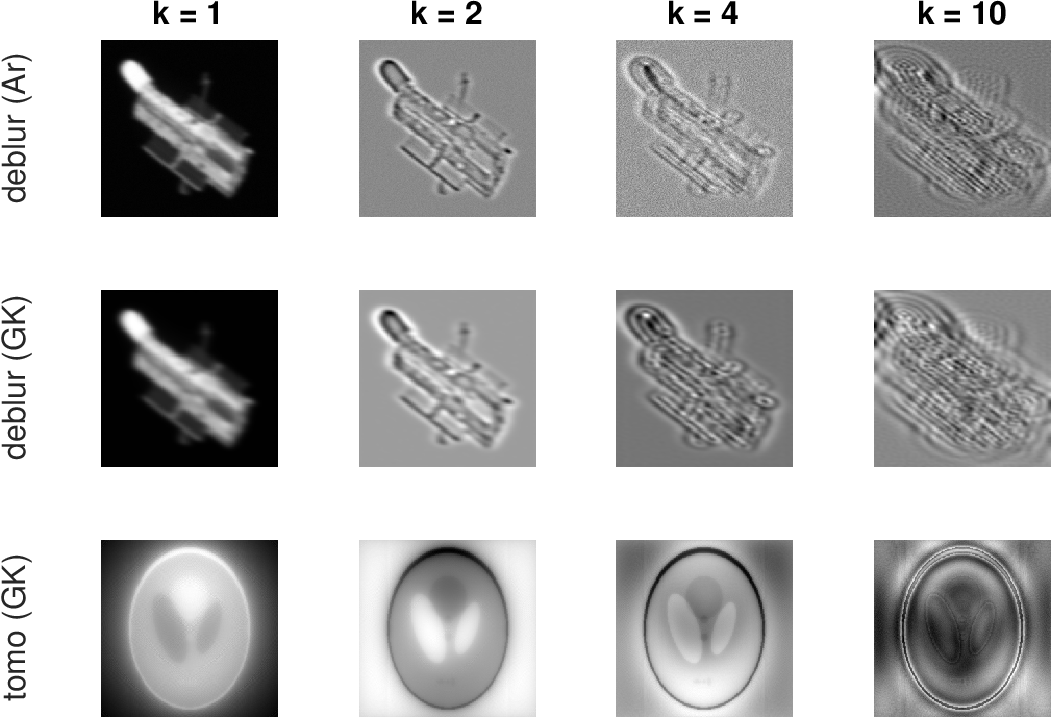}
	\caption{For the deblurring example and the tomography example, we provide basis images (i.e., the basis vectors $\mathbf{v}_k$ for the solution that have been reshaped into images) for iterations $k = 1, 2, 4, 10$. Notice that, since the coefficient matrix in the deblurring example is square, we can compare the Arnoldi and GKB basis vectors, while for the tomography example, we consider only GKB approaches.}
	\label{fig:subspace}
\end{figure}

In \cref{fig:subspace}, we show a few of the basis vectors for the solution generated by the Arnoldi and the GKB algorithms applied to the image deblurring and the tomography reconstruction problems described in \cref{sec:introduction}.  The displayed images are reshaped columns of $\bfV_k$.  For the deblurring example, we observe that, as expected, the first basis vector for the GKB process (i.e., $\bfA\t\bfb/\|\bfA\t\bfb\|$) is smoother or more blurred than that from the Arnoldi process (i.e., $\bfb/\|\bfb\|$).  However, we notice that the $4$th basis vector computed by Arnoldi may have more details but is more noisy than that for GKB.  For the tomography example, we provide some of the GKB solution basis vectors, for which we note that the basis vector for the first iteration is the scaled vector $\bfA\t \bfb$, which can be interpreted as an unfiltered backprojection image; see \cite{Nat01}.

\subsubsection{Solving the regularized, projected problem}
\label{sub:solvingprojected}
In \cref{sub:buildsubspace}, we described two iterative projection methods that can be used to generate and expand the projection subspace ${\rm ran}(\bfV_k)$ (step 2 in \cref{alg:hybridprojection}). Although projection onto subspaces of increasing dimension (e.g., Krylov subspaces) can have an inherent regularizing effect (recall the discussion in \cref{sub:iterative}), a key component of hybrid projection methods is the combination of iterative and variational regularization, i.e., the inclusion of a regularization term within the projected problem. We remark that Steps 3 and 4 in \cref{alg:hybridprojection} are closely intertwined and could easily be addressed together.   However, for clarity of presentation, methods to select the regularization parameter $\lambda_k$
and the stopping iteration $k$ will be addressed in detail in \cref{sec:param}, and this section will focus on solving the projected regularized problem \cref{eq:regprojproblem}. Notice that at step 4 of \cref{alg:hybridprojection}, we have computed $\bfV_k$ and $\lambda_k$, and the stopping criterion is not yet satisfied.

Computing a solution to \cref{eq:regprojproblem} can be done efficiently by exploiting components and relationships from the projection process (\cref{eq:Arnoldi} for Arnoldi and \cref{eq:GKB} for GKB). In this section, we denote the solution subspaces for Arnoldi and Golub-Kahan as $\bfV_k^\text{\tiny Ar}$ and $\bfV_k^\text{\tiny GK}$ respectively.

The \emph{hybrid GMRES iterate} at the $k$th iteration is given by
\begin{align}
	\bfx_k(\lambda_k)=\bfV_k^\text{\tiny Ar} \bfy_k(\lambda_k),\quad\mbox{where}\quad \bfy_k(\lambda_k)
	& = \argmin_{\bfy \in \bbR^k} \norm[]{\bfb - \bfA  \bfV_k^\text{\tiny Ar} \bfy}^2 + \lambda_k \norm[]{\bfy}^2 \nonumber\\
	& = \argmin_{\bfy \in \bbR^k} \norm[]{\bfV_{k+1}^\text{\tiny Ar} \norm[]{\bfb} \bfe_1 - \bfV_{k+1}^\text{\tiny Ar} \bfH_k\bfy}^2 + \lambda_k \norm[]{\bfy}^2 \label{eq:GMRESproj}\\
	& = \argmin_{\bfy\in \bbR^k} \ \norm[]{\norm[]{\bfb} \bfe_1 - \bfH_k \bfy}^2 + \lambda_k \norm[]{\bfy}^2.\nonumber
\end{align}
Note that, since $\norm{\bfV_k^\text{\tiny Ar} \bfy} = \norm[]{\bfy}$ if $(\bfV_k^\text{\tiny Ar})\t \bfV_k^\text{\tiny Ar} = \bfI_k$, we can equivalently write
\[
\bfx_k(\lambda_k) = \argmin_{\bfx \in {\rm ran}(\bfV_k^\text{\tiny Ar})} \norm[]{\bfb - \bfA \bfx}^2 + \lambda_k \norm[]{\bfx}^2
= \bfV_k^\text{\tiny Ar} \bfy_k(\lambda_k)\,.
\]
Oftentimes, the hybrid GMRES approach is referred to as an Arnoldi-Tikhonov approach \cite{calvetti2002gmres,lewis2009arnoldi,gazzola2015survey}.

Similar derivations can be done for the Golub-Kahan projection method.  In particular, the \emph{hybrid LSQR iterate} is computed as
\begin{align}
	\bfx_k(\lambda_k)=\bfV_k^\text{\tiny GK} \bfy_k(\lambda_k),\quad\mbox{where}\quad \bfy_k(\lambda_k)
	& = \argmin_{\bfy \in \bbR^k} \norm[]{\bfb - \bfA  \bfV_k^\text{\tiny GK} \bfy}^2 + \lambda_k \norm[]{\bfy}^2 \nonumber\\
	& = \argmin_{\bfy \in \bbR^k} \norm[]{\beta_1 \bfe_1 - \bfB_k \bfy}^2 + \lambda_k \norm[]{\bfy}^2
	\label{eq:LSQRproj}
\end{align}
or, equivalently,
\[
\bfx_k(\lambda_k) = \argmin_{\bfx \in {\rm ran}(\bfV_k^\text{\tiny GK})} \norm[]{\bfb - \bfA \bfx}^2 + \lambda_k \norm[]{\bfx}^2
= \bfV_k^\text{\tiny GK} \bfy_k(\lambda_k)\,.
\]
The \emph{hybrid LSMR iterate} \cite{chung2015hybrid} is computed as
\begin{align}
	\bfx_k(\lambda_k)=\bfV_k^\text{\tiny GK} \bfy_k(\lambda_k),\quad\mbox{where}\quad \bfy_k(\lambda_k)
	& = \argmin_{\bfy \in \bbR^k} \norm[]{\bfA\t(\bfb - \bfA  \bfV_k^\text{\tiny GK} \bfy)}^2 + \lambda_k \norm[]{\bfy}^2 \nonumber\\
	& = \argmin_{\bfy\in \bbR^k} \norm[]{\bar\beta_1 \bfe_1 -\begin{bmatrix} \bfB_k\t \bfB_k \\ \bar \beta_{k+1} \bfe_k\t \end{bmatrix} \bfy}^2 + \lambda_k \norm[]{\bfy}^2, \label{eq:LSMRproj}
\end{align}
where $\bar \beta_k = \alpha_k \beta_k.$

\cref{tab:methods} provides a summary of common methods with their defining subspace and corresponding subproblem.
\begin{table}[bthp]
\caption{Solution subspaces $\bfV_k$ and corresponding subproblems defining $\bfy_k(\lambda_k)$ for different hybrid projection methods.}
	\label{tab:methods}
	\begin{center}
	\begin{tabular}{|l|l|l|}
		\hline
		\textbf{method} & \textbf{subspace} & \textbf{subproblem} \\ \hline
hybrid GMRES & Arnoldi, ${\rm ran}(\bfV_k^\text{\tiny Ar}) $ & $\min\norm[]{ \beta_1\bfe_1 - \bfH_k \bfy}^2 + \lambda_k \norm[]{\bfy}^2$\\ \hline
		hybrid LSQR & Golub-Kahan, ${\rm ran}(\bfV_k^\text{\tiny GK}) $ & $\min\norm[]{ \beta_1\bfe_1 - \bfB_k \bfy}^2 + \lambda_k \norm[]{\bfy}^2$\\ \hline
		hybrid LSMR & Golub-Kahan, ${\rm ran}(\bfV_k^\text{\tiny GK}) $ & $\min \norm[]{ \bar\beta_1 \bfe_1 - \begin{bmatrix} \bfB_k\t \bfB_k \\ \bar\beta_{k+1} \bfe_k\t \end{bmatrix}
		\bfy}^2 + \lambda_k \norm[]{\bfy}^2$\\ \hline
	\end{tabular}
\end{center}
\end{table}

One important property to highlight is that, for these methods, the generated subspace is independent of the choice of the regularization parameters.  This is not always true (e.g., in the flexible methods presented in \cref{sub:flexible}).  A very desirable consequence of this is that one can avoid computing the regularized solution at each iteration.
For many of the regularization parameter selection methods, the choice of the parameter at each iteration does not depend on previous or later iterates.  Thus, it is possible to delay regularization parameter selection and solution computation until a solution is needed or some stopping criterion is satisfied.

Next, we draw some connections and distinctions to standard iterative methods for least-squares problems (i.e., for the case where $\lambda_k=0$ or $\lambda_k = \lambda$, fixed along the iterations). For problems where no regularization is imposed on the projected problem (i.e., $\lambda_k=0$), we recover standard iterative methods GMRES for Arnoldi (solving \cref{eq:GMRESproj}) and LSQR for GKB (solving \cref{eq:LSQRproj}).

For hybrid projection methods, a potential concern is the need to store all of the solution vectors $\bfV_k$.  For Arnoldi based methods for solving Tikhonov regularized problems, this requirement is the same as for standard iterative methods.  However, for Golub-Kahan based methods, this potential caveat warrants a discussion.  Indeed, it is well-known that the LSQR iterates can be computed efficiently using a three-term-recurrence property by exploiting an efficient QR factorization of $\bfB_k$.  As shown in \cite{PaSa82b}, such computational efficiencies can also be exploited for standard Tikhonov regularization if $\lambda$ is fixed a priori.  This can be done by exploiting the fact that the Krylov subspace is shift invariant with respect to $\lambda.$
However, semi-convergence will be an issue if $\lambda =0$ and the entire process must be restarted from scratch if a different $\lambda$ is desired.  This computational flaw is even more severe when using iterative methods such as LSQR for general-form Tikhonov regularization, since the solution subspace typically depends on the regularization parameter and the regularization matrix (although some strategies to ease this dependence are described in \cref{sub:generalform}). For hybrid projection methods based on GKB, storage of the solution vectors is the main additional cost associated with
the ability to select $\lambda$ adaptively during the hybrid projection procedure. For problems that require many iterations, potential remedies include developing a good preconditioner, compression and/or augmentation techniques \cite{jiang2021hybrid}.

\subsubsection{A unifying framework}\label{ssec:unify}
All the methods described so far can be expressed by this partial factorization,
\begin{equation}\label{eq:genfact}
\bfA\bfV_k = \bfU_{k+1}\bfG_k\,,
\end{equation}
where the columns of $\bfV_k$ are orthonormal and span the $k$-dimensional approximation subspace for the solution, the columns of $\bfU_{k+1}$ are orthonormal with $\bfu_1 = \bfb / \norm[]{\bfb}$ and span a $(k+1)$-dimensional subspace (e.g., associated to $\bfA\t$), and $\bfG_k$ is a $(k+1) \times k$ matrix that has some structure and represents the projected problem.
For hybrid GMRES, $\bfV_k = \bfV_k^\text{\tiny Ar}$, $\bfU_{k+1} = \bfV_{k+1} = \bfV_{k+1}^\text{\tiny Ar}$ and $\bfG_k = \bfH_k$ are given in \eqref{eq:Arnoldi}.  For hybrid LSQR and hybrid LSMR, $\bfV_k = \bfV_k^\text{\tiny GK}$, $\bfG_k = \bfB_k$, and $\bfU_{k+1}$ are given in \eqref{eq:GKB}. Let the SVD of the matrix $\bfG_k$ be given by
\begin{equation}
  \label{eq:SVD_G}
  \bfG_k = \bfU^\text{\tiny G} \bfSigma^\text{\tiny G} (\bfV^\text{\tiny G})\t,
\end{equation}
where $\bfU^\text{\tiny G} \in \bbR^{k+1 \times k+1}$ and $\bfV^\text{\tiny G}\in \bbR^{k \times k}$ are orthogonal, and $\bfSigma^\text{\tiny G} = \diag{\sigma_1^\text{\tiny G}, \ldots, \sigma_k^\text{\tiny G}} \in \bbR^{(k+1) \times k}$ contains the singular values of $\bfG_k$. Note that, for notational convenience, we have dropped the subscript $k$.

For hybrid GMRES and hybrid LSQR with standard Tikhonov regularization, the solution to the regularized projected problems \eqref{eq:GMRESproj} and \eqref{eq:LSQRproj} have the form,
\begin{align}
	\bfy_k(\lambda_k) & = \argmin_\bfy \norm[]{\bfb - \bfA \bfV_k \bfy}^2 + \lambda_k \norm[]{\bfy}^2
	\nonumber\\
	& = \argmin_\bfy \norm[]{\bfU_{k+1}\t\bfb - \bfG_k \bfy}^2 + \lambda_k \norm[]{\bfy}^2, \label{eq:proj_unified2}
\end{align}
and thus one may express the $k$th iterate of the hybrid projection method as
\begin{equation}
	\label{eq:hybridsolution}
	\bfx_k(\lambda_k) = \underbrace{\bfV_k(\bfG_k\t\bfG_k+\lambda_k\bfI_k)^{-1}\bfG_k\t\bfU_{k+1}\t}_{=:\bfAregpinv(\lambda_k,k)} \bfb.
\end{equation}
Note that $\bfU_{k+1}\t\bfb=\|\bfb\|\bfe_1$ by construction.

For many of the regularization parameter selection methods described in \cref{sec:param} and the theoretical results in \cref{sec:theory}, it will be helpful to define the so-called `influence matrix',
\begin{equation}\label{eq:influenceM}
\bfA\bfAregpinv(\lambda_k,k)=
\bfU_{k+1}\bfG_k(\bfG_k\t\bfG_k+\lambda_k\bfI_k)^{-1}\bfG_k\t\bfU_{k+1}\t\,.
\end{equation}
For purely iterative methods (where $\lambda_k = 0$), the influence matrix is given by
$$\bfA\bfAregpinv(k)=
\bfU_{k+1}\bfG_k^\dagger\bfU_{k+1}\t\,, \quad \mbox{where} \quad \bfAregpinv(k)=\bfV_k\bfG_k^\dagger\bfU_{k+1}\t.
$$
Also, it will be helpful to define the so-called `discrepancy' (or `regularized residual') at the $k$th iteration as
\begin{equation}\label{eq:discrep}
\bfr(\bfx_k(\lambda_k)) = \bfb - \bfA\bfx_k(\lambda_k) = (\bfI_m - \bfA\bfAregpinv(\lambda_k,k))\bfb\,,
\end{equation}
for hybrid methods; the same definition applies to iterative methods, with $\lambda_k = 0$. In the following, we will adopt the notation $\bfAregpinv(\lambda)$ and $\bfr(\bfx(\lambda))$ to denote the regularized inverse and the discrepancy associated to (direct) Tikhonov regularization, coherently to \cref{eq:hybridsolution} and \cref{eq:discrep}, respectively. When no regularization method is specified, $\bfxreg$, $\bfAregpinv$ and $\bfrreg$ denote a generic regularized solution, inverse, and residual, respectively.

\subsection{Regularization parameter selection methods}
\label{sec:param}
The success of any regularization method depends on the choice of one (or more) regularization parameter(s). This was illustrated in \cref{sec:background} for Tikhonov regularization and for iterative regularization methods; also, when the Tikhonov-regularized problem is solved using an iterative method, both the Tikhonov regularization parameter and the number of iterations should be accurately tuned (see \cref{sub:IterTikh}). Similarly, for hybrid projection methods, there are inherently two  regularization parameters to tune: (1) the number of iterations $k$ (i.e., the dimension of the projection subspace), and (2) the regularization parameter for the projected problem (e.g., $\lambda_k$ for Tikhonov). It is important to note that, when carelessly applied to hybrid projection methods where \emph{both} regularization parameters must be determined, standard regularization parameter choice strategies seldom produce good results. Instead, a two-pronged approach is typically adopted, where a well-established parameter choice strategy is used to determine $\lambda_k$ and one or more stopping rules are used to terminate the iterations, typically by monitoring the stabilization of some relevant quantities; see, for instance, \cite{ChNaOL08, GaNo14, ReVaAr17, CaGoRe99,KiOLe01}. Before describing specific parameter selection strategies, we comment on a few factors that are unique to selecting regularization parameters in a hybrid projection method, and that should be considered when determining which parameter selection methods to employ.

\paragraph{Selecting $\lambda_k$} Assume that $k$ is fixed and consider the case of standard form Tikhonov regularization.  First, for many Krylov subspace projection methods, the projection subspace is independent of the current value of the Tikhonov regularization parameter.  This is due to the shift-invariance property of the computed Krylov subspaces, and this has the effect that the approximate solution at the $k$th iteration only depends on the current regularization parameter $\lambda_k$.  Second, notice that, for a fixed $k$, the projected problem \eqref{eq:proj_unified2} is a standard Tikhonov-regularized problem.  A natural idea would be to directly apply well-established (e.g, SVD-based) parameter choice strategies that were developed for Tikhonov regularization; however, these methods are not self-contained in the hybrid setting, and care must be taken to ensure accurate results.
Third, since the number of iterations $k$ is significantly smaller than the size of the original problem (i.e., $k\ll \min\{m,n\}$) and the projected problem has order $k$ (i.e., the coefficient matrix $\bfG_k$ has size $(k+1) \times k$), computations with $\bfG_k$ can be performed efficiently. For many methods, there are computational advantages to using the SVD of $\bfG_k$ \cref{eq:SVD_G}. Indeed, the computational cost of obtaining \eqref{eq:SVD_G} is negligible compared to the computational cost of performing matrix-vector products with $\bfA$ (and possibly $\bfA\t$) to expand the approximation subspace.  We will see how some common parameter choice methods can be efficiently formulated using the SVD.

\paragraph{Selecting $k$}
Various rules have been developed for selecting the stopping iteration $k$, and these rules can be employed quite generally (and independently of the parameter choice rule for $\lambda_k$).  The main idea is to terminate iterations when a maximum number of iterations is achieved or when one or more of the following conditions are satisfied:
\begin{equation}\label{SCstab}
\begin{array}{lcl}
{|\lambda_k - \lambda_{k-1}|}&<&\tau_{\lambda}{\lambda_{k-1}}\,,\\
{\|\bfr(\bfx_k(\lambda_k))-\bfr(\bfx_{k-1}(\lambda_{k-1}))\|}&<&\tau_{r}{\|\bfr(\bfx_{k-1}(\lambda_{k-1}))\|}\,,\\
{\|\bfx_k(\lambda_k)- \bfx_{k-1}(\lambda_{k-1})\|}&<&\tau_{x}{\|\bfx_{k-1}(\lambda_{k-1})\|}\,,
\end{array}
\end{equation}
where $\tau_{\lambda}$, $\tau_{r}$, $\tau_{x}$ are positive user-specified thresholds.
The rationale behind these approaches is that, often and broadly speaking, when stabilization happens, the selected value for the regularization parameter is suitable for the full-dimensional Tikhonov problem and the approximated solution cannot significantly improve with more iterations. Although this argument is mainly empirical, in some cases it is supported by theoretical results; see \cref{sub:quadRule}. This property is also heavily exploited in  \cite{GoVMa91,ChKiOL15,renaut2017efficient}, where the regularization parameter (and other relevant quantities) for the original Tikhonov problem are estimated by projecting the original problem onto subspaces of smaller dimension.  Once the regularization parameter is estimated, it is fixed and any iterative method can be used to solve the resulting Tikhonov problem; therefore, these methods are not considered hybrid projection methods as defined in this manuscript; see \cref{sec:introduction} and \cref{sub:IterTikh}. The left frame of \cref{fig:param_adaptive} displays the typical behavior of adaptively chosen regularization parameters where, relevant to the first criterion in \cref{SCstab}, some stabilization is visible as the iterations  proceed.
Returning to the issue of selecting a stopping iteration for hybrid projection methods, we remark that, with a suitable choice of $\lambda_k$, hybrid methods can overcome semi-convergent behavior, as illustrated in \cref{fig:semiconvergence}.  Thus, an imprecise (over-)estimate of the stopping iteration does not significantly degrade the reconstruction quality.  In fact, one can typically afford a few more iterations without experiencing deterioration of the solution (on the contrary, the solution may improve because it is computed by solving a well-posed problem in a larger approximation subspace). This insight is also linked to the ability of the considered Krylov projection methods to `capture' the dominant (i.e., relevant) truncated right singular vector subspace information to reconstruct the solution; we present more details about this in \cref{sec:theory}.

Parameter choice rules represent a large and growing body of literature in the field of inverse problems, with papers ranging from theoretical developments of regularization methods to papers focused on methods specific to applications.  Although there are extensive survey papers describing parameter choice methods in the continuous setting (for instance, see \cite{Bauer1,OldNewP,Bauer2015}), we focus on parameter selection strategies for discrete inverse problems that have proven successful in conjunction with hybrid projection methods. We describe two main classes of parameter choice methods: (1) those that require knowledge of the noise magnitude in \cref{sub:paramNoise} and (2) those that do not require knowledge of the noise magnitude in \cref{sub:paramNoNoise}.
In almost all the considered strategies for hybrid methods, there are common quantities that must be monitored.  These include:
\begin{enumerate}
  \item the norm of the approximate solution, $\|\bfx_k(\lambda_k)\|$,
  \item the norm of the residual, $\|\bfr(\bfx_k(\lambda_k))\|$, and
  \item the trace of the influence matrix, $\trace{\bfA \bfAregpinv(\lambda_k,k)}$.
\end{enumerate}
Computing these quantities can be done very efficiently by monitoring the corresponding projected quantities and exploiting the orthogonal invariance of the 2-norm. This fact along with the relatively cheap computation of \eqref{eq:SVD_G} make the use of standard parameter choice strategies particularly appealing in the setting of hybrid methods. Simplified formulations for these values using the SVD can be found in \cref{sec:SVDformulations}.

\subsubsection{Methods that require knowledge of the noise magnitude}\label{sub:paramNoise}
The \emph{discrepancy principle} prescribes to select a regularized solution $\bfxreg$ satisfying
\begin{equation}\label{eq:discrP}
\|\bfb-\bfA\bfxreg\|=\eta\eps\,,
\end{equation}
where $\eps$ is an estimate of the norm of the noise $\|\bfe\|$ and $\eta>1$ is a safety factor (the larger $\eta$, the more uncertainty on $\eps$). The discrepancy principle typically works well if a good estimate of $\eps$ is available. Of course, the discrepancy principle also works well when $\bfxreg$ only depends on one regularization parameter: indeed, when considering (direct) Tikhonov regularization (\cref{sub:direct}), $\bfxreg = \bfx(\lambda)$, and \cref{eq:discrP} is a nonlinear equation in $\lambda>0$ (that should be solved by employing a zero-finder, typically of the Newton family \cite{ReichelAnzf}); when considering iterative regularization methods (\cref{sub:iterative}), $\bfxreg = \bfx_k$, and the iterative method should stop as soon as the $k$th iterate satisfies
\begin{equation}\label{DPstop}
\|\bfb-\bfA\bfx_k\|\leq\eta\eps\,.
\end{equation}
 Using the discrepancy principle, one can typically prove regularization properties of the kind $\bfxreg\rightarrow\bfx_\true$ as
 $\norm[]{\bfe}\rightarrow 0$; see \cite{engl1996regularization,hanke1995minimal,GMRESreg}.
 When considering hybrid solvers, the most common approach is to (fully) solve the nonlinear equation
\begin{equation}\label{eq:discrPhybr}
\|\bfb-\bfA\bfx_k(\lambda_k)\|=\eta\eps
\end{equation}
to determine the Tikhonov regularization parameter $\lambda_k$ to be employed at the $k$th iteration. If adopting this strategy, the discrepancy principle is satisfied at each iteration, and one cannot exploit \cref{DPstop} to also set a stopping criterion. In this case, one must resort to one (or more) stopping criterion of the form \cref{SCstab}. To be precise, and specifically for the discrepancy principle, a stopping criterion naturally arises, in that the equation \cref{eq:discrPhybr}
can be satisfied with respect to $\lambda_k$ typically only after a few iterations $k$ are performed: this property alone can act as a stopping criterion (i.e., one may stop as soon as the discrepancy principle \cref{eq:discrPhybr} can be satisfied). However, typically, the quality of the solutions improves if more iterations are performed and $\lambda_k$ is computed by solving \cref{eq:discrPhybr}, as the regularized solutions belong to richer approximation subspaces. This phenomenon is described at length in \cite{lewis2009arnoldi,reichel2012tikhonov}. Strategies that use the discrepancy principle for setting both $k$ and $\lambda_k$ have also been devised, see \cite{GaNo14,gazzola2021}. In other words, these approaches  update the regularization parameter for the projected problem in such a way that stopping by the discrepancy principle is ensured: this is typically achieved by performing only one iteration of a root-finder algorithm for \cref{eq:discrPhybr} at each iteration of a hybrid method; the approach derived in \cite{gazzola2021} and its underlying theory is described in more details in \cref{sub:quadRule}.

\emph{UPRE} (\emph{unbiased predictive risk estimation}) prescribes to choose the regularization parameter that minimizes the expectation of the predictive error, $\bbE(\bfp(\bfxreg))$, associated to the regularized solution $\bfxreg=\bfAregpinv\bfb$, where the predictive error is defined as
\[
\bfp(\bfxreg) = \bfA\bfxreg - \bfb_\true=\bfA\bfAregpinv\bfb - \bfb_\true = (\bfA\bfAregpinv - \bfI_n)\bfb_\true + \bfA\bfAregpinv\bfe\,,
\]
where $\bfA\bfAregpinv$ is assumed to be symmetric (this is the case for all the regularization methods considered so far). When considering (direct) Tikhonov regularization and when $\bfe$ represents Gaussian white noise with standard deviation $\sigma$, i.e., when
\begin{equation}\label{GaussWhite}
\bfe\sim\calN(\bfzero,\sigma^2\bfI_m)\,,
\end{equation}
one should compute
\[
\lambda^\ast \!\!=\!\argmin_{\lambda\in\bbR_+} \underbrace{\|\ (\bfA\bfAregpinv(\lambda) - \bfI_m)\bfb_\true \|^2 + \sigma^2\trace{(\bfA\bfAregpinv(\lambda))\t(\bfA\bfAregpinv(\lambda))}}_{=:\,\bbE(\|\bfp(\bfx(\lambda))\|^2)}\,.
\]
To circumvent the fact that the first term in the above expression of $\bfp(\bfx(\lambda))$ is unavailable, one should perform some algebraic manipulations and approximations
to get
\begin{equation}\label{UPREtwo}
\lambda^\ast = \arg \min_{\lambda\in\bbR_+} \underbrace{\| \bfr(\bfx(\lambda)) \|^2 + 2\sigma^2\trace{\bfA\bfAregpinv(\lambda)} - m\sigma^2}_{=:\,U(\lambda)}\,.
\end{equation}
Note that, in particular, $\bbE(U(\lambda))=\bbE(\|\bfp(\bfx(\lambda))\|^2)$. We refer to \cite{Vog02} for more details on the derivation of the UPRE method.
UPRE can also be used as a stopping rule for iterative methods, e.g., for nonnegativity constrained Poisson inverse problems \cite{Bar08a}.
When considering hybrid methods, applying UPRE to the projected problem is quite straightforward, and was first considered in \cite{ReVaAr17}: it is essentially a matter of replacing $\bfr(\bfx(\lambda))$ and $\bfAregpinv(\lambda)$ in \cref{UPREtwo} by $\bfr(\bfx_k(\lambda_k))$ and $\bfAregpinv(\lambda_k,k)$, respectively.
Since the influence matrix \cref{eq:influenceM} is still symmetric, performing some algebraic manipulations leads to the following projected UPRE functional
\begin{eqnarray}
{U}_k(\lambda)=&&\| (\bfG_k(\bfG_k\t\bfG_k+\lambda\bfI_k)^{-1}\bfG_k\t - \bfI_{k+1})\bfU_{k+1}\t\bfb\|^2
\nonumber\\
&+& 2\sigma^2\trace{\bfG_k(\bfG_k\t\bfG_k+\lambda\bfI_k)^{-1}\bfG_k\t} - (k+1)\sigma^2\,,\nonumber
\end{eqnarray}
which is minimized at each iteration of a hybrid method to find $\lambda_k$.
An interesting (and still partially open) question (common to other parameter choice strategies) is determining wether the regularization parameter $\lambda_k$ so obtained
is a good approximation of $\lambda$.
To answer this question, the authors of \cite{ReVaAr17} first consider the direct regularization method to be obtained by combining TSVD and Tikhonov regularization methods (sometimes referred to as `FTSVD', i.e., filtered TSVD), so that the variational regularized solution belongs to the subspace spanned by the dominant right singular vectors; they deduce that, if the $k$-dimensional projection subspace generated by the hybrid method captures the relevant (i.e., dominant) spectral information about the original problem, then
$\trace{\bfG_k(\bfG_k\t\bfG_k+\lambda\bfI_k)^{-1}\bfG_k\t}\approx \trace{\bfA\bfAregpinv}$ (the latter being specified for the FTSVD), and $\lambda_k\approx \lambda^\ast$.

We conclude this subsection by mentioning that, although it may seem to be a disadvantage that these methods require knowledge of the noise magnitude, there are actually various approaches for estimating the noise level from the data: here we describe a couple of them. A first approach uses statistical tools and performs quite well at estimating the variance \cite{donoho1995noising}. Assuming white noise \cref{GaussWhite}, an estimate $\hat{\sigma}$ of the standard deviation $\sigma$ can be obtained from the highest coefficients of the noisy data under some transformation (e.g., wavelet).
For instance, the following MATLAB code can be used
\begin{equation}\label{noiseestwav}
\begin{array}{lcl}
\text{\texttt{[\~{}, cD]}} & \text{\texttt{=}} & \text{\texttt{dwt(b,'db1');}}\\
\text{\texttt{sigmahat}} & \text{\texttt{=}} & \text{\texttt{median(abs(cD(:)))/.67;}}
\end{array}.
\end{equation}
In inverse problems, we often scale the noise, and use the following,
$$\bfb = \bfb_{\rm true} + \texttt{noiseLevel}* \frac{\norm[]{\bfb_{\rm true}}}{\norm[]{\bfe}} \bfe\,,$$
where $\bfe$ is a realization from a standard Gaussian and the noise level is given in \linebreak[4]\texttt{noiseLevel}. Thus, an estimate of \texttt{noiseLevel} obtained running \cref{noiseestwav} would be $\hat{\sigma} \sqrt{m}/\norm[]{\bfb}$.
Still assuming Gaussian white noise, a second approach to estimate relevant noise information (including $\|\bfe\|$) is to leverage some theoretical properties of the Krylov subspaces generated as projection subspaces for the solution: this is especially relevant in the context of hybrid regularization. To the best of our knowledge, the first detailed analysis of how the noise affects the approximation subspace generated by the GKB algorithm can be found in \cite{HnPlSt09}, where the authors show how to estimate at a negligible cost the (assumed) unknown amount of noise in the original data. By exploiting the connections between GKB and Gaussian quadrature rules (recalled in \cref{sub:quadRule}), theoretical estimates prove that the norm of the residual computed by LSQR stabilizes around the noise magnitude. This can lead to the construction of stopping criteria for the bidiagonalization process as well as to the application of any of the parameter choice rules described in this section, using the information gathered about the noise. Extensions of the analysis in \cite{HnPlSt09} to projection methods based on the Arnoldi algorithm can be found in \cite{GaNoRu14}.

\subsubsection{Methods that do not require knowledge of the noise magnitude}
\label{sub:paramNoNoise}
In many situations, assuming that an accurate estimate of $\eps=\|\bfe\|$ is available is unrealistic, and one cannot confidently apply the methods described in \cref{sub:paramNoise}. However, a number of strategies can be adopted when dealing with direct, iterative and hybrid methods: for the latter, most of these parameter choice rules can be regarded as the projected variants of their full-dimensional counterparts.

The \emph{$L$-curve} criterion was popularized by \cite{HOL93}: it relies on the intuition that a good choice of the regularization parameter should balance the contribution of the so-called perturbation error (i.e., the error due to overfitting the noisy data, which would lead to an under-regularized solution) and regularization error (i.e., the error due to replacing the original problem with a related one, which would lead to an over-regularized solution). The $L$-curve is a plot of the norm of the regularized solution $\|\bfxreg\|$ versus the norm of the regularized residual $\|\bfrreg\|$ for varying values of the regularization parameter, and it is named after the desirable shape of its graph.
When Tikhonov regularization is considered (i.e., when the $L$-curve is a parametric curve with respect to $\lambda$), one can prove that the $L$-curve is convex (see, for instance, \cite[Chapter 4]{hansen2010discrete}). The ideally steep or vertical part of the curve corresponds to small amounts of regularization, so that such solutions are dominated by perturbation error.  The ideally flat or horizontal part corresponds to too much regularization, so that such solutions are dominated by regularization error.  Therefore the corner represents the point on the $L$-curve where both errors are balanced. The $L$-curve is commonly plotted in logarithmic scale, e.g., $(\log_{10}(\|\bfrreg\|),\log_{10}(\|\bfxreg\|))$, to better highlight the corner (and also because of the large range of values of the plotted quantities). For the $L$-curve to be effective it is necessary to have monotonicity in $\|\bfrreg\|$ and $\|\bfxreg\|$ (this is not always the case, e.g., alternative $L$-curves have been devised for GMRES \cite{calvetti2002gmres}). Since for hybrid projection methods we need to select two parameters (namely, $\lambda_k$ and $k$), the $L$-curve can be interpreted as a surface, which can be challenging to analyze.  However, by exploiting connections to Gaussian quadrature rules (see \cref{sub:quadRule}), a variant of the $L$-curve called the `$L$-ribbon' has been considered that inexpensively constructs a ribbon-like region that contains the $L$-curve of the (direct) Tikhonov regularization applied to the full-dimensional problem \cite{calvetti2002curve, Calvetti2004,CaMoReSg00}. 

An approach related to the $L$-curve, which still aims at finding the right balance between the regularization error and the perturbation error, was described in \cite{bazan2008fixed,viloche2012generalization}, where the estimated regularization parameter $\lambda$ for (direct) Tikhonov regularization is obtained as the fixed point of the ratio between the residual norm $\|\bfr(\bfx(\lambda))\|$ and the solution norm $\|\bfx(\lambda)\|$.  This same strategy can be straightforwardly used to select a regularization parameter at each iteration of hybrid projection methods, namely, by computing the fixed point $\lambda_k$ of the function $\frac{\norm[]{\bfr(\bfx_k(\lambda))}}{\norm[]{\bfx_k(\lambda)}}$: this is done in \cite{bazan2010gkb} for a hybrid method based on GKB and Tikhonov regularization. These approaches can be regarded as generalizations of a parameter choice rule due to Regi{\'n}ska \cite{reginska1996regularization}.

The generalized cross-validation (GCV) method is another popular approach for selecting regularization parameters when the noise level is unknown.  The GCV method is a `leave-one-out' prediction method. That is, the basic idea behind GCV is that, if an arbitrary element of the observed data is left out, a good choice of the regularization parameters should be able to predict the missing observation \cite{GoHeWh79}.  For (direct) Tikhonov regularization applied to the full-dimensional  problem \cref{eq:Tikhonov}, the parameter computed by GCV is the one that minimizes the GCV function
\begin{equation}
  \label{eq:GCV_original}
  G(\lambda) = \frac{n \norm[]{\bfr(\bfx(\lambda))}^2}{\left(\trace{\bfI_m - \bfA\bfAregpinv(\lambda)}\right)^2}.
\end{equation}
In \cite{Bjo88} Bj{\" o}rck suggested using GCV in conjunction with TSVD for hybrid projection methods and using GCV to determine an appropriate stopping iteration as well.  However,
it was observed in \cite{ChNaOL08} that the GCV method tended to perform poorly when used within hybrid methods based on GKB and Tikhonov regularization, due to over-smoothing.
To remedy this, a weighted GCV (wGCV) method was introduced, which can be interpreted as a weighted `leave-one-out' approach; an adaptive approach to estimate the new weight parameter was also described.
At each iteration $k$, the weighted GCV functional for estimating $\lambda$ is given by
\begin{equation}
  \label{eq:WGCV}
  G_w(\lambda,k) = \frac{n \norm[]{\bfr(\bfx_k(\lambda))}^2}{\left(\trace{\bfI_m - w \bfA\bfAregpinv(\lambda,k)}\right)^2},
\end{equation}
where we get the standard GCV function for $w=1.$ By minimizing the above functional with respect to both $\lambda$ and $k$, it is possible to set both regularization parameters involved in a hybrid method. One can use an alternating approach to sample the 2D GCV surface \cref{eq:WGCV} as follows:  first, a projected version of the GCV functional \cref{eq:WGCV} is minimized for fixed $k$ to get $\lambda_k$ (this is sampling the GCV functional along lines); second, the GCV functional \cref{eq:WGCV} is minimized for fixed $\lambda_k$.

The noise cumulative periodogram (NCP) approach for selecting regularization parameters does not require an estimate of the noise level either. Indeed, NCP uses the residual components rather than the residual norm to estimate the regularization parameter
\cite{rust2008residual}.  NCP was used within a projection framework in \cite{hansen2006exploiting}, although this was not a hybrid projection method according to the criteria given in the present paper.

\begin{center}
\begin{figure}
\begin{tabular}{cc}
\small{\bf deblurring} & \small{\bf tomography}\\
\includegraphics[width=6cm]{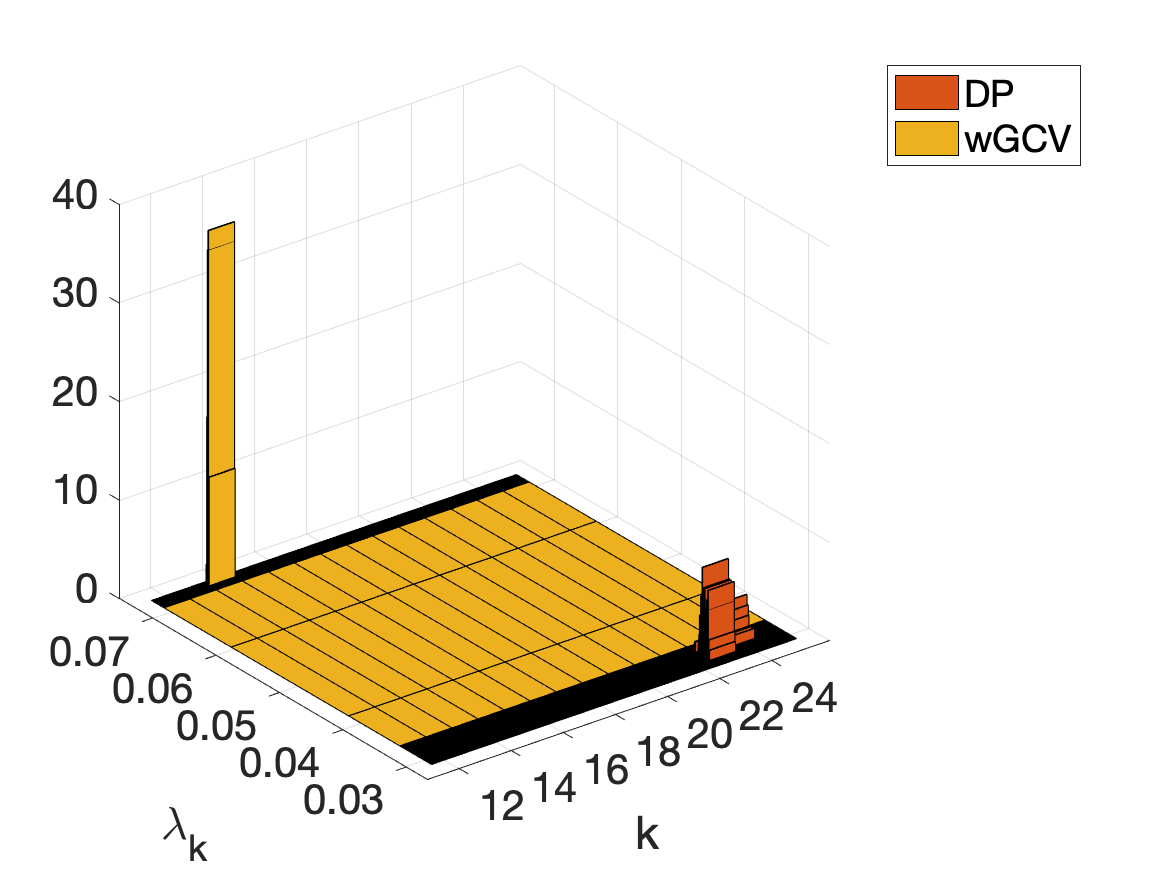} &
\includegraphics[width=6cm]{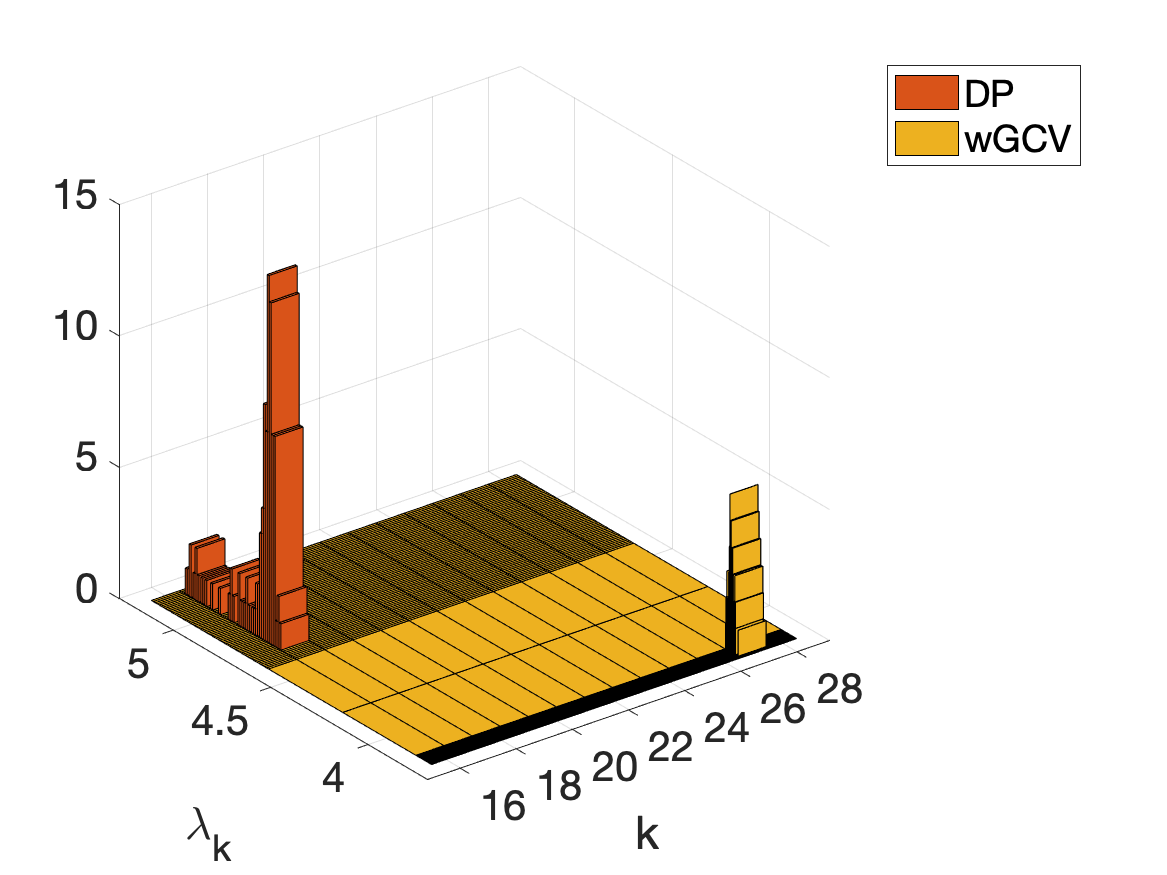}\\
\includegraphics[width=6cm]{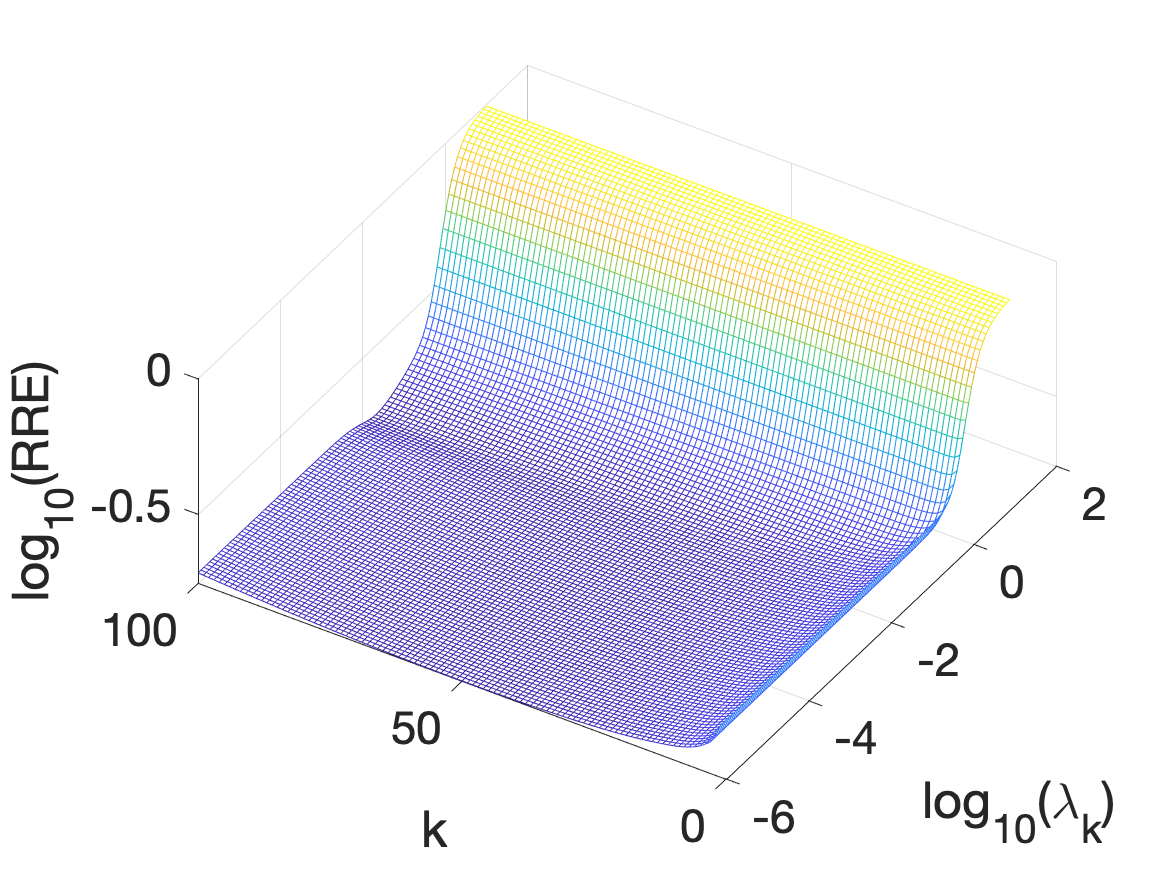} &
\includegraphics[width=6cm]{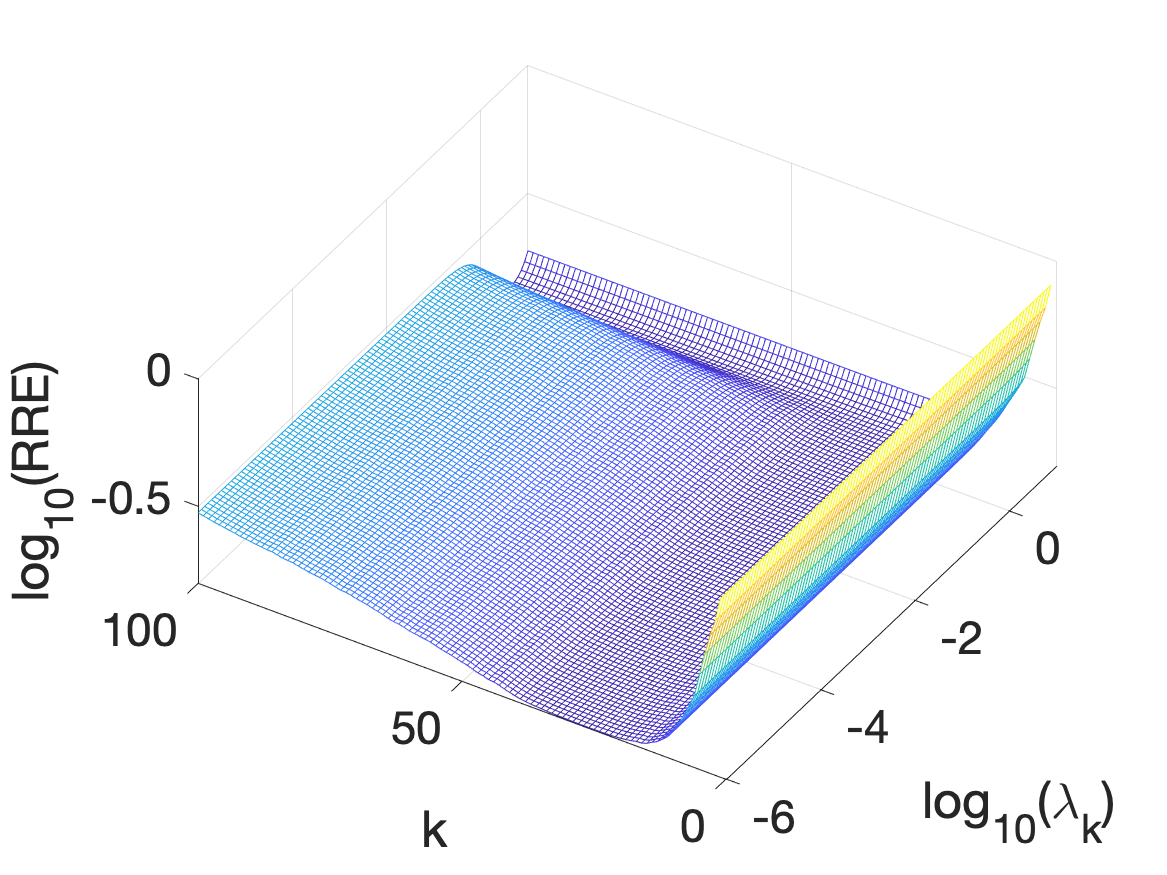}\\
\includegraphics[width=6cm]{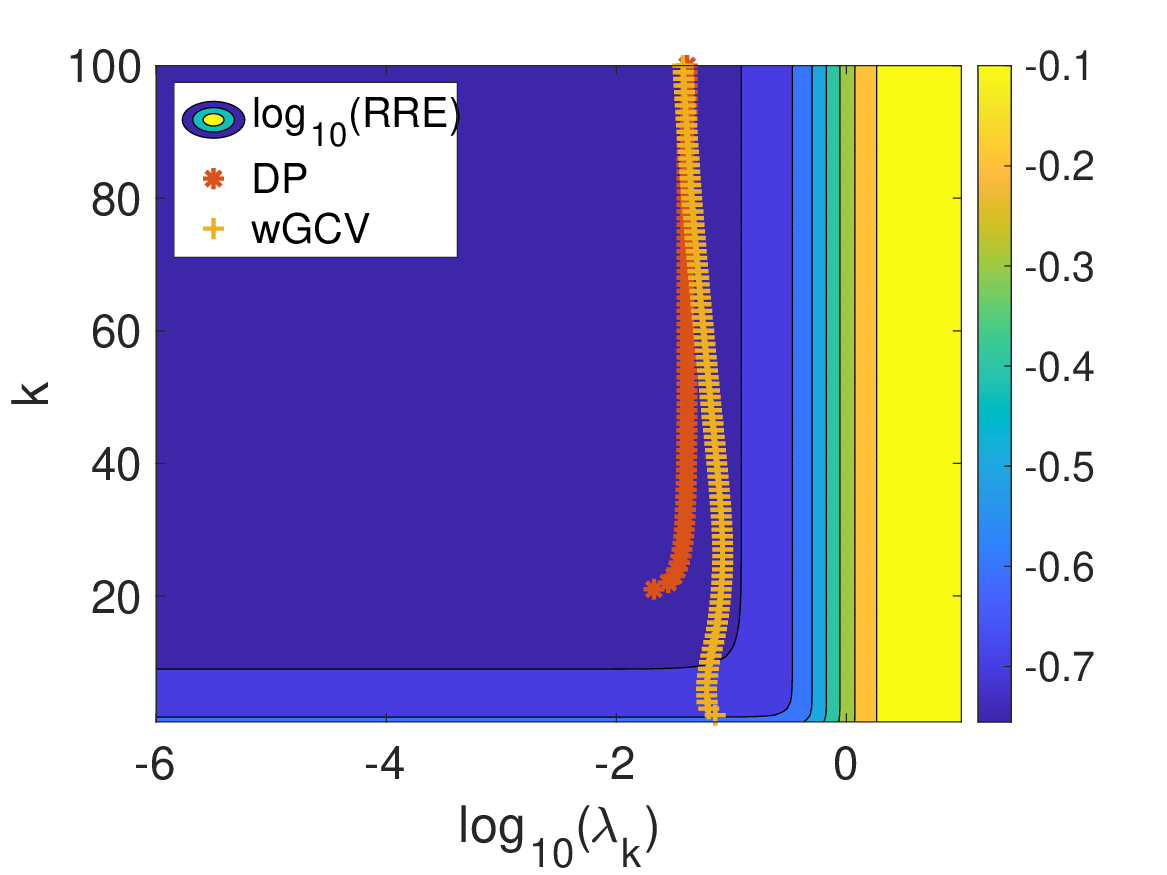} &
\includegraphics[width=6cm]{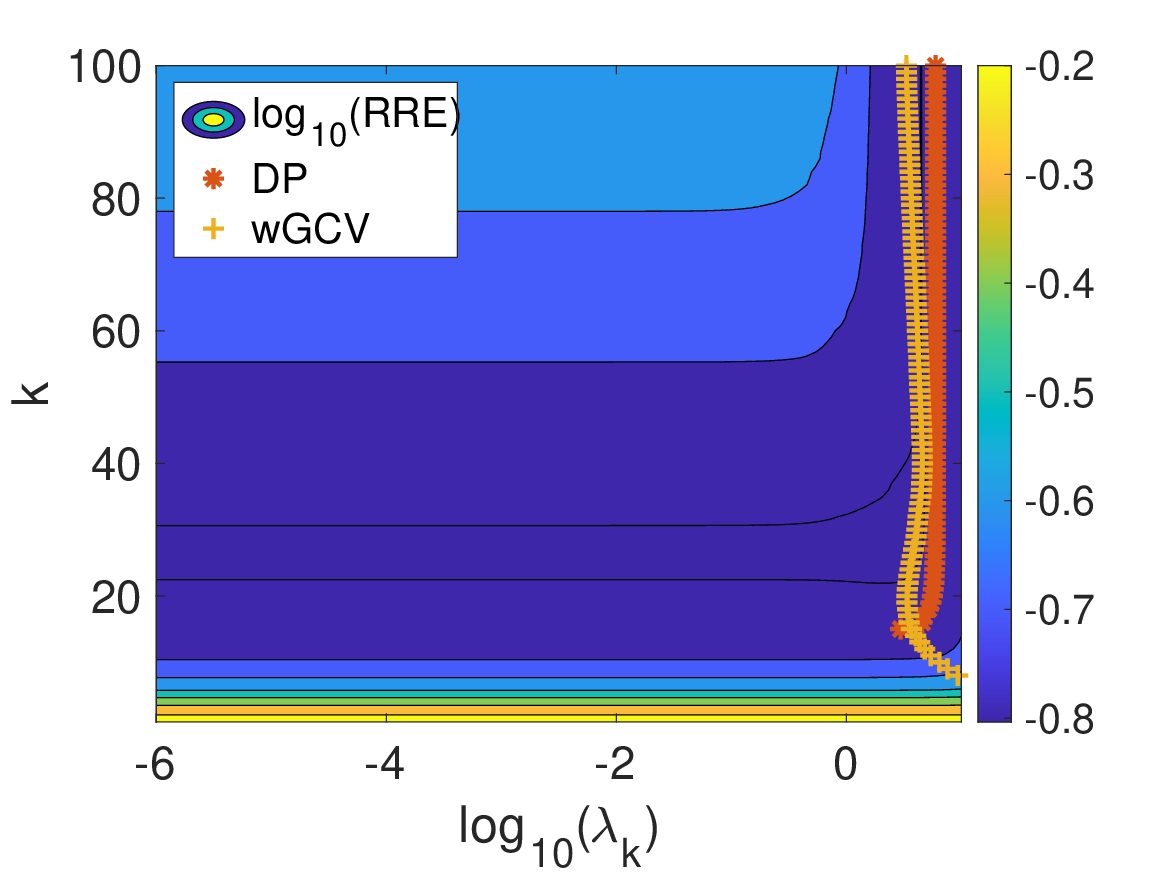}\\
\end{tabular}
\caption{A comparative study of the behavior of the GKB-based hybrid methods when the regularization parameter is chosen according to the discrepancy principle and the wGCV criterion. First row: the test problems are run 100 times, each with a different noise realization (noise level is always $10^{-2}$) and histogram plots of the stopping iteration $k$ and the corresponding chosen regularization parameter $\lambda_k$ are displayed. Second and third row: the relative error norm for each sampled $(k,\lambda_k)$ is recorded, to form the so-called `\emph{RRE surfaces}', which are displayed as 3D surface plots (second row) and as 2D contour plots, where the computed couples $(k,\lambda_k)$ obtained by applying the discrepancy principle and wGCV are highlighted by special markers.}\label{fig:regP}
\end{figure}
\end{center}

\paragraph{Illustration}
We use the image deblurring and tomography test problems from \cref{sec:background} to show how different regularization parameter selection methods can perform within hybrid projection methods. Relevant quantities are displayed in \cref{fig:regP}. For both test problems, looking at the frames displayed in the top row, we can infer that both the discrepancy principle and the wGCV method exhibit consistent behavior within different noise realizations; looking at the frames displayed in the middle and bottom rows, we can see that there are clear combinations of values of $k$ and $\lambda_k$ that deliver minimal relative reconstructions error norms, and that both the discrepancy principle and wGCV are able to compute couples $(k,\lambda_k)$ that eventually (for $k$ big enough) lay in such regions. We emphasize that every problem is different, and there is not one approach that will work for all problems. Thus, it is good to have a variety of methods that can guide one in selecting a suitable set of parameters.

\subsubsection{SVD formulations of quantities used for regularization parameter selection}
\label{sec:SVDformulations}
Notice that in many of the regularization parameter selection strategies, a nonlinear optimization or root finding problem must be solved to estimate $\lambda_k$.
For smaller problems (e.g., the projected problem in hybrid projection methods), the SVD can be exploited for efficient computation of the norm of the approximation solution, the norm of the residual, and the trace of the influence matrix, which all depend on $\lambda_k$.
We remark that these efficiencies are especially impactful when quantities need to be computed for many regularization parameters (e.g., to obtain a representative $L$-curve) or when the projected problem is significantly smaller than the size of the original problem (e.g., $k \ll n$).
In this section, we provide SVD formulations for these quantities for Tikhonov regularization, but remark that similar derivations can be made for any spectral filtering method (e.g., TSVD) used to solve the projected problem in a hybrid projection method.

Let the SVD of $\bfG_k\in \bbR^{(k+1) \times k}$ be given as in \cref{eq:SVD_G} and let
$$\widetilde\bfU = \bfU_{k+1} \bfU^\text{\tiny G} = \begin{bmatrix}\widetilde \bfu_1 & \cdots & \widetilde \bfu_{k+1} \end{bmatrix} \quad \mbox{and} \quad \widetilde\bfV = \bfV_k \bfV^\text{\tiny G}  = \begin{bmatrix}\widetilde \bfv_1 & \cdots & \widetilde \bfv_{k} \end{bmatrix}\,;$$  notice that these contain orthonormal columns. Then the solution of the regularized, projected problem \cref{eq:hybridsolution} has the form
\begin{align}
\bfx_k(\lambda_k)
& = \widetilde \bfV ((\bfSigma^\text{\tiny G})\t\bfSigma^\text{\tiny G}+\lambda_k\bfI_k)^{-1} (\bfSigma^\text{\tiny G})\t  \widetilde\bfU\t \bfb\nonumber\\
& = \sum_{i=1}^k \phi_i^\text{\tiny G} \frac{\widetilde\bfu_i\t \bfb}{\sigma_i^\text{\tiny G}} \widetilde\bfv_i, \quad \mbox{where} \quad \phi_i^\text{\tiny G}  = \frac{(\sigma_i^\text{\tiny G})^2}{(\sigma_i^\text{\tiny G})^2+\lambda_k}.\nonumber
\end{align}
The `influence matrix' \cref{eq:influenceM} can then be represented as
\begin{align}
\bfA\bfAregpinv (\lambda_k,k)
& = \widetilde \bfU \bfSigma^\text{\tiny G} ((\bfSigma^\text{\tiny G})\t\bfSigma^\text{\tiny G}+\lambda_k\bfI_k)^{-1} (\bfSigma^\text{\tiny G})\t  \widetilde\bfU\t.\nonumber
\end{align}

Using the above formulas, the squared norm of the solution can be represented as
  \[
    \|\bfx(\lambda_k)\|^2 = \sum_{i=1}^k \left(\phi_i^\text{\tiny G} \frac{\widetilde\bfu_i\t \bfb}{\sigma_i^\text{\tiny G}}\right)^2
  \]
and the trace of the influence matrix takes the form
\[
    \trace{\bfA \bfAregpinv(\lambda_k,k)} = \sum_{i=1}^k \phi_i^\text{\tiny G}.
\]

To obtain an expression for the norm of the residual, we exploit the fact that $\bfb = \beta_1 \bfU_{k+1} \bfe_1$ and that $\bfU_{k+1}$ contains orthonormal columns to get
    \begin{align}
      \|\bfr(\bfx(\lambda_k))\|^2
      & = \|\widetilde \bfU \bfSigma^\text{\tiny G} ((\bfSigma^\text{\tiny G})\t\bfSigma^\text{\tiny G}+\lambda\bfI_k)^{-1} (\bfSigma^\text{\tiny G})\t  \widetilde\bfU\t \bfb - \beta_1 \bfU_k \bfe_1\|^2 \nonumber\\
      & = \| (\bfSigma^\text{\tiny G} ((\bfSigma^\text{\tiny G})\t\bfSigma^\text{\tiny G}+\lambda_k\bfI_k)^{-1} (\bfSigma^\text{\tiny G})\t  -\bfI_k)  (\bfU^\text{\tiny G})\t \beta_1  \bfe_1\|^2.\nonumber
    \end{align}
Let $\bfp = (\bfU^\text{\tiny G})\t \beta_1  \bfe_1 \in \bbR^{k+1}$ with elements $[\bfp]_i$, then
\[
  \|\bfr(\bfx(\lambda_k))\|^2 = [\bfp]_{k+1}^2 + \sum_{i=1}^k (\phi_i^\text{\tiny G} - 1)^2 [\bfp]_i^2.
\]

\subsubsection{Methods that exploit the links between GKB and Gaussian quadrature rules}\label{sub:quadRule}

All the parameter choice rules for $\lambda$ presented so far involve at least one evaluation of a quadratic form of the matrix $\bfA\t\bfA$ or $\bfA\bfA\t$. These forms can be expressed as particular Riemann-Stieltjes integrals (written in terms of the SVD of $\bfA$), which can then be approximated by Gauss or Gauss-Radau quadrature rules. The latter are related to the symmetric Lanczos algorithm, which is itself related to the GKB 
algorithm; see \cref{sub:buildsubspace}. 
When dealing with large matrices $\bfA$ whose SVD cannot be computed, these relations are crucial to provide computationally affordable approximations of the original quadratic forms in terms of quantities (such as the upper bidiagonal $\bfB_k$) generated at the $k$th GKB iteration. Under some additional assumptions, it can be proven that such approximations are upper or lower bounds, which get increasingly accurate as $k$ increases.
Such connections are underpinned by a solid and elegant theory. In the following we will mention some basic facts without providing many motivations: these are fully unfolded in, e.g., \cite{MMQ}.
In the framework of hybrid projection methods, besides using $k$ GKB iterations to compute an approximation $\lambda_k$ of the regularization parameter $\lambda$ for the Tikhonov problem \cref{eq:Tikhonov} according to some rule, GKB is also employed to define the $k$th projected problem \cref{eq:LSQRproj} or \cref{eq:LSMRproj}.

Since all the parameter choice rules described for Tikhonov regularization in \cref{sub:paramNoise} and \cref{sub:paramNoNoise} involve the discrepancy norm $\|\bfr(\bfx(\lambda))\|$, and since
\begin{eqnarray}
\|\bfr(\bfx(\lambda))\|^2
&=& \lambda^2\bfb\t(\bfA\bfA\t+\lambda\bfI_m)^{-2}\bfb\nonumber\\
&=&\lambda^2\bfb\t\psi(\bfA\bfA\t)\bfb\,,\;\mbox{with}\; \psi(t)=(t+\lambda)^{-2}\,,\label{eq:phidef}
\end{eqnarray}
in the following we provide the details of the approximation of such a quadratic form and for an overdetermined $\bfA$ (i.e., $m\geq n$) only.
Considering the eigendecomposition $\bfA\bfA\t=\bfUA \bfSA(\bfSA)\t (\bfUA)\t$ given in terms of the SVD of $\bfA$ \cref{eq:svdA}, and setting
\[
\bfLA := \diag{\lambda_1^\text{\tiny A}, \ldots, \lambda_n^\text{\tiny A},0,\dots,0}=\bfSA(\bfSA)\t\in\bbR^{m\times m},\quad \bbfb:=\lambda(\bfUA)\t\bfb,
\]
it can be shown that
\begin{equation}\label{eq:int}
\lambda^2\bfb\t\psi(\bfA\bfA\t)\bfb=\bbfb\t\psi(\bfLA)\bbfb=\int_{0}^{\lambda_1^\text{\tiny A}}\psi(t)d\omega(t)=:I(\psi)\,,
\end{equation}
where the discrete measure $\omega(t)$ is a non-decreasing step function with jump discontinuities at the $\lambda_i^\text{\tiny A}$'s and at the origin. It is well-known that the above integral can be numerically approximated using a Gaussian quadrature rule, whose nodes are the zeros of a family of orthonormal polynomials with respect to the inner product induced by the measure $\omega(t)$. Such polynomials satisfy a three-term recurrence relation whose coefficients coincide with the entries of the tridiagonal matrix computed by the symmetric Lanczos algorithm applied to $\bfA\bfA\t$ with initial vector $\bfb$ (specifically, the entries on the $k$th column of such tridiagonal matrix define the $k$th orthonormal polynomial). Note that such decomposition can be conveniently obtained from the partial GKB decomposition, as shown in \cref{SymLanczos1}.

Consider the eigendecomposition $\bfT_{k,k}=\bfUB(\bfSB)^2(\bfUB)\t$ given in terms of the SVD $\bfB_{k,k}=\bfUB\bfSB(\bfVB)\t$. Then it is well-known that the $k$-point Gauss quadrature rule for approximating \cref{eq:int}
is given by
\begin{eqnarray}
\calG_k(\psi)&=&\lambda^2\|\bfb\|^2\sum_{j=1}^k\psi((\sigma_i^\text{\tiny B})^2)(\bfe_1\t\bfUB\bfe_k)^2=\lambda^2\|\bfb\|^2\bfe_1\t\bfUB\psi((\bfSB)^2)(\bfUB)\t\bfe_1\nonumber\\
&=&\lambda^2\|\bfb\|^2\bfe_1\t\psi(\bfT_k)\bfe_1=(\lambda\|\bfb\|\bfe_1)\t\psi(\bfB_{k,k}\bfB_{k,k}\t)(\lambda\|\bfb\|\bfe_1)\,,\label{Grule}
\end{eqnarray}
i.e., the quadratic form $\psi$ evaluated with respect to the projected matrices.
The $k$-point Gauss-Radau quadrature rule for approximating \cref{eq:int}, with one node at the origin, can be evaluated similarly, replacing the matrix $\bfB_{k,k}\bfB_{k,k}\t$ by the matrix $\bfB_{k-1}\bfB_{k-1}\t$, where $\bfB_{k-1}\in\bbR^{k\times (k-1)}$ is computed at the $(k-1)$st GKB iteration. That is,
\begin{equation}\label{GRrule}
\calR_k(\psi)=(\lambda\|\bfb\|\bfe_1)\t\psi(\bfB_{k}\bfB_{k}\t)(\lambda\|\bfb\|\bfe_1)\,.
\end{equation}
Since $\psi$ defined in \eqref{eq:phidef} is a $2k$-times differentiable function, the quadrature errors $\calE_{\calQ_k}(\psi):=I(\psi)-\calQ_k(\psi)$ associated with the $k$-point Gauss and Gauss-Radau quadrature rules (i.e., with $\calQ_k(\psi)=\calG_k(\psi)$ and $\calQ_k(\psi)=\calR_k(\psi)$, respectively), are given by
\begin{eqnarray*}\label{QError}
\calE_{\calG_k}(\psi)&=&\frac{\psi^{(2k)}(\zeta_{\calG_k})}{(2k)!}\int_{0}^{+\infty}\prod_{i=1}^k(t-\zeta_i)^2 d\omega(t)\,,\\
\calE_{\calR_k}(\psi)&=&\frac{\psi^{(2k-1)}(\bar{\zeta}_{\calR_k})}{(2k-1)!}\int_{0}^{+\infty}t\prod_{i=2}^k(t-{\bar{\zeta}_i})^2 d\omega(t)\,,
\end{eqnarray*}
where $\zeta_{\calG_k},\,\bar{\zeta}_{\calR_k}\in [0,\lambda_1^\text{\tiny A}]$, and the $\zeta_i$'s and the $\bar{\zeta}_i$'s denote the nodes of the Gauss and the Gauss-Radau quadrature rules, respectively. Since $\psi^{(2k-1)}(t)<0$ and $\psi^{(2k)}(t)>0$ for $t\geq 0$, the Gauss-Radau quadrature rule \cref{GRrule} is an upper bound for $\|\bfr(\lambda)\|^2$, while the Gauss quadrature rule \eqref{Grule} is a lower bound for $\|\bfr(\lambda)\|^2$.

The idea of exploiting Gaussian quadrature to approximate functionals used to set the Tikhonov regularization parameter was first presented by Golub and Von Matt \cite{golub1997generalized}, with derivations specific for GCV. This investigation can be regarded as a particular case of the broader explorations of the links between some Krylov methods and Gaussian quadrature conducted by Golub and collaborators; see again  \cite{MMQ} for a complete overview of the accomplishments in this area. An important remark concerning GCV is that a different treatment is needed to derive bounds for the numerator and denominator of the GCV function \cref{eq:GCV_original}: since a trace term appears in the denominator, a randomized trace estimator can be used to compute its approximation (see, e.g., \cite{golub1997generalized,randtrace}), while Gauss-quadrature bounds can be used to handle the numerator. Note however that, in \cite{golub1997generalized}, a similar approach is adopted to set the regularization parameter $\lambda$ for the full-dimensional problem \cref{eq:Tikhonov}, i.e., the value of $\lambda$ so approximated is then used to solve a Tikhonov-regularized problem. For this reason, this approach cannot be properly regarded as a hybrid solver, according to the framework adopted in this paper, although its adaptation would be quite natural. To the best of our knowledge,
Gaussian-quadrature-based regularization parameter choice strategies were first adapted to work within hybrid solvers in \cite{CaGoRe99}: here the algorithm involved in setting the regularization parameter according to GCV is identical to the one described in \cite{golub1997generalized}, and the computations performed to approximate the numerator of the GCV functional in \cref{eq:GCV_original} are then also employed to define an approximation subspace for the solution (of increasing dimensions), effectively making this approach a hybrid solver. A similar strategy, still tailored to GCV, is presented in \cite{Fenu16}.

Regularization parameter choice rules for hybrid methods that exploit the links between GKB and Gaussian quadrature rules are especially successful and well-understood when the $L$-curve is employed. This instance of the $L$-curve criterion is sometimes referred to as `$L$-ribbon'; see, e.g., \cite{calvetti2002curve, Calvetti2004}. Indeed, in this setting, Gauss and Gauss-Radau quadrature rules provide computationally inexpensive upper and lower bounds (i.e., boxes) for each point on the $L$-curve associated to the full-dimensional Tikhonov problem. Since the $L$-curve bounds are nested
(i.e., they become more and more accurate as the number of GKB iterations increases), when the $L$-ribbon is narrow it is possible to infer the approximate location of the `vertex' of the $L$-curve from its shape. For hybrid methods, the $L$-ribbon may first conveniently determine narrow bounds for the (direct) Tikhonov $L$-curve (thereby performing a suitable number of GKB iterations), and then easily locate the approximate vertex of the $L$-curve.

Most recently, the authors of \cite{gazzola2020survey} introduce a new principled and adaptive algorithmic approach for regularization, which provides reliable parameter choice rules by leveraging the framework of bilevel optimization, and the links between Gauss quadrature and GKB.

\subsection{Theoretical insights}
\label{sec:theory}
The goal of this section is to give the reader a flavor of the types of theoretical results that have been investigated and their significance for hybrid projection methods.

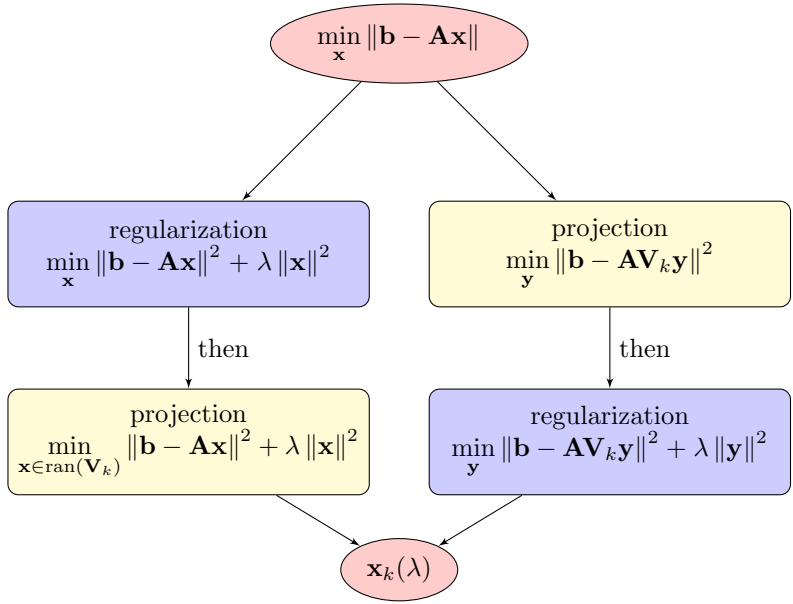
\begin{figure}
    \centering
    \begin{tikzpicture}[align=center, node distance = 2.8cm, auto]
        \node [cloud] (input) {$\displaystyle \min_\bfx \norm[]{\bfb - \bfA \bfx}$};
        \coordinate[below of = input] (c);
        \node [blockb, left of = c] (Reg1) {regularization \\$\displaystyle \min_\bfx \norm[]{\bfb - \bfA \bfx}^2 + \lambda \norm[]{\bfx}^2$};
        \node [blocky, right of = c] (Proj) {projection \\$\displaystyle \min_\bfy \norm[]{\bfb - \bfA \bfV_k \bfy}^2$};

        \node [blocky, below of = Reg1, node distance = 2.5cm] (Proj1) {projection \\ $\displaystyle \min_{\bfx \in {\rm ran}(\bfV_k)} \norm[]{\bfb - \bfA \bfx}^2 + \lambda \norm[]{\bfx}^2$};
        \node [blockb, below of = Proj, node distance=2.5cm] (Reg2) {regularization \\$\displaystyle \min_\bfy \norm[]{\bfb - \bfA \bfV_k \bfy}^2 + \lambda \norm[]{\bfy}^2$};
        \node [cloud, below of = input, node distance=7cm] (output) {$\bfx_k(\lambda)$};

        \path [line] (input) -- (Reg1);
        \path [line] (input) -- (Proj);
        \path [line] (Reg1) -- node {then} (Proj1);
        \path [line] (Proj) -- node {then} (Reg2);
        \path [line] (Proj1) -- (output);
        \path [line] (Reg2) -- (output);

    \end{tikzpicture}
\caption{For fixed regularization parameter $\lambda$ and in exact arithmetic, both approaches `first-regularize-then-project' and `first-project-then-regularize' result in the same solution approximation. A similar figure can be found in \cite{hansen2010discrete}.} \label{fig:theory}
\end{figure}

\paragraph{Project-then-regularize versus regularize-then-project}
For standard Tikhonov regularization, an insightful result deals with characterizing the iterates from a hybrid projection method. For a fixed regularization parameter,  there is an important result regarding the equivalence of iterations from two approaches: `first-regularize-then-project' and `first-project-then-regularize' \cite{hanke1993regularization,hansen2010discrete}. See \cref{fig:theory} and the following theorem.
\begin{theorem}
  Fix $\lambda >0$ and define $\bfx_k(\lambda)\in\bbR^n$
 to be the $k$th iterate of conjugate gradient applied to the Tikhonov problem,
 $$\min_\bfx \norm[]{\bfb - \bfA \bfx}^2 + \lambda \norm[]{\bfx}^2. $$  Let $\bfy_k(\lambda)\in\bbR^k$ be the exact solution to the regularized, projected problem
 $$\min_\bfy \norm[]{\beta_1\bfe_1 - \bfB_k \bfy}^2 + \lambda \norm[]{\bfy}^2, $$
where $\bfB_k$ and $\bfV_k$ are derived from GKB applied to the original problem.  Then $$\bfz_k(\lambda) = \bfV_k \bfy_k(\lambda) = \bfx_k(\lambda).$$
\end{theorem}
{The significance of this result is that, for large-scale problems where the `first-regularize-then-project' approach is not feasible (e.g., because obtaining a good regularization parameter is too expensive a priori, see \cref{fig:param_iterative}), the hybrid projection methods that follow the `first-project-then-regularize' approach produce, at the $k$th iteration, the same regularized solution (in exact arithmetic and with the same regularization parameter).} Extensions of this result to the case where TSVD regularization is used instead of Tikhonov can be found in \cite{KiOLe01}.

\paragraph{Projection results}
Since the approximation subspace for the solution plays a pivotal role in the success of both purely iterative and hybrid projection methods, there have been recent investigations on subspace approximation properties, including results quantifying how much the approximation subspace can be expanded, or how the noise in the data affects the approximation subspace. In particular, the authors of \cite{gazzola2015survey} prove that, when GKB is applied to a severely ill-posed problem with noiseless data, the product of the entries of the last column of the bidiagonal matrix $\bfB_k$ decays as $(k\sigma_k^\text{\tiny A})^2$; note that these entries are somewhat related to how much the Krylov subspace $\calK_k(\bfA\t\bfA,\bfA\t\bfb)$ is expanded at the $k$th iteration. Under similar assumptions, when employing the Arnoldi algorithm, the rate of decay of the last entry of the Hessenberg matrix $\bfH_k$ (i.e., $[\bfH_k]_{k+1,k}$) is comparable to $k\sigma_k^\text{\tiny A}$.
The authors of \cite{HnPlSt09} study how noise in the measurements (right-hand side of \cref{eq:linearmodel}) propagates to the projected problem and, in particular, how it enters the basis vectors spanning the approximation subspace.  This leads also to efficient estimates of the noise level that can be used to set the regularization parameter $\lambda_k$ for hybrid methods, or a stopping criterion (see \cref{sec:param}). Specifically, assuming that the Discrete Picard Condition holds, the noise amplification (i.e., the appearance of high frequency noise) is described as an effect of damping the smooth components due to convergence of Ritz values to large eigenvalues of $\bfA\t\bfA$ and the orthogonalization in the GKB algorithm. After the noise has `revealed' itself in the so called `noise-revealing' iteration, the LSQR relative residual can be used as a noise level estimator. This investigation is extended in \cite{hnvetynkova2017noise}, to consider the most popular iterative solvers based on GKB and to establish explicit relations between the noise-contaminated bidiagonalization vectors and the residuals of the considered regularization methods. In particular, it is shown that the coefficients of the linear combination of the computed bidiagonalization vectors reflect the amount of propagated noise in each of these vectors; influence of the loss of orthogonality is also discussed.
Some of these investigations (most relevantly, the analysis of the residual vector) have been also performed in \cite{GaNoRu14, GazzolaSilvia2016Iotd} for methods based on the Arnoldi algorithm, mainly using tools from approximation theory.
Finally, the propagation of noise in the approximation subspace and the residuals associated to purely iterative methods based on the GKB algorithm (such as LSQR) and the Arnoldi algorithm (such as GMRES and RRGMRES) is studied in \cite{JeHa07}.  The latter may not be successful in filtering the noise out, because of the so-called `mixing' of the SVD components in the approximation subspace; this issue may actually be irrelevant, if features of the solution can be reconstructed anyway thanks to the presence of favourable vectors in the Krylov approximation subspace and it can be mitigated using some `preconditioning' \cite{hansen2007smoothing}; see also \cref{sub:generalform}. Some `preconditioners' that can be used together with (RR)GMRES to correct for severe asymmetries in the discretized forward operator are proposed and analyzed in \cite{GazzolaSilvia2019AdGa}. Focusing on image deblurring problems, and exploiting the two-dimensional discrete cosine transform, \cite{HaJe08} also analyzes how noise from the data enters the solution computed by LSQR, GMRES, RRGMRES, MINRES, and MR-II,
concluding that the noise mainly affects the reconstructions in the form of low-pass filtered white noise, and discouraging the use of GMRES and MINRES for image deblurring.

Other classes of results pair the analysis of the approximation subspace for the solution with the analysis of the behavior of the projected problem, or focus solely on the latter, including results on the approximation of the singular value decomposition, which shed light on how the ill-conditioning of the original problem transfers to the projected problem, on the semi-convergence phenomenon, and on how properties such as the so-called discrete Picard condition are `inherited' by the projected problem. Bounds on the `residuals' associated to the SVD approximations obtained by the GKB and Arnoldi decompositions (i.e., the matrix 2-norm of the difference between the original coefficient matrix \cref{eq:svdA} and the approximated SVD obtained at each step of a Krylov subspace method) are presented in \cite{gazzola2015survey}, where it is concluded that, depending on the method and on the residuals, these could either be 0, or their norm may decay as the singular values of the original coefficient matrix. The authors of \cite{GazzolaSilvia2016Iotd} prove that, starting from a full-dimensional problem satisfying the discrete Picard condition (i.e., typically assuming noise-free data), the discrete Picard condition continues to hold when the problem is projected using the Arnoldi or the GKB algorithms: this is achieved by combining some SVD update rules and a backward induction argument (i.e., going from the full-dimensional problem to projected problems of decreasing dimension). Naturally, much analysis focuses on establishing theoretical and heuristic links between the performance of TSVD (which can be regarded as a least-squares solver projecting the full-dimensional problem \cref{LS} onto the subspace spanned by the dominant right singular vectors) and Krylov projection methods. Some early investigations, such as \cite{Hanke01}, find that, since the projection attained with the LSQR method is tailored to the specific right-hand side $\bfb$, this results in a more rapid convergence of the solver. Also, when it comes to computing TSVD-like approximations, the authors of \cite{Berisha2014} exploit a classical property of GKB (namely, the approximation quality depends on the relative distance between the singular values of $\bfA$ \cite{Saad80}) to conclude that, since for discrete inverse problems the relative gap between large singular values of $\bfA$ is generally much larger than the relative gap between its small singular values, the singular values of $\bfB_k$ converge very quickly to the largest singular values of $\bfA$. Also, links between the TSVD and Lanczos algorithms (for symmetric problems) are investigated in \cite{GaOnReRo16} where, thanks to the speed at which the nonnegative subdiagonal entries of $\bfT_k$ decay to zero, the solution of the discrete ill-posed problems can be efficiently and effectively (i.e., without significant, if any, reduction of the quality of the computed solution) expressed in terms of the Lanczos basis vectors rather than the eigenvectors of $\bfA$; this extends to the solution subspace determined by LSQR and the right singular vector subspace.
More recently, a body of papers by Jia and collaborators puts some renewed emphasis on the investigation of the relations between the TSVD solutions and the solutions computed by some popular purely iterative regularization methods (such as LSQR, LSMR, and MINRES). 
In particular, assuming distinct singular values, \cite{huang2017some} establishes bounds for the distance between the $k$-dimensional Krylov subspace associated to LSQR and the $k$-dimensional dominant right singular space associated to TSVD, and concludes that the former better `captures' the latter for severely and moderately ill-posed problems than for mildly ill-posed problems.  This implies that LSQR performs better as a stand-alone regularization method for severely and moderately ill-posed problems, but should be paired with additional regularization (i.e., in a hybrid framework) for successfully handling mildly ill-posed problems. Estimates for the accuracy of the rank-$k$ approximation generated by GKB are also provided. This topic is further investigated in \cite{JiaZhongxiao2020Tlra, JIA2020112786}, where an analysis of the approximation of the large singular values of $\bfA$ by the Ritz values is also performed to assess if this happens in natural order (i.e., largest singular values first). The paper \cite{JiaZhongxiao2020RpoL} extends the same investigations to the case of multiple singular values, with similar findings.

\paragraph{Continuous formulations}
Although the present survey paper considers discrete inverse problems, in many situations (especially when analyzing some convergence properties of solvers or inferring the behavior of the projected problems) looking at the continuous formulation (i.e., within a Hibert space $\calH$) can be beneficial. Indeed, it is well known that some discretization schemes (e.g., those based on boundary element methods) actually project the continuous problem onto finite dimensional vector spaces, and that such a projection has a regularizing effect; see \cite{hanke1993regularization, KiOLe01}. {Indeed, with `proper conversion', many of the presented derivations for Krylov methods in finite dimension can be translated to the continuous setting; see \cite{engl1996regularization}.}
 CG was historically among the first Krylov methods to be analyzed in a continuous setting: \cite{engl1996regularization} reports relevant results and references. More recently, the authors of \cite{GMRESreg} show that, when the error-free data vector $\bfb_\true$ lies in a finite-dimensional a finite-dimensional Krylov subspace, then the (continuous) GMRES method is a regularization method if the iterations are terminated by a stopping rule based on the discrepancy principle. Methods based on the Arnoldi algorithm (either FOM or GMRES) have been analyzed in a continuous setting \cite{No17}: here, some of the properties described (and somewhat heuristically justified) in a discrete setting in \cite{gazzola2015survey} are fully established (and sometimes strengthened) in a continuous setting. More specifically, assuming that we are dealing with Hilbert-Schmidt operators (like Fredholm integral equations of the first kind with square-integrable kernels), it can be shown that the rate of convergence of these solvers is related to the extendibility of the Krylov subspace $\calK_k(\bfA,\bfb)$, which is also comparable with the rate of decay of the singular values of $\bfA$. Moreover, it is proven that the dominant singular values can be approximated with improved accuracy as the iterations proceed, further supporting the statement that methods based on the Arnoldi algorithm are regularization methods. This investigation has an impact on the finite-dimensional setting, too: for instance, having proven convergence of the SVD approximation obtained by the Arnoldi algorithm, implies that many parameter choice rules that are intrinsically based on the success of this approximation (see \cref{sec:param}) are meaningful when used within hybrid methods.
An analogous investigation appeared in \cite{Caruso2019}, where the infinite-dimensional GKB algorithm and LSQR are considered.
The same author considers in \cite{NovatiP2018} the case of Krylov solvers applied to general-form Tikhonov regularization (still in a continuous setting), with a fixed regularization parameter; both solvers based on the Arnoldi algorithm and on the GKB algorithm are considered and the resulting scheme is dubbed ` Krylov-Tikhonov'. Denoting by $u^\dagger$ and $u_k$ the solution to the continuous problem and the Krylov solution, respectively, the main investigation in \cite{NovatiP2018} is concerned with proving that such Krylov methods are orthogonal projection methods for the linear operator equation associated to the normal equations for the continuous Tikhonov formulation: because of this, there exists a (semi)norm $E(\cdot)$ in $\calH$ such that $E(u_k-u^\dagger)$ converges to zero as $k$ goes to infinity.

\section{Extensions}
\label{sec:extensions}
The goal of this section is to describe some extensions and recent advancements in hybrid projection methods.
In particular, in \cref{sub:generalform} we describe various hybrid projection methods that have been developed for general-form Tikhonov regularization. In \cref{sub:recyclEnrich} we describe hybrid projection approaches that go beyond the standard projection subspaces by enrichment, augmentation and recycling of approximation subspaces for the solution.
Then in \cref{sub:flexible} we go beyond 2-norm regularization and consider the $\ell_p$ regularized problem.  In particular, we focus on recent advancements in flexible iterative methods (which are related to flexible preconditioning where the preconditioner changes during the iterative process) and extensions to hybrid frameworks.  In \cref{sub:bayesian} we describe the increasingly important role that hybrid iterative methods play in the field of large-scale computational uncertainty quantification.  We first give a brief introduction to draw connections between traditional variational regularization and statistical (Bayesian) inverse problems, and we describe various scenarios where hybrid projection methods have enabled researchers to go beyond point estimates and perform efficient UQ.  In \cref{sub:nonlinear}, we describe the role that hybrid projection methods have played in solving nonlinear inverse problems, where the forward model is nonlinear.

A common theme in all the following subsections is that extensions of hybrid projection methods beyond the standard approaches described in \cref{sec:hybrid}, in order to be effective in a large-scale setting, require the efficient computation of an appropriate solution subspace and adaptive regularization parameter choice for the projected problem (i.e., within the solution subspace), giving rise to hybrid formulations. Although the computed solution subspaces may not be as straightforward as standard methods, they can incorporate some prior knowledge or meaningful information about the solution or the kind regularization functional to be considered, eventually leading to superior reconstructions with respect to the ones delivered by standard methods. Most often this means that the standard GKB or Arnoldi decompositions have to be modified, for example, to formally work with variable preconditioning, to enforce orthogonality or optimality properties in a different norm, or to incorporate a low-dimensional subspace representing prior information. Also, most often, such tools have already been developed for use in other settings, e.g., for the solution of well-posed problems (e.g., PDEs). Therefore, a considerable part of this section is devoted to reviewing such tools, tailoring them to the problems at hand, and drawing possible analogies between them.

 \subsection{Beyond standard-form Tikhonov: Hybrid projection methods for general-form Tikhonov}
 \label{sub:generalform}
Thus far we have focused mainly on the standard-form Tikhonov problem \eqref{eq:Tikhonov}, where the regularizer $\|\bfx\|^2$ enforces an overall small norm for the solution; however, other smoothness properties may be desired. Many researchers have considered the \emph{general-form Tikhonov problem},
\begin{equation}
  \label{eq:genTikhonov}
  \min_\bfx \norm[]{\bfb - \bfA \bfx}^2 + \lambda \norm[]{\bfL\bfx}^2\,,
\end{equation}
where $\bfL\in \bbR^{p \times n}$ is called the regularization matrix.
We assume that the null spaces of $\bfA$ and $\bfL$ intersect trivially, so that $[\bfA\t\,, \bfL\t ]\t$ has full column rank and the solution to \eqref{eq:genTikhonov} is unique.
Typical choices of $\bfL$ include discretizations of a differential operator (e.g., the discrete first or second derivative).
An important observation is that the general-form Tikhonov solution can be written as a spectral filtered solution, where both the filter factors and the basis for the solution are determined by the generalized singular value decomposition (GSVD) of $\{\bfA, \bfL\}$; see, e.g., ~\cite[Chapter 8]{hansen2010discrete}. However computing the GSVD is not always
feasible, e.g.,
for large-scale unstructured problems. In these settings, the most immediate approach is to apply an iterative solver (e.g., any of the solvers listed in \cref{sub:iterative}) to the equivalent least-squares formulation of (\ref{eq:genTikhonov}): however, a suitable value of $\lambda$ must be fixed ahead of the iterative method.
Since this is typically not the case, one should apply an iterative solver to \cref{eq:genTikhonov} repeatedly (once for every tested value of $\lambda$), as dictated by some well-known parameter choice strategies (such as the ones listed in \cref{sec:param}): this eventually results in a costly strategy, and alternatives have been sought in the literature.

A common approach to handle the general-form Tikhonov problem is to transform it to standard form \cite{elden1977}. That is, one computes
\begin{equation}\label{eq:stdform}
\bar{\bfy}_{\bfL}(\lambda) =\argmin_{\bar{\bfy}\in\bbR^p} \|{\bfA}\bfL^{\dagger}_{\bfA}\bar{\bfy}-\bar{\bfb}\|^2+{\lambda}\|\bar{\bfy}\|^2\,,\;\mbox{where}\;
\begin{array}{lcl}
\bar{\bfb} \!&=&\! \bfb-\bfA \bfx_{0}^{\bfL} \\
\bfx_{\bfL}(\lambda) \!&=&\! \bfL^{\dagger}_{\bfA} \bar{\bfy}_{\bfL}(\lambda)+\bfx_{0}^{\bfL}
\end{array},
\end{equation}
and where $\bfx_0^{\bfL}$ is the component of the regularized solution $\bfx_{\bfL}(\lambda)$ in the null space of $\bfL$. Here $\bfL^\dagger_\bfA$ is the $\bfA$-weighted generalized inverse of $\bfL$, defined by \linebreak[4]$\bfL^{\dagger}_\bfA = (\bfI_n - ( \bfA (\bfI_n - \bfL^\dagger \bfL))^\dagger \bfA)\bfL^\dagger$; see~\cite[Section~2.3]{Han97}. {Note that, if $\bfL$ is invertible,
then $\bfL_\bfA^\dagger = \bfL^{-1}$, and, if $\bfL$ has full column rank, then $\bfL_\bfA^\dagger = \bfL^\dagger$. If
$\bfL$ is an underdetermined matrix, then $\bfL_\bfA^\dagger$ may not be the same as $\bfL^\dagger$.}  In the Bayesian inverse problems literature (see \cref{sub:bayesian}), this particular change of variables is referred to as \emph{priorconditioning}~\cite{calvetti2005priorconditioners,calvetti07bayesian,calvetti2007introduction,calvetti2014inverse}, because the matrix $\bfL$ is constructed from the covariance matrix of the solution modeled as a random variable (see also \cref{subsub:connection}).  The main idea is to interpret $\bfL$ as a preconditioner; however, while preconditioning for iterative methods is often discussed in the context of accelerating iterative methods (typically via clustering of the eigenvalues), priorconditioners in \cref{eq:stdform} change the subspace for the solution by including information from the prior.  Indeed, in earlier work that leveraged this idea without collocating it within the Bayesian framework, this approach was dubbed `smoothing norm preconditioning' and was first adopted in connection with CGLS \cite{hanke1993regularization}, and then GMRES \cite{hansen2007smoothing}. For regularization functionals designed to promote edges, Krylov subspace methods may have poor convergence whereas the transformed problem amplifies directions spanning the columns of the prior covariance matrix, thereby improving convergence \cite{arridge2014iterated}. Once formulation \cref{eq:stdform} is established, hybrid projection methods can be applied, which define the approximation subspace for the solution with respect to the matrix $\bfA \bfL^{\dagger}_\bfA$. For instance, in the case of GKB, the $k$th approximation subspace for the solution is defined as follows:
\[
\mathcal{K}_k((\bfA \bfL_\bfA^{\dagger})\t \bfA \bfL_\bfA^{\dagger},(\bfA \bfL_\bfA^{\dagger})\t \bar{\bfb}) = {\rm span} \{((\bfA \bfL_\bfA^{\dagger})\t \bfA \bfL_\bfA^{\dagger})^{i-1} (\bfA \bfL_\bfA^{\dagger})\t \bar{\bfb} \}_{i=0,\dots,k-1}\,.
\]
Correspondingly, the $k$th iteration of the LSQR algorithm updates the matrix recurrence
\[ \bfA \bfL_\bfA^{\dagger} \bfV_k = \bfU_{k+1} \bfB_k \]
where $\bfB_k$ is $(k+1) \times k$ {lower} bidiagonal,
$\bfV_k$ has $k$ orthonormal columns that span \linebreak[4]$\mathcal{K}_k((\bfA \bfL_\bfA^{\dagger})\t \bfA \bfL_\bfA^{\dagger},(\bfA \bfL_\bfA^{\dagger})\t \bfb)$,
and $\bfU_{k+1}$ has $k+1$ orthonormal columns with \linebreak[4]$\bfU_{k+1}(\beta e_1) = \bar{\bfb}$.
{Imposing $\bar{\bfy}\in\mathcal{K}_k((\bfA \bfL_\bfA^{\dagger})\t \bfA \bfL_\bfA^{\dagger},(\bfA \bfL_\bfA^{\dagger})\t \bfb)$ in (\ref{eq:stdform}), i.e., taking $\bar{\bfy}=\bfV_k\bar{\bfz}$, and exploiting} unitary invariance of the 2-norm, the objective function in  (\ref{eq:stdform}) can be equivalently rewritten as the $O(k)$-sized least-squares problem
\begin{equation}\label{eq:stdformProj}
{\bar{\bfz}_{\bfL}(\lambda) =\argmin_{\bar{\bfz}\in\bbR^k}  \| \bfB_k\bar{\bfz} - \beta \bfe_1\|^2 + \lambda\|\bar{\bfz} \|^2 ,}\qquad \bar{\bfy}_{\bfL}(\lambda) =\bfV_k\bar{\bfz}_{\bfL}(\lambda)\,.
\end{equation}
For fixed $k$, any suitable parameter choice strategy can be used on the projected problem (\ref{eq:stdformProj}) to choose $\lambda$. If $\bar{\bfy}_{\bfL}(\lambda)$ is computed at fixed pre-sampled values of $\lambda$ (e.g., when the $L$-curve criterion is adopted), this may be done with short term recurrences on the iteration $k$ (i.e., in $k$, for fixed $\lambda$),
without storage of $\bfV_k$; see \cref{sec:param} and \cite{KiOLe01} for details.
However, in general, $k$ should be assumed small enough, so that problem \cref{eq:stdformProj} can be efficiently generated and solved.
For some problems, multiplication with $\bfA \bfL^{\dagger}_\bfA$ can be done by solving a least-squares problem followed by an oblique projection; see \cite{Han94,hansen2007smoothing} for details. Furthermore, there are some regularization operators $\bfL$ for which sparsity or structure can be exploited so that multiplications with $\bfL^{\dagger}_\bfA$ can be done efficiently (e.g., banded matrices, circulant matrices, and orthogonal projections \cite{calvetti2005invertible,elden1982weighted,morigi2007orthogonal,reichel2009simple}).
However, in general, performing multiplications with $\bfL^{\dagger}_\bfA$ may be computationally difficult especially for large-scale problems (e.g., when $\bfL$ is defined as a product of matrices or as a sum of Kronecker products).

In the framework of hybrid projection methods, alternative approaches to approximate the solution of (\ref{sub:generalform}) have been devised without resorting to a standard-form transformation. For instance, Kilmer, Hansen and Espa\~nol \cite{KiHaEs06} developed a hybrid projection method that simultaneously bidiagonalizes both $\bfA$ and $\bfL$ (here $\bfL$ can be rectangular and does not need to be full rank), using a joint bidiagonalization algorithm inspired by Zha \cite{zha1996computing}. More precisely, a partial joint bidiagonalization for $\bfA$ and $\bfL$ is updated at each iteration: at the $k$th iteration, it reads
\begin{equation}\label{eq:jointbd}
\bfA \bfZ_k = \bfU^{\bfA}_{k+1} \bfB_k^{\bfA}, \qquad \bfL \bfZ_k = {\bfU}^{\bfL}_k {\bfB}_k^{\bfL},
\end{equation}
where $\bfZ_k = [\bfz_1,\ldots,\bfz_k]$, $\bfU_{k+1}(\beta \bfe_1) = \bfb$, and $\beta = \| \bfb \|$. Although this algorithm only requires that products with $\bfA, \bfL$ and their transposes can be computed, the cost of generating the partial factorization \cref{eq:jointbd} is somewhat difficult to quantify: indeed, the $k$th step requires one call to LSQR to produce orthogonal projections that are needed to update the partial joint bidiagonalization (\ref{eq:jointbd}), and each LSQR iteration requires four matrix-vector products (with each of $\bfA$, $\bfA\t$, $\bfL$, and $\bfL\t$). 
Taking $\bfx \in {\rm ran}({\bfZ}_k)$, i.e., $\bfx = \bfZ_k \bfw$, in \cref{eq:genTikhonov} results in the following equivalent problems
\begin{equation}
    \label{eq:KHEtransform}
   \min_{\bfw\in\bbR^k} \left\| \left[ \begin{array}{c} \bfA \bfZ_k \\ \sqrt{\lambda} \bfL\bfZ_k \end{array} \right] \bfw
   - \left[ \begin{array}{c} \bfb \\ \bfzero \end{array} \right]  \right\|
  = \min_{\bfw\in\bbR^k} \left\| \left[ \begin{array}{c} \bfB_k^{\bfA} \\ \sqrt{\lambda} {\bfB}^{\bfL}_k \end{array} \right] \bfw -
   \left[ \begin{array}{c} \beta \bfe_1 \\ \bfzero \end{array} \right] \right\|\,.
   \end{equation}
Since the last problem in the above equations is a Tikhonov problem of dimension $O(k)$ involving only bidiagonal matrices, it is significantly cheaper to solve than \cref{eq:genTikhonov}, especially for $k$ small.
Let us denote by $\bfw_{k}(\lambda)$ the minimizer of \cref{eq:KHEtransform}, then $\bfx_k(\lambda)=\bfZ_k\bfw_k{(\lambda)}$ is the minimizer of \cref{eq:genTikhonov}, constrained to ${\rm ran}(\bfZ_k)$.  For fixed $\lambda$, the joint bidiagonal structure leads to short recurrence updates (in $k$) for $\bfw_k{(\lambda)}$, for the residual norm $\|\bfA \bfx_k{(\lambda)} - \bfb \|$, and for the regularization norm term $\|\bfL \bfx_k{(\lambda)} \|$: this is relevant if the solution $\bfx_k{(\lambda)}=\bfZ_k\bfw_k{(\lambda)}$ has to be computed at fixed pre-sampled values of $\lambda$, e.g., when applying the $L$-curve criterion. The discrepancy principle to set $\lambda=\lambda_{k}$ adaptively at each iteration can also be efficiently employed.
It can be shown that the approximation subspace ${\rm ran}(\bfZ_k)$ is meaningful, as
$\bfx_k{(\lambda)}$ `resembles' a truncated GSVD regularized solution to \cref{eq:linearmodel}.

Reichel, Sgallari, and Ye \cite{reichel2012tikhonov} describe an iterative approach that simultaneously reduces square matrices $\bfA$ and $\bfL$ to generalized Hessenberg matrices using a generalized Arnoldi process, originally introduced in \cite{genArnoldi}. The formal full factorizations associated to this method consist of
{\begin{equation}\label{genArnoldi}
\begin{array}{lclclcl}
\bfA\bfQ &=&\bfQ\bfH^{\bfA}, &\quad\mbox{where}\quad& h^{\bfA}_{i,j}&=0&\quad\mbox{if}\; i\geq 2j+1\,,\\
\bfL\bfQ &=&\bfQ\bfH^{\bfL}, &\quad\mbox{where}\quad& h^{\bfL}_{i,j}&=0&\quad\mbox{if}\; i\geq 2j+2\,,
\end{array}
\end{equation}}
and where $\bfQ$ is orthogonal. In practice, $k$ iterations of the generalized Arnoldi algorithm can be computed with initial vector $\bfq_1=\bfb/\|\bfb\|=\bfb/\beta$ to generate the first $k$ columns of the generalized Hessenberg matrices $\bfH^{\bfA}$ and $\bfH^{\bfL}$ appearing above, and the corresponding columns of $\bfQ$. The columns of $\bfQ$ span a so-called generalized Krylov subspace obtained by multiplying $\bfq_1$ by $\bfA$ and $\bfL$ in a periodic fashion. Let $\bfH^{\bfA}_k\in\bbR^{(2k+1)\times k}$ and $\bfH^{\bfL}_k\in\bbR^{(2k+2)\times k}$ be the principal submatrices of the matrices $\bfH^{\bfA}$ and $\bfH^{\bfL}$, 
and let $\bfQ_k=[\bfq_1,\dots,\bfq_k]\in\bbR^{n\times k}$ be obtained by taking the first $k$ columns of $\bfQ$. Taking $\bfx \in {\rm ran}({\bfQ}_k)$, i.e., $\bfx = \bfQ_k \bfw$, in \cref{eq:genTikhonov} results in the following equivalent problems
\begin{eqnarray*}
   \min_{\bfw\in\bbR^k} \left\| \left[ \begin{array}{c} \bfA \bfQ_k \\ \sqrt{\lambda} \bfL\bfQ_k \end{array} \right] \bfw
   - \left[ \begin{array}{c} \bfb \\ \bfzero \end{array} \right]  \right\|
  = \min_{\bfw\in\bbR^k} \left\| \left[ \begin{array}{c} \bfH_k^{\bfA} \\ \sqrt{\lambda} {\bfH}^{\bfL}_k \end{array} \right] \bfw -
   \left[ \begin{array}{c} \beta \bfe_1 \\ \bfzero \end{array} \right] \right\|\,.
\end{eqnarray*}
Note that the approximation subspace for the solution of \cref{eq:genTikhonov} can alternatively be generated starting from $\bfq_1=\bfQ\bfe_1=\bfA\bfb/\|\bfA\bfb\|$, leading to minor adjustments in the above projected problem (analogous to the range-restricted approach described in \cref{sub:buildsubspace}). Similar to the method in \cite{KiHaEs06}, the present approach allows adaptive selection of the regularization parameter $\lambda=\lambda_k$ at the $k$th iteration, with a negligible computational overhead. Different from \cite{KiHaEs06}, the present method can work only if both $\bfA$ and $\bfL$ are square (but this restriction can be overcome, e.g., by zero-padding). Furthermore, the present method does not require the action of $\bfA\t$ and $\bfL\t$, and no (inner) calls to LSQR are needed to generate the factorization \cref{genArnoldi}. While these are clear advantages of the solver presented in \cite{reichel2012tikhonov} compared to that presented in \cite{KiHaEs06}, the downsides are that the structure of the generalized Hessenberg matrices $\bfH_k^{\bfA}$ and $\bfH_k^{\bfL}$ cannot be immediately exploited (as was possible for bidiagonal matrices $\bfB_k^{\bfA}$ and $\bfB_k^{\bfL}$ in \cref{eq:jointbd}); moreover, to the best of our knowledge, no heuristic insight is available to justify the choice of the approximation subspace ${\rm ran}({\bfQ_k})$ for the solution (in particular, no links with truncated GSVD have been established). Despite this, some specific features of the solution of \cref{eq:genTikhonov} can be enhanced by pairing the present method with so-called `selective regularization', i.e., by augmenting the approximation subspace for the solution associated to this specific method; generic (i.e., not specific to the generalized Arnoldi algorithm) schemes for achieving augmentation are described in \cref{sub:recyclEnrich}.

Lampe, Reichel, and Voss \cite{lampe2012large} propose an alternative projection method for approximating the solution (\ref{eq:genTikhonov}), which is initiated by picking an approximation subspace $\calV\subset \bbR^n$ of small dimension $k\ll n$ spanned by the orthonormal columns of a matrix $\bfV_k\in\bbR^{n\times k}$. This can be done, for instance, by running $k$ GKB iterations, so that ${\rm ran}(\bfV_k)=\calK_k(\bfA\t\bfA,\bfA\t\bfb)$. Denote by $\bfx_k{(\lambda)}$ the solution of the constrained Tikhonov problem obtained by imposing $\bfx\in\calV$ in \cref{eq:genTikhonov}, i.e., take \linebreak[4]$\bfx_k{(\lambda)}=\bfV_k\bfy_k{(\lambda)}$, where
\begin{eqnarray*}
\bfy_k{(\lambda)}=\argmin_{\bfy\in\bbR^k} \left\| \left[ \begin{array}{c} \bfA \bfV_k \\ \sqrt{\lambda} \bfL\bfV_k \end{array} \right] \bfw
   - \left[ \begin{array}{c} \bfb \\ \bfzero \end{array} \right]  \right\|
  = (\bfV_k\t\underbrace{(\bfA\t\bfA+\lambda\bfL\t\bfL)}_{=:\bfT(\lambda)}\bfV_k)^{-1}\bfV_k\t\bfA\t\bfb
\end{eqnarray*}
Note that $\bfy_k(\lambda)$ can be computed efficiently by using a `skinny' QR factorization of $[(\bfA\bfV_k)\t, (\bfL\bfV_k)\t]\t$, which can be updated as $k$ increases. A suitable value of the regularization parameter $\lambda=\lambda_k$ can be determined using any of the strategies described in \cref{sec:param}. Once $\bfx_{k}(\lambda)$ is computed, the original search space $\calV\in\bbR^{n\times k}$ is expanded by the normalized gradient of the functional in \cref{eq:genTikhonov} evaluated at $\bfx_{k}(\lambda)$. Namely, one considers
\[
\bfV_{k+1}=[\bfV_k, \bfv_{new}]\,,\quad \mbox{where}\quad \bfv_{new}=\left(\bfT(\lambda_k)\bfx_k(\lambda_k)-\bfA\t\bfb\right)/\|\bfT(\lambda_k)\bfx_k(\lambda_k)-\bfA\t\bfb\|\,.
\]
Although $\bfv_{new}$ so defined is orthogonal to ${\rm ran}(\bfV_k)$, this can be enforced numerically using reorthogonalization techniques, if loss of orthogonality is a concern. In general, the approximation subspace ${\rm ran}(\bfV_k)$ for the solution of \cref{eq:genTikhonov} is not a Krylov subspace and it is therefore called a `generalized Krylov subspace': note however that  it is yet a different generalized Krylov subspace than the one used in \cite{reichel2012tikhonov}. The process of solving problem \cref{eq:genTikhonov} constrained to the updated approximation subspace $\calV$ is repeated, and expansions of $\calV$ are generated until a stopping criterion is satisfied.

Hochstenbach and Reichel \cite{hochstenbach2010iterative} use a partial GKB built with respect to $\bfA$ and $\bfb$ to generate, at the $k$th iteration, the solution subspace $\calK_k(\bfA\t\bfA,\bfA\t\bfb)$, onto which they project the regularization matrix $\bfL$ in order to compute an approximate solution.
A similar idea was used in \cite{viloche2014extension,gazzola2015survey}.
Note that, contrary to all the strategies summarized so far in the present section, when using these approaches, the matrix $\bfL$ does not enter the approximation subspace for the solution.

\subsection{Beyond standard projection subspaces: Enrichment, augmentation, and recycling}\label{sub:recyclEnrich}
The success of Krylov projection methods, and in particular the regularizing properties of Krylov methods, depends on the ability of the associated Krylov subspace to capture the most salient features of the solution.
In the previous section, we surveyed approaches to solve the general-form Tikhonov problem \cref{eq:genTikhonov}, where the choice of the matrix $\bfL \neq \bfI_n$ may affect the solution subspace for the iterative method (e.g., via subspace preconditioning or priorconditioning with a change of variables). Such transformations result in a Krylov subspace that is `better suited' for the problem \cite{calvetti2005priorconditioners,hansen2007smoothing}. This means that the approximate solution vectors can capture the essence or important features of the solution, which can lead to faster convergence properties since fewer iterations are needed to obtain an accurate solution.  However, in some applications, the regularized solution can be further improved by extending or enhancing the approximation subspace via enrichment, augmentation, and recycling techniques.

We begin with some motivating scenarios where such techniques may be desired.  First, experienced practitioners may have solution vectors containing important information about the desired solution (e.g., a low-dimensional subspace) from previous experiments or theoretical analyses.  By incorporating these solution basis vectors, the solution accuracy can be \textit{significantly} improved, depending on the quality of the provided vectors.  Second, for very large inverse problems where many iterations are needed or where there are many unknowns, one of the main computational disadvantages of most hybrid methods compared to standard iterative methods is the need to store the basis vectors for solution computation, namely $\bfV_k$ (using the same notations as in \cref{eq:genfact}). Furthermore, for problems with slow convergence, implicit restarts and augmentation can speed up convergence without affecting memory.  Third, and this is related to \cref{sub:nonlinear}, hybrid methods may be embedded within larger frameworks, e.g., optimal experimental design \cite{haber2003learning,HaberE2008Nmfe} or nonlinear inverse problems \cite{haber2000optimization}.  A sequence of inverse problems must be solved, e.g., where the forward model may be parameterized such that the change in the model from one problem to the next is relatively small, or the goal may be to compute and update solutions from streaming data. Rather than start each solution computation from scratch, recycling techniques can be used to improve the given subspace and to compute a regularized solution efficiently in the improved subspace.
In the context of solving inverse problems, we describe extensions that go beyond the standard projection subspaces.

We split our discussion into two main classes: (1) enrichment methods and (2) augmentation or recycling approaches. For both classes of methods, we assume that we are given a low-dimensional subspace that represents prior information or expert knowledge. That is, we assume that we are given a suitable set of basis vectors, denoted $\bfW_p \in \bbR^{n \times p}$, such that a well-chosen subspace for the solution is given by $\calW_p = {\rm ran}(\bfW_p)$.
We call this the `prior solution subspace' (not to be confused with priors or priorconditioning methods).
The main assumption is that the solution has a significant component in the given subspace $\calW_p$, and the goal is to incorporate the provided subspace judiciously in order to obtain improved regularized solutions.

\paragraph{\textbf{Enrichment methods}}
We begin with enrichment methods, which we define as methods that, at the $k$th iteration, compute regularized solutions in a solution subspace that is a direct sum of the prior solution subspace and a standard $k$-dimensional Krylov subspace, i.e.,
\begin{equation}
  \label{eq:enhancedsolnspace}
  \bfS_{p,k} = \calW_p \oplus \calK_k = \left\{\bfx +\bfy \mid \bfx \in \calW_p \mbox{ and } \bfy \in \calK_k\right\}.
\end{equation}
That is, hybrid projection methods with enrichment construct iterates defined as
\begin{equation}
  \label{eq:enhanceopt}
 \bfx_k(\lambda_k) = \argmin_{\bfx \in \bfS_{p,k}} \norm[]{\bfb - \bfA \bfx}^2 + \lambda_k \norm[]{\bfx}^2.
 \end{equation}
 Enrichment methods can improve the solution accuracy by incorporating information about the desired solution into the solution process, i.e., by enriching the solution subspace. Methods for enriching Krylov subspaces have been described in the context of both Arnoldi and Lanczos projection methods.

For Arnoldi based methods, the Regularized Range-Restricted GMRES method (R3GMRES) \cite{dong2014r3gmres} extends the RRGMRES approach (see \cref{sub:buildsubspace}) to include a subspace that represents prior information about the solution where the approach is to compute as $\bfx_k(\lambda_k)$ as in \cref{eq:enhanceopt} with $\bfS_{p,k} = \calK_k(\bfA, \bfA \bfb)$.  Efficient implementations are described in \cite{dong2014r3gmres}.

For Lanczos-based projection methods, enrichment methods seek solutions of the form \cref{eq:enhanceopt} where the solution space is given by
\cref{eq:enhancedsolnspace} with $\calK_k(\bfA\t \bfA, \bfA\t \bfb)$.
A modification of CGLS method (sometimes also denoted in the literature with the acronym CGNR) was described in \cite{calvetti2003enriched}, where Krylov subspaces are combined with vectors containing important information about the desired solution (e.g., a low-dimensional subspace). Given an enrichment subspace $\calW_p$, the enriched CGLS method determines the $k$th iterate in the subspace $\calW_p \oplus \calK_k(\bfA\t \bfA, \bfA\t \bfb).$ The enriched CGLS method was used for Tikhonov problems (e.g., penalized least-squares problems), where several linear systems were solved for different choices of the regularization parameter as described in Frommer and Maass \cite{FrMa99}.  More specifically, a sequence of Tikhonov problems for regularization parameters $\lambda = \lambda_j = 2^{-j}$ for $j=0,1,\ldots$ are solved, until for some $\lambda_j,$ the corresponding solution satisfies the discrepancy principle. A hybrid enriched bidiagonalization (HEB) method that stably and efficiently augments a `well-chosen enrichment subspace'
with the standard Krylov basis associated with LSQR is described in \cite{hansen2019hybrid}. The idea is to use the Lanczos bidiagonalization algorithm to compute an orthonormal basis, and augment it by a user-defined low-dimensional subspace that captures desired features of the solution in each step of the algorithm.

Finally, we remark that a subspace-restricted SVD was described in \cite{hochstenbach2010subspace}, where a modification of the singular value decomposition permits a specified linear subspace to be contained in the singular vector subspaces. That is, the user can prescribe some of the columns of the $\bfU$ and $\bfV$ matrices for all truncations, and truncated versions as well as Tikhonov problems can be solved. For inverse problems, such approaches can give more accurate approximations of the solution when compared to standard TSVD and Tikhonov.

\paragraph{\textbf{Augmentation and recycling techniques}}
The development of augmentation schemes is deeply intertwined with the development of recycling techniques, in particular, Krylov subspace recycling methods for solving sequences of linear systems with slowly changing coefficient matrices, multiple right hand sides, or both.  We point the interested reader to recent surveys \cite{soodhalter2014krylov,soodhalter2020survey} and references therein.
Subspace recycling methods can exhibit accelerated convergence by enabling effective reuse of subspace information that is typically generated during a Krylov method on one or more of the systems.
For example, for well-posed least-squares problems that require many LSQR iterations, augmented LSQR methods have been described in \cite{baglama2013augmented,baglama2005augmented} that augment Krylov subspaces using harmonic Ritz vectors that approximate singular vectors associated with the small singular values, which can reduce computational cost by using implicit restarts for improved convergence.
However, when applied to ill-posed inverse problems, the augmented LSQR method without an explicit regularization term exhibits semi-convergent behavior.
Thus, in this setting, the goal is not so much to recycle spectral information for faster convergence, but rather to  augment a solution space to include certain known features about the solution.

Next, we provide an overview of augmented Krylov methods for solving ill-posed problems with a specific focus on hybrid methods.
Suppose we are given a general $p$-dimensional subspace $\calW_p={\rm ran}(\bfW_p)$, where $\bfW_p$ has orthonormal columns, and define
$\bfP_{\bfW_p} = \bfW_p\bfW_p\t$ and $\bfP_{\bfW_p}^\perp = \bfI_n - \bfW_p\bfW_p\t$ to be the orthonormal projectors onto $\calW_p$ and $\calW_p^\perp$ respectively. The decomposition approach described in \cite{baglama2007decomposition} and the augmented Krylov methods split the solution space such that the solution at the $k$th iteration is partitioned as
\[
  \bfx_k = \hat \bfx_k + \tilde \bfx_k, \quad \mbox{where} \quad \hat \bfx_k = \bfP_{\bfW_p} \bfx_k \quad \mbox{and}\quad  \tilde \bfx_k = \bfP_{\bfW_p}^\perp \bfx_k.
\]
It is assumed that the augmented space is small (e.g., $1 \leq p\leq 3$) and does not require regularization, thus a direct method can be used to compute $\hat \bfx_k$  To compute $\tilde \bfx_k$, an iterative method is used.  More specifically, let
\begin{equation}
  \label{eq:AWQR}
  \bfA \bfW_p = \bfV_p \bfR
\end{equation}
be the `skinny' QR factorization, where $\bfV_p \in \bbR^{n \times p}$ contains orthonormal columns and $\bfR \in \bbR^{p \times p}$ is upper triangular.
Then the approximate solution for the augmented GMRES approach is obtained by solving the constrained least-squares problems,
\[
  \min_\bfx \norm[]{\bfA \bfx - \bfb}^2, \quad \bfx \in \calW_p \oplus \calK_k(\bfP_{\bfV_p}^\perp \bfA, \bfP_{\bfV_p}^\perp \bfb)\,;
\]
similarly, the augmented RRGMRES method solves
\[
  \min_\bfx \norm[]{\bfA \bfx - \bfb}^2, \quad \bfx \in \calW_p \oplus \calK_k(\bfP_{\bfV_p}^\perp \bfA, \bfP_{\bfV_p}^\perp \bfA \bfb)\,.
\]
More details are given in  \cite{baglama2007decomposition,baglama2007augmented}.
Augmented LSQR methods have been described in \cite{baglama2013augmented,baglama2005augmented}.

A general framework to combine a recycling GKB (recycling GKB) process with tools from compression in a hybrid context are described in \cite{jiang2021hybrid}. More specifically, the approach consists of three steps.
First, a suitable set of orthonormal basis vectors $ \bfW_{p}$ is provided (e.g., from a related problem or from expert knowledge) or may need to be determined (e.g., via compression of previous solutions). This prior solution subspace may even include an initial guess $\bfx^{\ast}$ for the solution. The second step is to use a recycling GKB process to generate the columns of $\widetilde{\bfV}_{\ell}\in \bbR^{n \times \ell}$ that span the Krylov subspace
$$\calK_k(\bfA\t\bfP_{\bfV_p}^\perp \bfA, \bfA\t \bfP_{\bfV_p}^\perp \bfr)\,,$$
where $\bfr = \bfb - \bfA \bfx^{\ast}$ and $\bfV_p$ is defined in \cref{eq:AWQR}.
The third step is to find a suitable regularization parameter $\lambda_k$ and compute a solution to the regularized projected problem in the extended solution space.
The recycled, hybrid iterate is given by
\[
   \bfx_k(\lambda_k) =\argmin_\bfx \norm[]{\bfA \bfx - \bfb}^2 + \lambda_k \norm[]{\bfx}^2, \quad \bfx \in \calW_p \oplus \calK_k(\bfA\t\bfP_{\bfV_p}^\perp \bfA, \bfA\t \bfP_{\bfV_p}^\perp \bfr).
\]
This process can be repeated in a cyclic fashion, and it is successful in solving a broad class of inverse problems (e.g., solving very large problems where the number of basis vectors becomes too large for memory storage, solving a sequence of regularized problems (e.g., changing regularization terms or nonlinear solvers), and solving problems with streaming data).

\paragraph{\textbf{Let's compare and contrast}}
For all of the described methods, the improvement in solution accuracy \textit{significantly} depends on the quality of the provided prior solution subspace vectors $\bfW_p$.  Thus, we reiterate that the major advantage as well as the major concern for all enrichment, augmentation, and recycling methods is that the augmented subspace needs to be chosen appropriately, i.e., problem-dependent or prior information about the solution needs to be known.  The performance of these methods compared to standard methods will depend on the amount of relevant and important information about the solution that is encoded in the given subspace.
Some suggestions for subspaces are provided in \cite{hochstenbach2010subspace,jiang2021hybrid}.  For problems with unsuitable augmentation spaces, an adaptive augmented GMRES was considered in \cite{kuroiwa2008adaptive}, which automatically selects a suitable subspace from a set of user-specified candidates.

The main difference between enrichment methods and augmentation methods lies in how the solution space is constructed.  Iterates of enrichment methods belong to a direct sum of the prior solution subspace and a standard Krylov subspace, whereas the generated subspace for augmented methods is restricted to the orthogonal complement of $\calV_p.$
It may be argued that since augmentation and recycling methods generate subspace vectors that incorporate the provided prior solution subspace, these methods can \emph{improve} on the generated space.  However, this observation is problem dependent. Further investigations regarding the relationship between spaces $\calW_p$ and $\calV_p$ can be found in \cite{dong2014r3gmres}.

We also mention that there have been various works on developing hybrid methods for solving systems with multiple right hand sides.  In \cite{OlSi81}, the authors suggest using the block Lanczos algorithm  \cite{GoLuOv81}.  Some more recent work on extending generalized hybrid methods for multiple right hand sides exploit global iterative methods \cite{toutounian2006global} and block reduction methods \cite{alqahtani2021block}.  Such approaches have been considered for solving dynamic inverse problems \cite{chung2018efficient}. Furthermore, there is recent work on developing hybrid projection methods based on tensor-tensor products and applications to color image and video restoration \cite{reichel2021tensor,reichel2021tensorAT}.

\subsection{Beyond the 2-norm: Sparsity-enforcing hybrid projection methods for $\ell_p$ regularization}
\label{sub:flexible}
Although widely used, it is well known that standard-form Tikhonov regularization \cref{eq:Tikhonov} and even its general-form counterpart \cref{eq:genTikhonov} can be rather restrictive and that other regularization terms can yield better approximations of $\bfx_\true$ in \cref{eq:linearmodel}. For example, if $\bfx_\true$ is known to be sparse (i.e., when most of its entries $[\bfx_\true]_i$, $i=1,\dots,n$, are expected to be zero), the variational formulation
\begin{equation}\label{msl:eq:originalpbm}
 \min_\bfx \|\bfA\bfx-\bfb \|_2^2+\lambda \|\bfx\|_p^p,\quad\mbox{with $0<p\leq 1$},
\end{equation}
is commonly considered in the literature, especially since the advent of the compressed sensing theory; see, e.g., \cite{CaRoTa06}. Indeed, while sparse vectors have a small $\ell_0$-`norm', replacing the $\ell_p$-norm in \cref{msl:eq:originalpbm} by an $\ell_0$-`norm' would yield an NP-hard optimization problem (see, for instance, \cite{Fornasier2011}). A well-established practice is to consider $\ell_p$ (instead of $\ell_0$) regularization terms with $0<p\leq 1$ to enforce sparsity, noting that the objective function in (\ref{msl:eq:originalpbm}) is non-differentiable at the origin and it is nonconvex when $0 < p < 1$ (see, e.g., \cite{huang2017majorization, lanza2015generalized,lanza2017nonconvex}). Note that, if sparsity of the solution in a different domain (e.g., wavelets, discrete cosine transform, gradient) is desired, a regularization term of the form $\|\bfPsi\bfx\|_p^p$ should be considered, where $\bfPsi\in\bbR^{n\times n}$ is a so-called `sparsity' transform; this situation will be considered in \cref{ssec:transf}. Also, sometimes, an $\ell_q$-norm with $q\neq 2$ has to be used to evaluate the fit-to-data term in \cref{msl:eq:originalpbm}: this is appropriate, for instance, when the noise term in \cref{eq:linearmodel} is modeled as impulse noise; see, e.g., \cite{impulsenoise}.

The $\ell_2$-$\ell_p$ regularized problem \cref{msl:eq:originalpbm} can be solved using
a variety of general nonlinear solvers and optimization methods; see, for instance, \cite{beck2009fast,goldstein2009split,wright2009sparse} and the references therein. Here we focus on solvers that approximate the regularization term in (\ref{msl:eq:originalpbm}) by a sequence of weighted $\ell_2$ terms (i.e., iteratively reweighted schemes).
This approach intrinsically relies on the interpretation of problem (\ref{msl:eq:originalpbm}) as a nonlinear weighted least-squares problem of the form
\begin{equation}\label{msl:eq:originalpbm2}
 \min_\bfx \|\bfA\bfx-\bfb \|_2^2+\lambda \|\bfx\|_p^p= \min_\bfx \|\bfA\bfx-\bfb \|_2^2+\lambda \|\bfW^{(p)}(\bfx)\bfx\|_2^2\,,
\end{equation}
where the diagonal weights in $\bfW^{(p)}(\bfx)$ are defined as
\begin{equation}\label{msl:eq:weights}
\bfW^{(p)}(\bfx)=\text{diag} \left( (|[\bfx]_i|^{\frac{p-2}{2}})_{i=1,...,n}\right)\,.
\end{equation}
Since we are considering $ 0 < p\leq 1$, division by zero might occur if $[\bfx]_i=0$ for any $i \in \{1,...,n\}$ and, in fact, this is common situation in the case of sparse solutions. For this reason, instead of (\ref{msl:eq:weights}), we can consider the following closely related weights
\begin{equation}\label{msl:eq:weights2}
    \widetilde{\bfW}^{(p,\tau)}(\bfx)=\text{diag} \left( (([\bfx]_i^2 + \tau^2)^{\frac{p-2}{4}})_{i=1,...,n}\right)\,,
\end{equation}
where $\tau$ is a parameter that, to keep the derivations simple, we will consider fixed and chosen ahead of the iterations (although some authors adaptively choose it). {Problem (\ref{msl:eq:originalpbm2}) is then replaced by}
\begin{equation}\label{msl:eq:smoothpbm}
   \min_\bfx \underbrace{\|\bfA\bfx-\bfb \|_2^2+\lambda \|\widetilde{\bfW}^{(p,\tau)}(\bfx)\bfx\|_2^2}_{T^{(p,\tau)}(\bfx)},
\end{equation}
{where $\tau\neq 0$} ensures that $T^{(p,\tau)}(\bfx)$
is differentiable at the origin for $p > 0$.
Note that problem (\ref{msl:eq:originalpbm2}) can be recovered from problem (\ref{msl:eq:smoothpbm}) by setting $\tau=0$.

A well-established framework to solve problem (\ref{msl:eq:smoothpbm}) is the local approximation of $T^{(p,\tau)}$ by a sequence of quadratic functionals  $T_k(\bfx)$ that gives rise to a sequence of quadratic problems of the form
\begin{equation}\label{msl:eq:approxpbm}
 \bfx_{k,\ast} = \argmin_\bfx \underbrace{\|\bfA\bfx-\bfb \|_2^2+ \lambda \|\bfW_{k} \bfx\|_2^2+c_k}_{=:T_k(\bfx)}.
\end{equation}
where  {$\bfW_k=\widetilde{\bfW}^{(p,\tau)}(\bfx_{k-1,\ast})$}, $c_k$ is a constant (with respect to $\bfx$) term for the $k$th problem in the sequence, and $\lambda$ has absorbed other possible {multiplicative} constants so that $T_k(\bfx)$ in (\ref{msl:eq:approxpbm}) corresponds to a quadratic tangent majorant of $T^{(p,\tau)}(\bfx)$ in (\ref{msl:eq:smoothpbm}) at $\bfx=\bfx_{k-1,\ast}${, i.e., $T_k(\bfx)\geq T^{(p,\tau)}(\bfx)$ for all $\bfx\in\bbR^n$, $T_k(\bfx_{k-1,\ast})=T^{(p,\tau)}(\bfx_{k-1,\ast})$, and $\nabla T_k(\bfx_{k-1,\ast}) = \nabla T^{(p,\tau)}(\bfx_{k-1,\ast})$}; see also \cite{huang2017majorization,wohlberg2008lp}.
Because of this, iteratively reweighted schemes (\ref{msl:eq:approxpbm}) can be regarded as particular instances of majorization-minimization (MM) schemes, and convergence to a stationary point of the objective function in problem (\ref{msl:eq:smoothpbm}) can be guaranteed; see \cite{huang2017majorization, MM2014, lanza2015generalized}.

The vector $\bfx_{k,\ast}$ denotes the solution of (\ref{msl:eq:approxpbm}). For moderate-sized problems, or for large-scale problems where $\bfA$ has some exploitable structure, $\bfx_{k,\ast}$  may be obtained by applying a direct solver to (\ref{msl:eq:approxpbm}). However, for large unstructured problems, iterative solvers can be used in different fashions to approximate the solution of (\ref{msl:eq:approxpbm}): this leads to either schemes based on nested (inner-outer) iteration cycles for the sequence of problems (\ref{msl:eq:smoothpbm}) or schemes that adaptively incorporate updated weights at each iteration of a single iteration cycle. The former are detailed in \cref{ssec:inout}, while the latter are detailed in \cref{ssec:GKSlplq,ssec:flexible}.

\subsubsection{Strategies based on inner-outer iterations}\label{ssec:inout}
Iteratively Reweighted Least Squares (IRLS, \cite[Chapter 4]{bjorck1996numerical}, \cite{daubechies2010})
or Iteratively Reweighted Norm (IRN, \cite{IRNekki, wohlberg2008lp}) methods are very popular schemes that can handle the smoothed reformulation \cref{msl:eq:smoothpbm} of problem \cref{msl:eq:originalpbm}. Let $\bfx_{k,l}$ denote the approximate solution at the $l$th iteration (of an inner cycle of iterations) for the $k$th problem of the form \cref{msl:eq:approxpbm} (i.e., at the $k$th iteration of an outer cycle of iterations); one typically takes $\bfx_{k,\ast}=\bfx_{k,l}$ when a stopping criterion is satisfied for the inner iteration cycle. IRLS or IRN methods based on an inner-outer iteration scheme are very popular  and have been used in combination with different inner solvers, such as steepest descent and CGLS \cite{fornasier2016cg, wohlberg2008lp}.

If $\bfW_k$ is square and invertible (note that this can be assumed {when} the weights are defined as in \cref{msl:eq:weights2} with $\tau>0 $ for any fixed $p > 0$), problem \cref{msl:eq:approxpbm} can be easily and conveniently transformed into standard form as follows
\begin{equation}\label{msl:eq:pbm_stdform}
 \bar{\bfx}_{k,\ast} = \argmin_{\bar{\bfx}} \|\bfA \bfW_{k}^{-1} \bar{\bfx} -\bfb \|_2^2+\lambda \| \bar{\bfx}\|_2^2\,, \quad \text{so that} \quad \bfx_{k,\ast} = \bfW_{k}^{-1}\bar{\bfx}_{k,\ast}.
\end{equation}
The interpretation of the matrix $\bfW^{-1}_k$ as a right preconditioner for problem \cref{msl:eq:approxpbm} can be exploited under the framework of priorconditioning; see \cref{sub:generalform} and \cite{calvetti07bayesian}. Similar to the approaches mentioned for \cref{msl:eq:approxpbm}, Krylov methods can be straightforwardly adopted to compute an approximate solution of the equivalent $k$th problem of the form \cref{msl:eq:pbm_stdform} (i.e., at the $k$th iteration of an outer cycle of iterations) through an inner cycle of iterations. In particular, variants based on hybrid projection methods (such as hybrid GMRES or hybrid LSQR) have been successfully considered, with the added benefit that adaptive regularization parameter choice strategies (based, e.g., on the discrepancy principle, UPRE, and GCV) can be considered; see \cite{ReVaAr17}. Specifically, at the $l$th iteration of a Krylov method for (\ref{msl:eq:pbm_stdform}), let $\bfV_{k,l} \in \mathbb{R}^{n \times l}$ be the matrix whose columns, span a (preconditioned) Krylov Subspace $\mathcal{K}_{k,l}$ of dimension $l$; the precise definition of $\mathcal{K}_{k,l}$  depends on the selected projection method. Problem (\ref{msl:eq:pbm_stdform}) can be projected and solved in $\mathcal{K}_{k,l}$ by {computing}
\begin{equation}\label{msl:eq:ykm1}
 \bar{\bfy}_{k,l} = \argmin_{\bar{\bfy}} \|\bfA \underbrace{\bfW_{k}^{-1} \bfV_{k,l}}_{=:\bfZ_{k,l}} \bar{\bfy} -\bfb \|_2^2+\lambda \| \bfV_{k,l} \bar{\bfy}\|_2^2\,,
\end{equation}
so that $\bar{\bfx}_{k,l} = \bfV_{k,l} \, \bar{\bfy}_{k,l}$, and
\[
\bfx_{k,l} = \bfW_{k}^{-1}\, \bar{\bfx}_{k,l} = \bfW_{k}^{-1} \,\bfV_{k,l} \, \bar{\bfy}_{k,l}=\bfZ_{k,l} \, \bar{\bfy}_{k,l}\,.
\]

Once $\bfx_{k,l}$ has been computed and once the weights $\bfW_{k+1}=\widetilde{\bfW}^{(p,\tau)}(\bfx_{k,l})$ have been updated according to  (\ref{msl:eq:weights2}), the $(k+1)$st problem of the form \cref{msl:eq:pbm_stdform}, with $k$ replaced by $k+1$, needs to be solved: to achieve this, one needs to compute a new Krylov subspace from scratch, defined with respect to the new coefficient matrix $\bfA\bfW_{k+1}^{-1}$. For this reason, this approach may become
computationally expensive.

\subsubsection{Strategies based on generalized Krylov methods}\label{ssec:GKSlplq} 
One way to avoid inner-outer iterations when applying the IRN method described in \cref{ssec:inout} is to exploit generalized Krylov subspaces (GKS), like the ones described in \cref{sub:generalform}: this approach underlies the so-called GKSpq \cite{lanza2015generalized} and MM-GKS \cite{huang2017majorization} algorithms.

Given the initial subspace $\calK_h(\bfA^\top \bfA,\bfA^\top \bfb)={\rm ran}(\bfV_0)$ generated by the GKB algorithm, where $h$ is small (typically $h\leq 5$) and $\bfV_0\in\bbR^{n\times h}$ has orthonormal columns, at the $k$th iteration of the MM-GKS method, a GKS $\calV_k$ spanned by the orthonormal columns of $\bfV_k\in\bbR^{n\times (h+k)}$ is computed, and it is used to define the following projection of problem \cref{msl:eq:approxpbm},
\[
\bfy_k=\argmin_{\bfy\in\bbR^{h+k}}\|\bfA\bfV_k\bfy-\bfb \|_2^2+ \lambda \|\bfW_{k}\bfV_k \bfy\|_2^2 .
\]
Equivalently, introducing the `skinny' QR factorizations
\begin{equation}\label{eq:QRgks}
\bfA\bfV_k = \bfQ^{\bfA}\bfR^{\bfA}\,,\quad \bfW_k\bfV_k = \bfQ^{\bfW}\bfR^{\bfW}\,,
\end{equation}
the above problem can be reformulated as the low-dimensional (if $h+k\ll\min\{m,n\}$) problem
\[
\bfy_k=\argmin_{\bfy\in\bbR^{h+k}}\|\bfR^{\bfA}\bfy-(\bfQ^{\bfA})\t\bfb \|_2^2+ \lambda \|\bfR^{\bfW} \bfy\|_2^2\,.
\]
After the $k$th approximation to problem \cref{msl:eq:approxpbm} is formed by taking $\bfx_k=\bfV_k\bfy_k$, one computes the residual of the normal equations associated to \cref{msl:eq:approxpbm}
\[
\bfr_{k}= (\bfA\t\bfA + \lambda\bfW_k^2)\bfx_k - \bfA\t\bfb\,.
\]
The subspace $\calV_{k+1}={\rm ran}(\bfV_{k+1})$ is then expanded by adding the vector $\bfr_k/\|\bfr_k\|_2$;
reorthogonalization of the columns of $\bfV_{k+1}$ may be advisable to enforce accuracy.
An initial guess $\bfx_0$ for $\bfx_\true$ is needed to define the first weights $\bfW_1$, and the authors of \cite{huang2017majorization} suggest to take $\bfx_0=\bfb$, although this only works when $\bfA$ is square. We emphasize that the regularization parameter $\lambda$ can be adaptively set during the GKS iterations; see \cite{impulsenoise, buccini2019dp}. Also, an expression for the quadratic  tangent  majorant \cref{msl:eq:approxpbm} different from $T_k(\bfx)$ can be devised, which allows for a cheaper GKS computation, in that smart updates of the QR factorization \cref{eq:QRgks} are performed: we refer to \cite{huang2017majorization} for the details. Finally, such strategies based on GKS can also handle generalizations of problem \cref{msl:eq:originalpbm} where weights are also included in the fit-to-data term.

\subsubsection{Strategies based on flexible Krylov methods}\label{ssec:flexible}
Another approach that avoids inner-outer iterations when applying the IRN method described in \cref{ssec:inout} is to exploit flexible Krylov subspaces (FKS): this approach underlies hybrid methods based on the flexible Arnoldi algorithm \cite{gazzola2014generalized} and the flexible Golub-Kahan algorithm \cite{chung2019}, as well as the more recent IRN-FKS methods \cite{FKSIRW}. The starting point of this class of methods is the formulation \cref{msl:eq:pbm_stdform} of the $k$th reweighted problem. Flexible Krylov methods, which are classically employed to handle iteration-dependent preconditioning \cite{notay2000flexible, Saad1993, simoncini2007recent}, provide a natural framework to update the inverted weights $\bfW_k^{-1}$ as soon as a new approximation of $\bfx_\true$ is computed. That is, at each iteration, the updated inverted weights are immediately incorporated within the approximation subspace for the solution to the next problem in the sequence.

Formally, both the flexible Arnoldi and the flexible Golub-Kahan algorithms look very similar to their standard counterparts \cref{eq:Arnoldi} and \cref{eq:GKB}, respectively. If $\bfA$ is square, taking $\bfv_1=\bfb/\norm[2]{\bfb}$ and defining $\bfW_1$ using an initial guess $\bfx_0$ ($\bfW_1=\bfI_n$ if $\bfx_0=\bfzero$), the $i$th iteration of the flexible Arnoldi algorithm \cite{Saad1993} applied to problem \cref{msl:eq:pbm_stdform} computes
\begin{equation}\label{eq:fArnoldi1}
\bfz_i=\bfW_i^{-1}\bfv_i\,,\quad \bfv = \bfA\bfz_i\,,\; \bfv=(\bfI_n - \bfV_i\bfV_i\t)\bfv\,,\quad \bfv_{i+1}=\bfv/\norm[2]{\bfv}\,,
\end{equation}
where $\bfW_i=\widetilde{\bfW}^{(p,\tau)}(\bfx_{i-1})$ is defined as in (\ref{msl:eq:weights2}), $\bfx_{i-1}$ being the solution obtained at the previous iteration, and where $\bfV_i=[\bfv_1, ..., \bfv_i]\in\bbR^{n\times i}$ has orthonormal columns.
$k$ iterations of the flexible Arnoldi algorithm can be equivalently expressed as the partial matrix factorization
\begin{equation}\label{eq:fArnoldi2}
\bfA\bfZ_{k}=\bfV_{k+1}\bfH_{k}\,,\quad\mbox{where}\quad
\begin{array}{lcll}
\bfZ_{k} &=&[\bfW_1^{-1} \bfv_1, ..., \bfW_k^{-1} \bfv_k] &\in\bbR^{n\times k}\\
\bfH_k &=&\bfV_{k+1}\t\bfA\bfZ_k &\in\bbR^{(k+1)\times k}
\end{array}.
\end{equation}
Here the columns of $\bfZ_{k}$ span the approximation subspace for the solution, $\bfV_{k+1}$ is defined as in \cref{eq:fArnoldi1} and, by construction, $\bfH_k$ is upper Hessenberg.
If $\bfA$ is rectangular, taking $\bfu_1=\bfb/\norm[2]{\bfb}$, the $i$th iteration of the flexible Golub-Kahan algorithm \cite{chung2019} applied to problem \cref{msl:eq:pbm_stdform} computes
\begin{equation}\label{eq:partialFGKB1}
\begin{array}{llll}
&\bfv= \bfA\t\bfu_i\,,\; &\bfv=(\bfI_n - \bfV_{i-1}\bfV_{i-1}\t)\bfv\,,\; &\bfv_{i}=\bfv/\norm[2]{\bfv},\\
\bfz_i=\bfW_i^{-1} \bfv_i\,,\; & \bfu= \bfA\bfz_i\,,\; &\bfu=(\bfI_m - \bfU_i\bfU_i\t)\bfu\,,\; &\bfu_{i+1}=\bfu/\norm[2]{\bfu},
\end{array}
\end{equation}
where the inverted weights $\bfW_i^{-1}$ are updated as in the flexible Arnoldi case, and both $\bfV_i=[\bfv_1, ..., \bfv_i]\in\bbR^{n\times i}$ and $\bfU_i=[\bfu_1, ..., \bfu_i]\in\bbR^{m\times i}$ have orthonormal columns.
$k$ iterations of the flexible Golub-Kahan algorithm can be equivalently expressed as the partial matrix factorizations
\begin{equation}\label{eq:partialFGKB}
\begin{array}{lcl}
 \bfA \bfZ_k\!\!\!\!\!\!&=&\!\!\!\!\bfU_{k+1} \bfM_{k}\\
\bfA\t \bfU_{k+1}\!\!\!\!\!\!&=&\!\!\!\!\bfV_{k+1} \bfS_{k+1}
\end{array},
\:\mbox{where}\:
\begin{array}{lcll}
\bfZ_{k}\!\!\!\!&=&\!\!\!\![\bfW_1^{-1} \bfv_1, ..., \bfW_k^{-1} \bfv_k] &\!\!\!\!\in\bbR^{n\times k}\\
\bfM_k\!\!\!\!&=&\!\!\!\!\bfU_{k+1}\t\bfA\bfZ_k &\!\!\!\!\in\bbR^{(k+1)\times k}\\
\bfS_{k+1}\!\!\!\! &=&\!\!\!\!\bfV_{k+1}\t\bfA\bfU_{k+1} &\!\!\!\!\in\bbR^{(k+1)\times k}
\end{array}.
\end{equation}
Similar to \cref{eq:fArnoldi2}, the columns of $\bfZ_{k}$ span the approximation subspace for the solution $\bfx_k$, $\bfV_{k+1}$ and $\bfU_{k+1}$ are defined as in \cref{eq:partialFGKB1} and, by construction, $\bfM_k$ and $\bfS_{k+1}$ are upper Hessenberg and upper triangular, respectively. We emphasize that, in both \cref{eq:fArnoldi2} and \cref{eq:partialFGKB}, ${\rm ran}(\bfZ_k)$ is not necessarily a Krylov subspace (see \cite{simoncini2007recent}), and interpreting them as standard preconditioned Krylov subspaces is not straightforward (see \cite{notay2000flexible}). As in \cref{sec:hybrid}, we assume that both \cref{eq:fArnoldi2} and \cref{eq:partialFGKB} are breakdown-free, i.e., at iteration 
$k \leq \min \{ m,n\}$,
${\rm ran}(\bfZ_{k})$ has dimension $k$. The approximate solution $\bfx_k$ computed at the $k$th iteration is used to update the matrix $\bfW_{k+1}$ to be employed to expand the solution subspace at the $(k+1)$st iteration, as in \cref{eq:fArnoldi1} and \cref{eq:partialFGKB1}.

The available strategies \cite{chung2019, gazzola2014generalized, FKSIRW} that approximate the solution of problem (\ref{msl:eq:pbm_stdform}) through flexible Krylov methods all compute, at the $k$th iteration,
\begin{equation}\label{eq:flexiblelpproj}
\bfy_k=\argmin_{\bfy}\|\bfA\bfZ_k\bfy-\bfb\|_2^2 + \lambda\|\bfP_k\bfy\|_2^2,\quad \bfx_k=\bfZ_k\bfy_k\in{\rm ran}(\bfZ_k)\,,
\end{equation}
and essentially differ in the choice of the regularization matrix $\bfP_k\in\bbR^{p\times k}$. Specifically,
\begin{itemize}
\item[(i)] The method in \cite{gazzola2014generalized} and one of the methods in \cite{chung2019} (dubbed `hybrid-I') take $\bfP_k=\bfI_k$.
\item[(ii)] One of the methods in \cite{chung2019}  (dubbed `hybrid-R') takes $\bfP_k=\bfR_k\in\bbR^{k\times k}$, where $\bfR_k$ is the upper triangular factor in the `skinny' QR factorization of the basis vectors, i.e., $\bfZ_k=\bfQ_k\bfR_k$, which can be efficiently updated for consecutive $k$.
\item[(iii)] The methods in \cite{FKSIRW} take $\bfP_k=\bfR^{\bfW}_k\in\bbR^{k\times k}$, where $\bfR^{\bfW}_k$ is the upper triangular factor in the `skinny' QR factorization of the transformed basis vectors, i.e., $\bfW_k\bfZ_k=\bfQ^{\bfW}_k\bfR^{\bfW}_k$. Note that, unless $\bfW_k$ is a multiple of the identity matrix, the QR factorization cannot be efficiently updated for consecutive values of $k$, but it is still feasible to compute if $k\ll \min\{m,n\}$.
\end{itemize}
The three options above reveal that, when using hybrid projection methods based on flexible Krylov subspaces, even if problem \cref{msl:eq:approxpbm} is transformed into standard form \cref{msl:eq:pbm_stdform}, projecting and regularizing are not interchangeable anymore. Indeed, option (i) corresponds to a `first-project-then-regularize' scheme (i.e., it penalizes the 2-norm of the projected solution $\bfy_k$) while option (iii) corresponds to a `first-regularize-then-project' scheme;  option (ii) penalizes the 2-norm of the full-dimensional solution $\bfx_k$. When option (iii) is used, exploiting the majorization-minimization framework allows to prove that the sequence of approximate solutions $\{\bfx_k\}_{k\geq 1}$ to \cref{eq:flexiblelpproj} converges to a stationary point of problem (\ref{msl:eq:smoothpbm}). Note that, by exploiting the properties of the matrices appearing in \cref{eq:fArnoldi2} and \cref{eq:partialFGKB}, the fit-to-data term is equivalently expressed as $\|\bfH_k\bfy - \|\bfb\|_2\bfe_1\|_2^2$ in the flexible Arnoldi case, or $\|\bfM_k\bfy - \|\bfb\|_2\bfe_1\|_2^2$ in the flexible Golub-Kahan case.
Although the regularization parameter $\lambda$ in \cref{eq:flexiblelpproj} is displayed as independent of $k$, a heuristical adaptive choice is possible employing any of the strategies described in \cref{sec:param}. A suitable approximation to a solution of \cref{msl:eq:originalpbm} is efficiently computed if, for small $k$, the columns of $\bfZ_k$ can be used to capture the main features thereof.

\subsubsection{Sparsity under transform}\label{ssec:transf}

Many of the methods described above for problem \cref{msl:eq:originalpbm} can be generalized to approximate a solution to the more general variational regularization problem,
\begin{equation}\label{msl:eq:originaltransfpbm}
 \min_\bfx \|\bfA\bfx-\bfb \|_2^2+\lambda \calR(\bfx)\,,\quad\mbox{for many choices of the functional $\calR(\cdot)$}.
\end{equation}
All the methods based on inner-outer iterations described in \cref{ssec:inout} and the methods based on GKS described in \ref{ssec:GKSlplq} can handle regularization terms of the form $\calR(\bfx)=\|\bfL\bfx\|_p^p$: this can be done efficiently if matrix-vector products with $\bfL$ are cheap to compute. Note that the discrete (anisotropic) total variation regularization functional,
\begin{equation}\label{def:tv}
\TV(\bfx)=\|\bfD \bfx\|_p^p,\quad
\begin{array}{l}
\mbox{where $\bfD$ represents a discretization of the}\\
\mbox{(magnitude of the) gradient operator}
\end{array}
\end{equation}
can be expressed in this framework. When $0<p\leq 1$, such regularization terms can be interpreted as enforcing sparsity of the unknown $\bfD\bfx$. An IRN method for handling TV-like regularizers was proposed in \cite{wohlberg2007tv} that employs CGLS during the inner iterations, while \cite{arridge2014iterated} uses preconditioned LSQR applied to a problem that is effectively transformed into standard form. The authors of \cite{gazzola2020edges} propose an inner-outer iterative method that enhances edges through multiplicative updates of the weights, and which exploits the hybrid projection method based on joint bidiagonalization method \cref{eq:jointbd}. The methods based on flexible Krylov subspaces described in \cref{ssec:flexible} can be as well extended to handle problem \cref{msl:eq:originaltransfpbm}, but the strategies adopted to achieve this are regularizer-dependent: in other words, no universal approach is possible, even for the special case $\calR(\bfx)=\|\bfL\bfx\|_p^p$. For instance, \cite{chung2019} describes an approach that works when $\bfL$ is an orthogonal wavelet transform for methods based on both the flexible Arnoldi and the flexible Golub-Kahan algorithms. Also, \cite{GaSL18} describe an extension of FGMRES that handles total variation \cref{def:tv} by first transforming into standard form. Finally, hybrid Krylov projection methods based on either an inner-outer iterative scheme or a flexible scheme, have been used for nuclear norm regularization $\calR(\bfx)=\|\bfx\|_{\ast}$, which is meaningful when $\bfx$ is a vectorization of a low-rank 2D quantity (such as a low-rank 2D image); see \cite{gazzola2020krylov}.

\subsection{Beyond deterministic inversion: Hybrid projection methods in a Bayesian setting}
\label{sub:bayesian}

Recently, hybrid projection methods have found utility in the framework of statistical inverse problems, within a Bayesian approach. Our goal in this section is to draw important connections between deterministic and statistical inverse problems, so that the reader can understand the role that hybrid projection methods can play in the Bayesian setting. For more detailed descriptions, we refer the interested reader to excellent books and reviews on Bayesian inverse problems \cite{kaipio2006statistical,calvetti2007introduction,calvetti2018inverse, calvetti2018bayes} and computational uncertainty quantification \cite{bardsley2018computational}, and references therein.

Traditionally, regularization can be interpreted as a \emph{practical} approach to replace an ill-posed problem by a nearby well-posed one. For example, Tikhonov regularization was originally motivated as a tool to stabilize the solution process by computing one single deterministic solution of a modified problem \cite{Tik63a,Tik63b}; see also \cref{sec:background}.
However, the Bayesian framework also provides a systematic way for solving and analyzing inverse problems by modelling the unknown as a random variable. Although Bayesian inverse problems can be analyzed in a continuous setting \cite{stuart2010inverse}, here we focus on discrete inverse problems, where the underlying mathematical models have been discretized before applying Bayesian analysis.

In \cref{subsub:connection} we will draw connections between regularization and Bayesian inversion.  In particular, we will show that, under some specific assumptions on the prior, regularized (e.g., Tikhonov) problems arise in the Bayesian framework. One of the main benefits of the Bayesian approach is that different priors can be easily incorporated. Then, in \cref{subsub:UQ}, we describe some computational tools for UQ that can exploit hybrid projection methods. These tools bridge developments in the numerical linear algebra community with those in the UQ community for more efficient and practical inversion that can benefit a wide range of applications.

\subsubsection{Connections between regularization and Bayesian inversion}
\label{subsub:connection}
First some preliminary notation. Let $\bfmu \in \bbR^n$ and $\bfGamma \in \bbR^{n\times n}$ be symmetric positive definite (SPD); then, an $n$-variate normally distributed random variable $\bfy \sim \calN(\bfmu, \bfGamma)$ has a density function given by
\[\pi(\bfy) = \left(\frac{1}{2\pi |\bfGamma|}\right)^{n/2}\exp\left(-\frac{1}{2} \|\bfy - \bfmu\|_{\bfGamma^{-1}}^2 \right) \]
where $| \cdot |$ denotes the determinant of a matrix and $\| \bfz \|_\bfM^2 = \bfz\t \bfM \bfz$ for any $\bfz \in \bbR^n$ and $\bfM\in\bbR^{n\times n}$ SPD.

Consider the stochastic extension of \eqref{eq:linearmodel} where $\bfb, \bfx$, and $\bfe$ are random variables. For notational clarity, note that we have dropped the subscript on $\bfx$. Assume that $\bfx$ and $\bfe$ are independent and normally distributed,
\begin{equation}
  \label{eq:assumptions}
\bfx \sim \calN(\bfmu, \alpha^2\bfQ) \quad\text{and}\quad
\bfe \sim \calN(\bfzero, \sigma^2\bfR)
\end{equation}
where $\bfmu \in \bbR^n$ is the prior mean and $\bfQ \in \bbR^{n \times n}$ and $\bfR \in \bbR^{m \times m}$ are SPD matrices for the prior and noise respectively. Here we assume that $\alpha$ and $\sigma$ are known.  Although it is typically not necessary to include these parameters (i.e., they can be absorbed in the definition of $\bfQ$ and $\bfR$), we include them for clarity and to draw connections to hybrid projection methods in \cref{subsub:UQ}.

The probability density function for $\bfb$ given $\bfx$ is given by
 \[\pi_{\rm like}(\bfb \mid \bfx) = \left(\frac{1}{2\pi \sigma^2 |\bfR|}\right)^{m/2}\exp\left(-\frac{1}{2\sigma^2} \| \bfb - \bfA \bfx\|_{\bfR^{-1}}^2\right),\]
which is also called the \emph{likelihood function}.  Note that the estimate that maximizes the likelihood function (dubbed the maximum likelihood estimator or MLE) is the solution to an unregularized \emph{weighted} least-squares problem, i.e.,
\begin{align}
  \bfx_{\rm MLE} & =\arg\max \pi_{\rm like}(\bfb \mid \bfx)  =
 \argmin - \ln \pi_{\rm like}(\bfb \mid \bfx) \nonumber\\
 & =  \argmin \frac{1}{2\sigma^2} \| \bfb - \bfA \bfx\|_{\bfR^{-1}}^2  =  \argmin \frac{1}{2\sigma^2} \| \bfL_\bfR(\bfb - \bfA \bfx)\|_2^2 \nonumber
\end{align}
where $\bfR^{-1} = \bfL_\bfR^\top \bfL_\bfR$ is a symmetric (e.g., Cholesky) factorization.
Note that if $\bfR = \bfI_m$, then we get a standard least-squares problem.
The likelihood function is defined solely based on the assumptions from the noise model and extends naturally for nonlinear problems, e.g., if one were to consider a nonlinear forward model, \linebreak[4]$\bfb = F(\bfx)+\bfe$ where $F(\cdot):\bbR^n \to \bbR^m$ and $\bfe \sim \calN(\bfzero,\bfR)$, then the MLE estimate is the solution to the nonlinear least-squares problem, $\min_\bfx \frac{1}{2\sigma^2} \| \bfb - F(\bfx)\|_{\bfR^{-1}}^2$.  Additive Gaussian noise is the most common noise model used in the literature, but it is worth mentioning that in many applications, such as medical and atmospheric imaging, a Poisson noise model is appropriate. Specifically, in the linear case and assuming independence, $[\bfb]_i \sim Poisson([\bfA \bfx]_i)$; the likelihood model takes the form
\[
  \pi_{\rm like}(\bfb \mid \bfx) = \prod_{i=1}^m \frac{[\bfA \bfx]_i^{[\bfb]_i}}{[\bfb]_i !} \exp(-[\bfA \bfx]_i).
\]
Various noise models can be incorporated, leading to different likelihood functions.  For general data-fit terms, hybrid projection methods described in \cref{sub:flexible} can be used to approximate the MLE. For the remainder of this section, we consider additive Gaussian noise.

Next, we consider assumptions on prior information about $\bfx$, as described in \eqref{eq:assumptions}, where the prior density function is given by
\[\pi_{\rm prior}(\bfx) = \left(\frac{1}{2\pi \alpha^2 |\bfQ|}\right)^{n/2}\exp\left(-\frac{1}{2\alpha^2}\| \bfx - \bfmu \|_{\bfQ^{-1} }^2 \right).\]
 Using
Bayes' Theorem, we can derive the posterior density function
\begin{equation}
\begin{split}
\pi_{\rm post}(\bfx \mid \bfb) &= \frac{\pi_{\rm like}(\bfb \mid \bfx)\pi_{\rm prior}(\bfx)}{\pi (\bfb)}\\
&\propto \exp\left(-\frac{1}{2\sigma^2} \| \bfb - \bfA \bfx\|_{\bfR^{-1}}^2 -\frac{1}{2\alpha^2}\| \bfx - \bfmu \|_{\bfQ^{-1} }^2 \right) \label{eq:post_RQ}
\end{split}
\end{equation}
where $\propto$ stands for `proportional to'.
\emph{A crucial point to note is that for Bayesian inverse problems, the posterior density function is the solution to the inverse problem.}  The posterior distribution provides the full information about the distribution of parameters in $\bfx$, given the data $\bfb$. However, for practical interpretation and data analysis, it is necessary to describe various characteristics of the posterior distribution \cite{tenorio17}; thus leading to the field of \emph{uncertainty quantification}.
For example, a common goal is to compute the maximum a posteriori (MAP) estimator,
i.e., the point estimate that corresponds to the maximum of \eqref{eq:post_RQ}.  We denote this as
\begin{align}
  \bfx_{\rm MAP} &  = \arg \max_\bfx \pi_{\rm post}(\bfx \mid \bfb) \label{eq:MAP}\\
  & = \argmin_\bfx \frac{1}{2\sigma^2} \| \bfb - \bfA \bfx\|_{\bfR^{-1}}^2 + \frac{1}{2\alpha^2}\| \bfx - \bfmu \|_{\bfQ^{-1} }^2 \nonumber\\
  & = \left(\tfrac{1}{\sigma^2}\bfA^\top\bfR^{-1}\bfA + \tfrac{1}{\alpha^2}\bfQ^{-1}\right)^{-1}\left(\tfrac{1}{\sigma^2} \bfA^\top \bfR^{-1} \bfb + \tfrac{1}{\alpha^2} \bfQ^{-1} \bfmu \right), \label{eq:closedform}
\end{align}
where the closed form solution in \eqref{eq:closedform} can be obtained by setting the derivative to zero.

Furthermore, it can be shown that, under the above assumptions (linear problem with Gaussian noise and prior), the posterior is indeed Gaussian, where the posterior covariance matrix and the posterior mean (i.e., the MAP estimate) are given as
\begin{equation}\label{eq:postcovmean}
\bfGamma \equiv (\lambda\bfQ^{-1} + \bfA\t\bfR^{-1}\bfA)^{-1} \quad\text{and}\quad \bfx_{\rm MAP} =  \bfGamma(\bfA\t\bfR^{-1} \bfb + \lambda\bfQ^{-1} \bfmu),
\end{equation}
respectively, with $\lambda = \frac{\sigma^2}{\alpha^2}$ \cite{kaipio2006statistical}.

Next we discuss some examples of priors and show how the posterior distribution and, in particular, the MAP estimate is defined. \emph{Gaussian priors} are the most common priors used in the literature.  We split the discussion into three scenarios. First, for cases where the symmetric factorization of the precision matrix (i.e., the inverse of the covariance matrix) is computationally feasible, i.e., $\bfQ^{-1} = \bfL_\bfQ^\top \bfL_\bfQ$, the MAP estimate is the solution to the optimization problem,
\[
    \min_{\bfx} \frac{1}{2\sigma^2} \| \bfb - \bfA \bfx \|_{\bfR^{-1}}^2 + \frac{1}{2\alpha^2}\| \bfL_\bfQ (\bfx - \bfmu) \|_2^2.
\]
With a simple change of variables, $\bfy = \bfx-\bfmu$ and $\bfd = \bfb - \bfA \bfmu$, the problem reduces to solving
\[
  \min_\bfy \frac{1}{2\sigma^2} \| \bfd - \bfA \bfy \|_{\bfR^{-1}}^2 + \frac{1}{2\alpha^2}\| \bfL_\bfQ \bfy \|_2^2.
\]
Let $\bfR = \bfI_m$: then, this is general-form Tikhonov problem \eqref{eq:genTikhonov}, and hybrid projection methods to solve this problem were discussed in \cref{sub:generalform}.
For this scenario, a prevalent approach in the literature models $\bfL_\bfQ$ as a sparse discretization of a differential operator (for example, the Laplacian or the biharmonic operator). This choice corresponds to a covariance matrix that represents a Gauss-Markov random field.  Since for these choices, the precision matrix is sparse, working with $\bfL_\bfQ$ directly has obvious computational advantages.

Second, for the case where application of $\bfQ^{-1}$ is feasible (e.g, it is a sparse matrix) but computing the symmetric factorization is too expensive, a factorization-free preconditioned LSQR approach called MLSQR was proposed in~\cite{arridge2014iterated}. The approach is analytically equivalent to LSQR applied to a preconditioned least-squares problem with the preconditioner $\bfL_\bfQ^{-1}$, but avoids the factorization by using weighted inner products to implicitly apply the preconditioner. The authors considered applications in the context of nonlinear regularizers.

Third, there are many prior models where the precision matrix $\bfQ^{-1}$ and its factorization are not available, but matrix-vector multiplications with $\bfQ$ can be done efficiently.
For example, for Gaussian random fields, entries of the covariance matrix are computed directly as $\bfQ_{ij} = \kappa(\bfx_i,\bfx_j)$, where $\{\bfx_i\}_{i=1}^n$ are the spatial points in the domain and $\kappa$ is a covariance kernel function (e.g., $\gamma$-exponential, or Mat\'ern class).  Although there are many modeling advantages, the main challenge is that the resulting prior covariance matrices are often very large and dense; explicitly forming and factorizing these matrices is prohibitively expensive. For such prior models, efficient matrix-free techniques (e.g., FFT embedding and $\calH$-matrix approaches) can be used, for example, to compute matrix-vector products with the prior covariance matrix $\bfQ$. Generalized hybrid projection methods were proposed in \cite{ChSa17} to compute Tikhonov regularized solutions (i.e., MAP estimate \eqref{eq:MAP}) effectively.  The approach relies on a change of variables and uses a generalized GKB algorithm for projection \cite{arioli2013generalized}; note that this is a fundamentally different algorithm than the one presented in \cref{ssec:GKSlplq} for $\ell_2-\ell_p$ regularization that is based on GKB (for initialization) and some generalization of Krylov projection methods.

 Non-Gaussian priors are also common in the literature, but are more challenging to handle.  One class of priors contains the so-called \emph{sparsity-promoting priors} or $\ell_p$-priors, which have the density function,
\begin{equation}\label{eq:sparseprior}
    \pi_{\rm prior}(\bfx) = \exp(-\alpha \| \bfx \|_p^p), \quad \mbox{for} \quad 0 < p \leq 1.
\end{equation}
  One of the main challenges in considering \cref{eq:sparseprior} is that the posterior distribution is no longer Gaussian.  However, one could still compute the MAP estimate, which is the solution to optimization problem,
    \[
      \min_{\bfx} \frac{1}{2\sigma^2} \| \bfb - \bfA \bfx \|_{\bfR^{-1}}^2 + \alpha \|  \bfx \|_p^p.
    \]
Note that $p=1$ corresponds to a prior where the components of $\bfx$ are independent and follow the univariate Laplace$(0, \alpha^{-1})$ distribution~\cite[Chapter 4.3]{bardsley2018computational}; the Laplacian prior enforces sparsity in $\bfx$.
On a related note, \emph{Besov priors} can preserve sharp discontinuous interfaces in Bayesian inversion, and when combined with a wavelet-based approach, they can promote sparsity in the MAP estimates \cite{bui2015scalable,vanska2009statistical,dashti2012besov} but, to the best of our knowledge, hybrid projection methods have not been used in this setting. For many non-Gaussian priors, the methods described in \cref{sub:flexible} may be used to approximate the MAP estimate.

\subsubsection{Exploiting hybrid projection methods for UQ}
\label{subsub:UQ}
Recall that the solution to a Bayesian inverse problem is the posterior distribution, and the goal is to describe and explore it. For linear inverse problems with Gaussian priors, if the regularization parameters $\sigma$ and $\alpha$ are set in advance, the posterior is Gaussian with covariance matrix \eqref{eq:postcovmean}, and various tools can be used for UQ.
However, if the regularization parameters are not known in advance, one would consider a fully Bayesian framework, where all unknown parameters (including the regularization parameters) are treated as random variables, and use hierarchical models \cite{bardsley2018computational}, also defining hyperpriors for these hyperparameters. For example, for Gaussian priors, since $\sigma$ and $\alpha$ are unknown, we can assume that they are random variables with hyperpriors $\pi(\sigma)$ and $\pi(\alpha)$, respectively.  Then, using again Bayes' law, we have the full posterior density,
\begin{equation}
\pi_{\rm post}(\bfx, \sigma, \alpha | \bfb) \propto \pi_{\rm like}(\bfb|\bfx)\pi_{\rm prior}(\bfx) \pi(\sigma) \pi(\alpha).
 \label{eq:post_hyper}
\end{equation}
For linear inverse problems with Gaussian priors, common choices are Gamma hyperpriors, i.e.,
\begin{align}
   \pi(\sigma) & \propto \sigma^{\beta_\sigma-1} \exp(-\gamma_\sigma \sigma),\nonumber \\
   \pi(\alpha) & \propto \alpha^{\beta_\alpha-1} \exp(-\gamma_\alpha \alpha),\nonumber
 \end{align}
where $\beta_\sigma, \gamma_\sigma, \beta_\alpha,$ and $\gamma_\alpha$ are parameters.
In addition to requiring the user to pre-define these additional parameters for the hyperpriors, a potential computational disadvantage of this approach is that the posterior \cref{eq:post_hyper} is no longer Gaussian, so more sophisticated sampling techniques (e.g., hierarchical Gibbs sampling) should be used. However, recall that the main connection between hybrid projection methods and Bayesian inverse problems is that hybrid projection methods can be used to compute MAP estimates efficiently while selecting the regularization parameter automatically (often based on some statistical tools, see \cref{sec:param}). Therefore, a simple \emph{practical} alternative to adopting a full Bayesian paradigm is to use a hybrid projection method to estimate the regularization parameter and, once this is done, perform standard UQ with fixed hyperparameters \cite{saibaba2020efficient}.

Handling unknown hyperparameters is not the only task where hybrid projection methods can be beneficial for efficient UQ. In the remaining part of this section, we consider other more \emph{practical} scenarios that combine hybrid projection methods with UQ. We still focus on the Gaussian, linear case where the posterior is given in \cref{eq:post_RQ}. We begin by remarking that there has been plenty of work at the intersection of Krylov methods and UQ.  For example, the idea of using low-rank perturbative approximations for the posterior covariance matrix previously appeared in~\cite{flath2011fast,bui2012extreme, bui2013computational, spantini2015optimal}.  Along these lines, there has also been some work on using randomized approaches to efficiently compute a low-rank approximation for use in UQ settings \cite{saibaba2015fastc}; however, such methods are not ideal for inverse problems where the decay of the singular values is not sufficiently rapid (e.g., in tomography applications).
Previous work on Lanczos methods for sampling from Gaussian distributions {can be found in, e.g.,}~\cite{parker2012sampling,schneider2003krylov,simpson2008krylov,chow2014preconditioned}, but these algorithms are meant for sampling from generic Gaussian distributions and do not exploit the structure of the posterior covariance matrix \cref{eq:post_RQ}. We point the interested reader to \cite{bardsley2018computational,calvetti2018inverse} for more examples where Krylov methods are used for UQ. Here we consider the role played by hybrid projection methods, and list some further approaches for combining hybrid projection methods and UQ. Specifically:

\begin{itemize}
    \item Hybrid projection methods can be used to generate an approximation to the posterior covariance matrix. Uncertainty measures can then be computed by storing bases for the Krylov subspaces. This approach was considered for dynamic inverse problems using generalized Golub-Kahan approaches in \cite{chung2018efficient}.
    \item Low-rank approximations for proposal sampling are commonly used for efficient UQ. In \cite{brown2018low}, the authors consider a fully Bayesian approach with assigned hyperpriors and use MCMC methods with Metropolis-Hastings independence sampling with a proposal distribution based on a low-rank approximation of the prior-preconditioned Hessian. In \cite{saibaba2019efficient}, these ideas are combined with marginalization.
\end{itemize}
Theoretical results on approximating the posterior distribution and posterior covariance matrix can be found in \cite{saibaba2020efficient,simon2000low}.

\subsection{Beyond linear forward models: Hybrid projection methods for nonlinear inverse problems}
\label{sub:nonlinear}

In this section, we consider the role that hybrid projection methods have played in the context of solving nonlinear inverse problems. We remark that the literature on nonlinear inverse problems is vast, and it is not our intention to survey this topic.  We point the interested reader to books on this topic, see e.g., \cite{engl2005nonlinear,mueller2012linear,seo2012nonlinear}, and we focus on frameworks and methodologies for solving nonlinear inverse problems that have successfully utilized or incorporated hybrid projection methods.

We are concerned with nonlinear inverse problems, where the forward model depends nonlinearly on the desired parameters.  Consider \emph{discrete} nonlinear inverse problem,
\begin{equation}
  \label{eq:nonlinearmodel}
  \bfb = F(\bfx_\true) + \bfe
\end{equation}
where $F:\bbR^n \to \bfR^m$ is a transformation representing the forward model. Notice that the linear model \cref{eq:linearmodel} is a special case of \eqref{eq:nonlinearmodel} with $F(\bfx_\true) = \bfA \bfx_\true$.

Regularization methods for nonlinear inverse problems typically follow a variational approach, where the goal is to solve an optimization problem of the form,
\begin{equation}
  \label{eq:nonlinopt}
  \min_\bfx \calJ(\bfb,F(\bfx)) + \lambda \calR(\bfx)
\end{equation}
where similar to \eqref{eq:directreg}, $\calJ$ is some loss function, $\calR$ is a regularization operator, and $\lambda \geq 0$ is a regularization parameter.

There are two main computational difficulties in solving problems like \eqref{eq:nonlinopt}, especially for large-scale problems. First, the regularization parameter is unknown, and estimating it may require a significant computational effort to solve the same nonlinear optimization problem multiple times for various parameter choices, which can be a very expensive and time-consuming task; see \cite{engl1996regularization}.
Second, since the optimization functional to be minimized is nonlinear, gradient-based iterative techniques (even Newton-type approaches) are needed, but computing derivatives can be costly \cite{haber2000optimization}. We remark that similar challenges arise for linear inverse problems with nonlinear regularization terms (see, e.g., \cref{msl:eq:originalpbm2}), where sophisticated nonlinear optimization schemes are required; we point the reader to the discussion in 
\cref{sub:flexible} for some efficient approaches to deal with these kinds of nonlinearities.
Alternatively, a two-stage method could also be used that splits the inversion process. First the misfit $\calJ$ is reduced to some target misfit value.  In the second stage, the target misfit is kept constant and the regularization term $\calR$ is reduced. Although very popular in practice, this approach is not guaranteed to converge (in fact it diverges in some cases) and appropriate safety steps and ad hoc parameters must be used; see \cite{parker1994geophysical,constable1987occam}.

Specifically for nonlinear least-squares problems, other solvers can be used, such as the one described in Haber and Oldenburg \cite{haber2000gcv}, which  combines a damped Gauss-Newton method for local regularization with a GCV method for adaptively selecting the global regularization parameter. To the best of our knowledge, this is the first paper that exploits hybrid projection methods for selecting the regularization parameter for solving large-scale nonlinear inverse problems.  Thus, we further describe the approach. Consider problem \eqref{eq:nonlinopt} where $\calJ(\bfb,F(\bfx)) = \norm[2]{F(\bfx)-\bfb}^2$ and $\calR(\bfx) = \norm[2]{\bfx}^2$, such that we get the Tikhonov-regularized nonlinear problem,
\begin{equation}
  \label{eq:nonlinLS}
  \min_\bfx \norm[2]{F(\bfx)-\bfb}^2 + \lambda \norm[2]{\bfx}^2.
\end{equation}
Let's begin by assuming that $\lambda$ is fixed.  If we differentiate the objective function in \eqref{eq:nonlinLS} with respect to $\bfx$ and set the result to zero, we get
\begin{equation}
  \label{eq:gradzero}
  \bfg(\bfx) = 2(\bfJ(\bfx)\t (F(\bfx)-\bfb) + \lambda \bfx) = \bfzero,
\end{equation}
where $\bfg(\bfx)$ is the gradient and $\bfJ(\bfx) = \frac{\partial F}{\partial \bfx} \in \bbR^{m \times n}$ is the Jacobian matrix that contains the partial derivatives or sensitivities of the forward model.
If we find a solution to \eqref{eq:gradzero}, we have the desired solution to \cref{eq:nonlinLS} for fixed $\lambda$.
Rather than using Newton's method, which requires calculating the second derivative of $F$ with respect to $\bfx$ (i.e., differentiation of $\bfg(\bfx)$), a damped Gauss-Newton approach can be used.  That is,
let $F(\bfx + \delta \bfx) = F(\bfx) + \bfJ(\bfx) \delta \bfx + R(\bfx,\delta\bfx), $ be a linearization of $F$.
Since $R(\bfx,\delta\bfx)$ is assumed small, we aim to solve the Gauss-Newton equations,
\begin{equation}
  \label{eq:normaleq}
   (\bfJ(\bfx)\t \bfJ(\bfx) +\lambda \bfI_n)\delta \bfx  =  - \bfJ(\bfx)\t (F(\bfx) -\bfb) - \lambda \bfx.
\end{equation}
or equivalently the least-squares problem,
\begin{equation}
  \label{eq:LSeq}
   \min_{\delta \bfx} \left\| \begin{bmatrix} \bfJ(\bfx) \\ \sqrt{\lambda} \bfI_n \end{bmatrix} \delta \bfx - \begin{bmatrix} \bfb - F(\bfx) \\ \sqrt{\lambda} \bfx\end{bmatrix} \right\|.
\end{equation}
The Gauss-Newton method is an iterative approach where given an initial guess $\bfx_0,$ the update at the $k$th iteration is computed by solving \eqref{eq:LSeq} for the perturbation $\delta \bfx$ and updating the solution estimate as $\bfx_{k+1} = \bfx_k + \alpha_k \delta \bfx$, where $\alpha_k$ is a line-search parameter.
Notice that, for a given $\lambda,$ Krylov iterative methods could be used to efficiently estimate the solution to \eqref{eq:LSeq}. However, the authors of \cite{haber2000gcv} propose a reformulation that enables the direct use of hybrid projection methods described in Section \ref{sec:hybrid}.  The astute observation is that at each iteration of the Gauss-Newton method, we have the option to solve directly for the step $\delta \bfx$ or to solve for the updated solution $\bfx_{k+1}$.  To see the latter option, substitute $\bfx_{k+1} = \bfx_k + \delta \bfx$ into \eqref{eq:normaleq} to get least-squares problem,
\begin{equation}
  \label{eq:reformulateLS}
   \min_\bfx \norm[2]{\begin{bmatrix} \bfJ(\bfx_k) \\ \sqrt{\lambda} \bfI_n \end{bmatrix} \bfx - \begin{bmatrix} \bfb - F(\bfx_k) - \bfJ(\bfx_k)\bfx_k \\ \bfzero \end{bmatrix}}^2.
\end{equation}
At this point, \cref{eq:reformulateLS} is just a standard-form, linear Tikhonov problem, and hybrid projection methods can be used in each nonlinear step to
simultaneously estimate the regularization parameter $\lambda$ and compute an approximate solution. Although the hybrid projection method solves the linearized problem with automatic regularization parameter selection, a major concern is that the objective function changes from iteration to iteration because the regularization parameter is changing.  Furthermore, we are not guaranteed to reduce the nonlinear function that we seek to minimize. A remedy is to consider the solution to \eqref{eq:reformulateLS} as a proposal $\bfx_{k+1}^p$ and to take the step direction to be a weighted version of the perturbation $\delta \bfx = \bfx_{k+1}^p - \bfx_k$. The general damped Gauss-Newton approach with a hybrid Krylov solver is summarized in Algorithm \ref{alg:GaussNewton}.

\begin{algorithm}
\caption{Gauss-Newton with hybrid projection for nonlinear problem \eqref{eq:nonlinLS}}
\label{alg:GaussNewton}
\begin{algorithmic}[1]
	\REQUIRE{$F$, $\bfb$, $\bfx_0$, $k=0$}
\WHILE {stopping criterion not satisfied}
\STATE Calculate $\bfJ(\bfx_k)$ and $\bfr(\bfx_k)$
\STATE Compute $\bfx_{k+1}^p$ by solving \eqref{eq:reformulateLS} using a hybrid projection method
\STATE Calculate the perturbation $\delta \bfx =\bfx_{k+1}^p - \bfx_k$
 \STATE Update $\bfx_{k+1} = \bfx_k + \alpha \delta \bfx$ where $\alpha$ ensures reduction of the nonlinear function
 \STATE $k = k+1$
\ENDWHILE
\ENSURE{$\bfx_k$}
\end{algorithmic}
\end{algorithm}
Another adaptive approach to select the regularization parameter when solving nonlinear inverse problems is based on continuation or cooling approaches \cite{haber2000optimization}.  By reformulating the Tikhonov problem as a constrained optimization problem where  the objective function is given in terms of the regularizer, and the constraints are given in terms of the data misfit, and exploiting the connection between the regularization parameter and the inverse of the Lagrange multiplier for the constraint, the regularization parameter can be selected to be a large value initially, and then gradually reduced to an appropriate value.

Finally, in the context of solving large-scale  nonlinear inverse problems, there is another class of nonlinear inverse problems that have greatly benefited from recent developments in hybrid projection methods.  These are \emph{separable nonlinear inverse problems}, where the unknown parameters can be separated into two distinct components such that the forward model is linear in one set of parameters and nonlinear in another set of parameters. This situation corresponds to problem
\eqref{eq:nonlinearmodel} with
$$F(\bfx_\true) = \bfA\left(\bfx_\true^{(nl)}\right)\ \bfx_\true^{(l)},$$
where $\bfx_\true = \begin{bmatrix}\bfx_\true^{(nl)} &\bfx_\true^{(l)} \end{bmatrix}\t$, with $\bfx_\true^{(nl)}\in \bbR^\ell$ and $\bfx_\true^{(l)} \in \bbR^{n-\ell}$ and $\bfA$ is a linear operator defined by a set of nonlinear parameters $\bfx_\true^{(nl)}$.  Therefore $F$ is linear in $\bfx_\true^{(l)}$ and nonlinear in $\bfx_\true^{(nl)}.$  Furthermore, it is often the case that $\ell \ll n-\ell$.  Such problems arise in imaging applications such as blind deconvolution, super-resolution imaging, and motion correction \cite{ChHaNa06,ChNa10}, but also in other inverse problems; see \cite{GolubVP}. Although one could treat separable nonlinear inverse problems using standard optimization techniques, approaches that exploit their separable structure have been successfully used to improve convergence. For example, a decoupled or alternating optimization approach is a simple approach where, for fixed $\bfx^{(nl)}$, a hybrid projection method could be used to solve the resulting linear inverse problem, and, for fixed $\bfx^{(l)}$, a nonlinear optimization scheme can be used to solve the smaller dimensional optimization problem.  However, a potential disadvantage of these alternating approaches is slow convergence and sensitivity to initializations.  An alternative  approach is to use the variable projection method \cite{golub1973differentiation,GoPe03,GolubVP,kaufman1975variable,o2013variable}, where the linear variables are implicitly eliminated and a reduced optimization problem is solved. Consider the full Tikhonov problem,
\begin{equation}
\label{eq:NLLS}
\min_{\bfx^{(nl)}, \bfx^{(l)}}  \|\bfA\left(\bfx^{(nl)}\right)\,\bfx^{(l)} - \bfb\|_2^2 + \lambda\|\bfx^{(l)}\|_2^2 = \left\| \left[\begin{array}{c} \bfA\left(\bfx^{(nl)}\right) \\ \sqrt{\lambda} \bfI_n \end{array} \right] \bfx^{(l)} -
    \left[ \begin{array}{c} \bfb \\ \bfzero \end{array} \right]
\right\|_2^2.
\end{equation}
A variable projection approach applied to \eqref{eq:NLLS} would seek the solution to the reduced optimization problem,
\begin{equation}
\label{eq:reduced}
  \min_{\bfx^{(nl)}}  \|\bfA(\bfx^{(nl)})\,\bfx^{(l)}(\bfx^{(nl)}) - \bfb\|_2^2,
\end{equation}
where
\begin{equation}
  \label{eq:linearsubprob}
  \bfx^{(l)}(\bfx^{(nl)}) = \argmin_{\bfx^{(l)}} \left\| \left[\begin{array}{c} \bfA(\bfx^{(nl)}) \\ \sqrt{\lambda} \bfI_n \end{array} \right] \bfx^{(l)} -
    \left[ \begin{array}{c} \bfb \\ \bfzero \end{array} \right]
\right\|_2^2.
\end{equation}
Each iteration of a nonlinear optimization method to solve \eqref{eq:reduced} would require solving linear subproblems \eqref{eq:linearsubprob}.  Since this is a linear inverse problem, hybrid projection methods can be used (especially in the large-scale setting), and the regularization parameter can be estimated automatically, giving rise to an inner-outer iterative scheme \cite{ChHaNa06,ChNa10}. A more recent hybrid regularization method, which avoids inner-outer iterations by exploiting inexact Krylov solvers, is presented in \cite{gazzola2021}. For solving inverse problems with coupled variables, an approach called Linearize and Project (LAP) was described in \cite{nagy2018lap} as an alternative to variable projection methods. LAP can support different regularization strategies as well as equality and inequality constraints, and thus is more broadly applicable than variable projection. Hybrid regularization methods are used in LAP to simultaneously compute the search direction and automatically select an appropriate regularization parameter at each iteration.

Although hybrid projection methods have found practical uses in various nonlinear scenarios, there are still open questions and yet-to-be-explored uses of hybrid projection methods for solving nonlinear inverse problems. Recent advancements in flexible hybrid methods and iteratively reweighted approaches for solving linear inverse problems with nonlinear regularization terms have opened the door to new approaches for handling nonlinear optimization problems.  Furthermore, the training of many state-of-the-art deep neural networks can be interpreted as a separable nonlinear inverse problems \cite{newman2020train}, and there is great potential for the use of hybrid projection methods for training deep neural networks \cite{kan2020avoiding}.

\section{Software}
\label{sec:software}

The most recent and increasingly popular 
software package for hybrid projection methods is called \textsc{IR Tools} \cite{IRtools}, and it runs in MATLAB.  This toolbox provides implementations of a range of iterative solvers, including many of the hybrid projection methods discussed in this paper and extending the GKB-based hybrid solver originally implemented in \cite{NaPaPe04}. Various examples (including image deblurring and tomography reconstruction) are also provided as test problems. For specific details, we point the reader to the \textsc{IR Tools} software package and documentation. However, we point out a few important design objectives and considerations for the implementation of hybrid projection methods for inverse problems:
\begin{itemize}
\item The implementation should work for the common scenarios where the coefficient matrix is only available as a (sparse) matrix, a function handle, or an object. Matrix-vector multiplications with the coefficient matrix represent a core component (and oftentimes the computational bottleneck) of iterative projection methods.  High performance or distributed computing tools can be used to accelerate the task of matrix-vector multiplication.  For example, hybrid projection methods were used with message passing tools for cryo-EM reconstruction \cite{ChStYa09} and were included in a high-performance java library in \cite{wendykier2010parallel}.
\item The software should enable automatic regularization parameter selection and automatic stopping criteria, but at the same time provide functionality for users to change algorithm parameters in an options structure.
\item Compared to standard iterative methods, one of the main disadvantages of hybrid methods is the need to store the basis vectors for solution computation. Thus, it is often assumed that solutions can be captured in relatively few iterations or that appropriate preconditioning can be used. Thus, it is desirable that implementations of hybrid projection methods can incorporate preconditioning techniques or exploit recycling techniques.  Another potential issue is the computational costs associated with reorthogonalization \cite{barlow2013reorthogonalization,BjGrDo94}.
\item The return structure should provide relevant information (e.g., for subsequent UQ or prediction).
\end{itemize}

The illustrations for this paper were made using the IR Toools software, and MATLAB scripts used to produce the illustrations and experiments reported herein are available at the website: \small{\url{https://github.com/juliannechung/surveyhybridprojection}}.
\normalsize
Although \textsc{IR Tools} contains implementations of many state-of-the-art hybrid projection methods, it is by no means exhaustive and the authors encourage further extensions.  For some recent work on generalized hybrid projection methods, subsequent UQ, and hybrid methods with recycling, see \small{\url{https://github.com/juliannechung}}\normalsize.  For getting familiar with inverse problems (e.g., some toy problems), see \textsc{Regularization Tools} \cite{Han94}.  For image deblurring problems, see \textsc{Restore Tools} \cite{NaPaPe04}. For tomographic reconstruction problems, see \textsc{AIR Tools II} \cite{hansen2018air}.  For computational UQ, see accompanying MATLAB codes for \cite{bardsley2018computational}.

\section{Concluding remarks and future outlook}
\label{sec:outlook}
This survey provides an overview of hybrid projection methods for solving large-scale inverse problems.  Our aim is to provide a gentle introduction to the area for budding scientists and interested researchers working in related fields. We explain the general principles underlying hybrid projection methods for standard Tikhonov regularization (Algorithm \ref{alg:hybridprojection}), and carefully describe the most common algorithmic choices (including subspace projection approaches and parameter choice techniques), and finish with a summary of extensions that go \emph{beyond} standard approaches and uses.

The field of hybrid projection methods has gained significant momentum in the past few years, and there are many new areas of exploration.  Large-scale inverse problems continue to motivate the development of more efficient and accurate solvers, which often leverage tools from different areas of numerical linear algebra, optimization, and statistics.
We envisage that a lot of scientific research will be devoted to using hybrid projection methods to incorporate more sophisticated priors and for performing subsequent UQ, e.g., for further prediction and forecasting, in forthcoming years.
There is still need for the development and analysis of new parameter selection rules, and strategies that can incorporate supervised learning techniques in optimal experimental design and deep learning frameworks seem promising \cite{afkham2021learning}. Efficient implementations of hybrid methods will be needed for advanced distributed computing architectures and sampling techniques for massive or streaming data problems.  With the plethora of examples of inverse problems in areas ranging from biomedical imaging to atmospheric monitoring for threat detection, it is inevitable that many scientific and engineering applications will benefit from recent and forthcoming advancements in hybrid projection methods for inverse problems.

\section{Acknowledgements} The authors would like to thank Carola-Bibiane Sch\"{o}n\-lieb for encouraging us to write this review article. This work was partially supported by the National Science Foundation (NSF) under grant DMS-1654175 and 1723005 (J. Chung) and by the Engineering and Physical Sciences Research Council (EPSRC) under grant EP/T001593/1 (S. Gazzola).

\bibliographystyle{abbrv}
\bibliography{7references}
\end{document}